\documentclass[a4paper,onecolumn,11pt]{article}  

\usepackage{graphicx}

\usepackage{amsmath,amsfonts,setspace}
\usepackage{amssymb, epsfig,subfigure}
\usepackage{color}

\usepackage{soul} 
\usepackage[applemac]{inputenc}

\usepackage[english]{babel}
\usepackage{fancyhdr}

\newcommand{\x}{\mathbf{x}}

\newcommand{\vN}{\mathbf{N}}
\newcommand{\n}{\mathbf{n}}
\newcommand{\vV}{\mathbf{V}}
\newcommand{\vR}{\mathbf{R}}

\newcommand{\vb}{\mathbf{b}}
\newcommand{\vd}{\mathbf{d}}
\newcommand{\vM}{\mathbf{M}}

\newcommand{\R}{\mathbb{R}}

\newcommand{\grass}[1]{\boldsymbol{#1}} 
\newcommand{\vomega}{\grass{\omega}}

\newtheorem{theorem}{Theorem}[section]

\newcommand{\vA}{\mathbf{A}}

\newcommand{\ve}{\mathbf{e}}

\begin{document}

\title{Analysis and approximation of some Shape--from--Shading models for non-Lambertian surfaces}

\author{Silvia Tozza\footnote{Dipartimento di Matematica, Sapienza - Universit\'a di Roma,
             Piazzale Aldo Moro, 5 - 00185 Roma,
              \it{tozza@mat.uniroma1.it}         }      
               \and
        Maurizio Falcone\footnote{Dipartimento di Matematica, Sapienza - Universit\'a di Roma, 
              Piazzale Aldo Moro, 5 - 00185 Roma,
              \it{falcone@mat.uniroma1.it}     }
}

\date{January 27, 2016}
\maketitle

\begin{abstract}
The reconstruction of a 3D object or a scene is a classical inverse problem in Computer Vision. In the case of a single image this is called the Shape--from--Shading (SfS) problem and it is known to be ill-posed even in a simplified version like the vertical light source case. A huge number of works deals with the orthographic SfS problem based on the Lambertian reflectance model, the most common and simplest model which leads to an eikonal type equation when the light source is on the vertical axis. 
In this paper we want to study non-Lambertian models since  they are more realistic and suitable whenever one has to deal with different kind of surfaces, rough or specular. We will present a unified mathematical formulation of some popular orthographic non-Lambertian models, considering vertical and oblique light directions as well as different viewer positions.  These models lead to more complex stationary nonlinear partial differential equations of Ha\-milton-Jacobi type which can be regarded as the generalization of the classical eikonal equation corresponding to the Lambertian case.  However, all the equations corresponding to the models considered here (Oren-Nayar and Phong) have a similar structure so we can  look for weak solutions to this class in the viscosity solution framework. Via this unified approach, we are able to develop a semi-Lagrangian approximation scheme for the Oren-Nayar and the Phong model and to prove a general convergence result. Numerical simulations on synthetic and real images will illustrate the effectiveness of this approach and the main features of the scheme, 
also comparing the results with previous results in the literature.\\

\end{abstract}

\section{Introduction} \label{intro}
The 3D reconstruction of an object starting from one or more images is a very interesting inverse problem with many applications.
In fact, this problem appears in various fields which range from the digitization of curved documents \cite{CCDG07} to the reconstruction of archaeological artifacts \cite{FDSB04}. More recently, other applications have been considered  in astronomy to obtain a characterization of properties of planets and other astronomical entities \cite{SKWMM11,Grumpe2014,Lohse2006} and in security where the same problem has been applied to the facial recognition of individuals.\\
In real applications, several light sources can appear in the environment and the object surfaces represented in the scene can have different
reflection properties because they are made by different materials, so it would be hard to imagine a scene which can satisfy the classical assumptions of the 3D reconstruction models. In particular, the typical Lambertian assumption often used in the literature has to be weakened.
Moreover, despite the fact that the formulation of the Shape-from-Shading problem is rather simple for a single light source and under Lambertian assumptions, its solution is hard and requires rather technical mathematical tools as the use of weak solutions to nonlinear partial differential equations (PDEs). From the numerical point of view the accurate approximation of non regular solutions to these nonlinear PDEs is still a challenging problem. 
In  this paper we want to make a step forward in the direction of a mathematical formulation of non-Lambertian models in the case 
of orthographic projection with a single light source located far from the surface. In this simplified framework, we present a unified approach to two popular models for non-Lambertian surfaces proposed by  Oren-Nayar \cite{ON_SIGGRAPH94,ON_ECCV94,ON_CVPR93,ON_IJCV95}  and by Phong \cite{Phong75}. We will consider  light sources placed in oblique directions with respect to the surface and we will use that unified formulation to develop a general numerical approximation scheme which is able to solve the corresponding nonlinear partial differential equations arising in the mathematical description of these models. \\
To better understand the contribution of this paper, let us start from the classical SfS problem where the goal is to reconstruct the surface from a single image.
In mathematical terms,  given  the shading informations contained in a single two-dimensional gray level digital image $ I(\x)$, where $\x := (x,y)$, we look for a surface $z = u(\x)$ that corresponds to its shape (hence the name \emph{Shape from Shading}). This problem is described in general by the image irradiance equation introduced by Bruss \cite{Bruss81}
\begin{equation}\label{general_irradiance_eq}
I(\x) = R(\vN(\x)),
\end{equation}
where the normalized brightness of the given grey-value image $I(\x)$ is put in relation with the function $R(\vN(\x))$ that represents the reflectance map giving the value of the light reflection on the surface as a function of its orientation (i.e., of the normal $\vN(\x)$) at each point  $(\x,u(\x))$. 
Depending on how we describe the function $R$, different reflection models are determined. 
In the literature, the most common representation of $R$ is based on the Lambertian model (the L--model in the sequel) which takes into account only the angle between the outgoing normal to the surface $\vN(\x)$ and the light source $\vomega$, that is
\begin{equation}
I(\x) = \gamma_D(\x) \vN(\x) \cdot \vomega,
\end{equation}
where $\cdot$ denotes the standard scalar product between vectors  and $\gamma_D(\x)$ indicates the diffuse albedo, i.e. the diffuse reflectivity or reflecting power of a surface. It is the ratio of reflected radiation from the surface to incident radiance upon it. Its dimensionless nature is expressed as a percentage and it is measured on a scale from zero for no reflection of a perfectly black surface to 1 for perfect reflection of a white surface. The data are the grey-value image $I(\x)$, the direction of the light source represented by the unit vector $\vomega$ and the albedo $\gamma_D(\x)$. The light source $\vomega$ is a unit vector, hence $|\vomega| = 1$. In the simple case of a vertical light source, that is when the light source is in the direction of the vertical axis, this gives rise to an eikonal equation. 
Several questions arise, even in the simple case: is a single image sufficient to determine the surface? If not, which set of additional informations is necessary to have uniqueness? How can we compute an approximate solution? Is the approximation accurate? It is well known that for Lambertian surfaces there is no uniqueness and other informations are necessary to select a unique surface (e.g. the height at each point of local maximum for $I(\x)$). However, rather accurate schemes for the classical eikonal equation are now available for the approximation. Despite its simplicity, the Lambertian assumption is very strong and does not match with many real situations that is why we consider in this paper some non-Lambertian models trying to give a unified mathematical formulation for these models.\\
 In order to set this paper into a mathematical perspective, we should mention that the pioneering work of Horn \cite{Horn_PhD1970,H75} and his activity with his collaborators at MIT  \cite{HB86,HB89} produced the first formulation of  the Shape from Shading (SfS) problem via a partial differential equation (PDE) and a variational problem. These works have inspired many other contributions in this research area as one can see looking at the extensive list of references in the two surveys \cite{ZTCS99,DFS08}. Several approaches have been proposed, we can group them in two big classes (see the surveys \cite{ZTCS99,DFS08}): methods based on the resolution of PDEs and optimization methods based on a variational approximation. In the first group the unknown is directly  the height of the surface $z=u(\x)$,  one can find here rather old papers based on the method of characteristics \cite{JVD51,Rindfleisch_PE1966,H75,Oliensis91a,Oliensis91b,BCK92,Kozera97} where one typically looks for classical solutions.  More recently,  other contributions were developed in the framework of weak solutions in the viscosity sense starting from the seminal paper by  Rouy and Tourin \cite{RT92} and, one year later, by Lions-Rouy and Tourin \cite{LRT93} (see e.g.  \cite{IR95,CF96,CS99,FS97,CG00,PFR,Barles94,Sagona01,KO01,FSS03,PF03}). The second group  contains the contribution based on minimization methods for the variational problem where the unknown are the partial derivatives of the surface, $p=u_x$  and $q=u_y$ (the so-called normal vector field. See e.g. \cite{HB86,DD00,Strat79,FC88,Szeliski91,IH81,BH85,WH99}). It is important to note that in this approach one has to couple the minimization step to compute the normal field with a local reconstruction for $u$ which is based usually on a path integration. This necessary step has also been addressed by several authors (see \cite{DurouEMMCVPR2009} and references therein). 
We should also mention that a continuous effort has been made by the scientific community to take into account more realistic reflectance models \cite{BY94,AF07,RH05,VBW08,VVBW09}, different scenarios including perspective camera projection \cite{OD97,TSY05,PF_ICCV2003,AF06,VBLW09,Courteille_ICPR2004} 
and/or multiple images of the same object \cite{WNJ10,Yoon_IJCV2010}. The images can be taken from the same point of view but with different light sources as in the photometric stereo method \cite{Woo80,Kozera91,MF13,Tankus2005,MT13} or from different points of view but with the same light source as in stereo vision \cite{Cha94}. 
Recent works have considered more complicated scenarios, e.g. the case when light source is not at the optical center under perspective camera projection \cite{JTBBK13}. 
It is possible to consider in addition other supplementary issues, as the estimation of the albedo \cite{ZC91,BAC09,SH06,SH05} or of the direction of the light source  that are usually considered known quantities for the model but in practice are hardly available for real images. The role of boundary conditions which have to be coupled with the PDE is also a hard task. Depending on what we know, the model has to be adapted leading to a calibrated or uncalibrated problem (see \cite{Yoon_IJCV2010,WNJ10,FP12,QLD15} for more details). In this work we will assume that the albedo and the light source direction are given. \\
Regarding the modeling of non-Lambertian surfaces we also want to mention the important contribution of Ahmed and Farag  \cite{AF06}. These authors have adopted the modeling for SfS proposed by Prados and Faugeras \cite{PF_CVPR2005,PCF06} using a perspective projection where the light source is assumed to be located at the optical center of the camera instead at infinity and the light illumination is attenuated by a term $1/r^2$ ($r$ represents here the distance between the light source and the surface). They have derived the Hamilton-Jacobi (HJ) equation corresponding to the Oren-Nayar model, developed an approximation via the Lax-Friedrichs sweeping method. 
They gave there an experimental evidence that the non-Lambertian model seems to resolve the classical  concave/convex ambiguity in the perspective case if one includes the attenuation term $1/r^2$. In \cite{AF07} they extended their approach for various image conditions under orthographic and perspective projection, comparing their results for the orthographic L--model shown in \cite{Samaras_TPAMI2003} and in \cite{ZTCS99}. 
 Finally, we also want to mention the paper by Ragheb and Hancock \cite{Ragheb_3DPVT2004} where they treat a non-Lambertian model via a variational approach, investigating the reflectance models described by Wolff and by Oren and Nayar \cite{Wolff_JOSAA1994,WON98,ON_IJCV95}.\\

\vspace{-0.6cm}\noindent
\noindent \textbf{Our Contribution.}\\
In this paper we will adopt the PDE  approach in the framework of weak solutions for Hamilton-Jacobi equations. As we said, we will focus our attention on a couple of non-Lambertian reflectance models: the Oren-Nayar and the Phong models \cite{ON_SIGGRAPH94,ON_ECCV94,ON_CVPR93,ON_IJCV95,Phong75}. Both models are considered  for an orthographic projection and a single light source at infinity,  so no attenuation term is considered here. 
We are able to write the Hamilton-Jacobi equations in the same fixed point form useful for their analysis and approximation and, using the exponential change of variable introduced by Kru\v zkov  in \cite{Kru75}, we obtain natural upper bound for the solution. Moreover, we propose a semi-Lagrangian approximation scheme which can be applied to both the models, prove a convergence result for our  scheme that can be applied to this class of  Hamilton-Jacobi equations, hence to both non-Lambertian models. Numerical comparisons will show that our approach  is more accurate also for the 3D reconstructions of non-smooth surfaces. \\ 
A similar formulation for the Lambertian SfS problem with oblique light direction has been studied in \cite{FSS03} and here is extended to non-Lambertian models. We have reported some preliminary results just for the Oren-Nayar model in \cite{TF14}.\\

\vspace{-0.3cm}\noindent
\noindent \textbf{Organization of the paper.} \\
The paper is organized as follows.
After a formulation of the general model presented in Section \ref{sec:formulation_problem}, we present the SfS models starting from the classical Lambertian model (Section \ref{sec:Lambertian_model}). In Sections \ref{sec:ON_model} and \ref{sec:PH_model} we will give details on the construction of the nonlinear partial differential equation which corresponds respectively to the Oren-Nayar and the Phong models. Despite the differences appearing in these non-Lambertian models, we will be able to present them in a unified framework showing that the Hamilton-Jacobi equations for all the above models share a common structure. Moreover, the Hamiltonian appearing in these equations will always be convex in the gradient $\nabla u$. Then, in Section \ref{sec:general_convergence_theorem}, we will introduce our general approximation scheme which can be applied to solve this class of problems. 
In Section \ref{sec:numerical_tests} we will apply our approximation to a series of benchmarks based on synthetic and real images. We will discuss some  issue like accuracy, efficiency and the capability to obtain the maximal solution showing that the semi-Lagrangian approximation  is rather effective even for real images where several parameters are unknown. 
Finally, in the last section we will give a summary of the contributions of this work with some final comments and future research directions.

\section{Formulation of the general model} \label{sec:formulation_problem} 
\begin{figure}[h!]
\centering
\includegraphics[width=6.5cm]{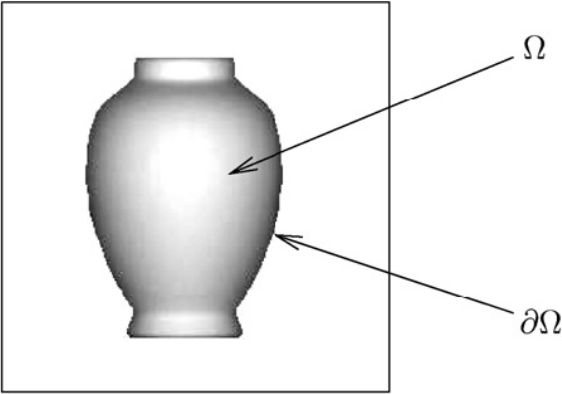}
\caption{An object on a flat background: $\Omega$ indicates the region inside the silhouette, $\partial \Omega$ the boundary of it.}
\label{fig:object}
\end{figure}

We fix a camera in a three-dimensional coordinate system (\emph{Oxyz}) in such a way that \emph{Oxy} coincides with the image plane and \emph{Oz} with the optical axis.
Let $\grass{\omega} = (\omega_1,\omega_2,\omega_3) = (\grass{\tilde{\omega}},\omega_3)\in\mathbb{R}^3$ (with $\omega_3 > 0$) be the unit vector that represents the direction of the light source (the vector points from the object surface to the light source); 
let $I(\x)$ be the function that measures the gray-level of the input image at the point $\x := (x,y)$. $I(\x)$ is the datum in the model since it is measured at each pixel of the image, for example in terms of a greylevel (from $0$ to $255$). In order to construct a continuous model, we will assume that $I(\x)$ takes real values in the interval $[0,1]$, defined in a compact domain $\overline{\Omega}$ called ``reconstruction domain" (with $\Omega\subset\mathbb{R}^2$ open set),  
$I:\overline{\Omega}\rightarrow[0,1]$, where the points with a value of $0$ are the dark point (blacks), while those with a value of $1$ correspond to a completely reflection of the light (white dots, with a maximum reflection). \\ 
We consider the following \emph{assumptions}: 
\begin{description}
 \item [A1.] there is a single light source placed at infinity in the direction $\grass{\omega}$ (the light rays are, therefore, parallel to each other);
 \item [A2.] the observer's eye is placed at an infinite distance from the object you are looking at (i.e. there is no perspective deformation);
 \item [A3.] there are no autoreflections on the surface.
\end{description}
In addition to these assumptions, there are other hypothesis that depend on the different reflectance models (we will see them in the description of the individual models).\\
Being valid the assumption $(A2)$ of orthographic projection, the visible part of the scene is a graph $z = u(\x)$ and the unit normal to the regular surface at the point corresponding to $\x$ is given by:
\begin{equation} \label{normal_N}
\vN(\x) = \frac{\n(\x)}{|\n(\x)|} = \frac{(- \nabla u(\x), 1)}{\sqrt{1 + |\nabla u(\x)|^2}},
\end{equation} 
where $\n(\x)$ is the outgoing normal vector.\\
We assume that the surface is standing on a flat background so the height function, which is the unknown of the problem, will be non negative, $u:\Omega \rightarrow[0,\infty)$. We will denote by $\Omega$ the region inside the silhouette and we will assume (just for technical reasons) that $\Omega$ is an open and bounded subset of $\R^2$ (see Fig. \ref{fig:object}).
It is well known that the SfS problem is described by the image irradiance equation \eqref{general_irradiance_eq} and depending on how we describe the function $R$ different reflection models are determined. We describe below some of them.
To this end, it would be useful to introduce a representation of the brightness function $I(\x)$ where we can distinguish different terms representing the contribution of ambient, diffuse reflected and specular reflected light. We will write then
\begin{equation}
\label{I_decomposition}
I(\x) = k_A I_A(\x)+k_D I_D(\x)+k_S I_S(\x),
\end{equation}
where $I_A(\x)$, $I_D(\x)$ and $I_S(\x)$ are respectively the above mentioned components and $k_A$, $k_D$ and $k_S$ indicate the percentages of these components such that their sum is equal to 1 (we do not consider absorption phenomena). Note that the diffuse or specular albedo is inside the definition of $I_D(\x)$ or $I_S(\x)$, respectively. 
In the sequel, we will always consider $I(\x)$ normalized in $[0,1]$. 
This will allow to switch on and off the different contributions depending on the model. 
Let us note that the ambient light term $I_A(\x)$ represents light everywhere in a given scene.  
In the whole paper we will consider it as a constant and we will neglect its contribution fixing $k_A = 0$. Moreover, for all the models presented below we will suppose uniform diffuse and/or specular albedo and we will put them equal to $1$, that is all the points of the surface reflect completely the light that hits them. We will omit them in what follows. 
As we will see in the following sections, the intensity of diffusely reflected light in each direction is proportional to the cosine of the angle $\theta_i$ between surface normal and light source direction, without taking into account the point of view of the observer, but another diffuse model (the Oren--Nayar model) will consider it in addition. The amount of specular reflected light towards the viewer is proportional to $(\cos\theta_s)^\alpha$, where $\theta_s$ is the angle between the ideal (mirror) reflection direction of the incoming light and the viewer direction, $\alpha$ being a constant modelling the specularity of the material. 
In this way we have a more general model and, dropping the ambient and specular component, we retrieve the Lambertian reflection as a special case.

\section{The Lambertian model (L--model)} \label{sec:Lambertian_model}
For a {\em Lambertian surface}, which generates a purely diffuse model, the specular component does not exist, then in \eqref{I_decomposition} we have just the diffuse component $I_D$ on the right side. 
Lambertian shading is view independent, hence the irradiance equation \eqref{general_irradiance_eq} becomes
\begin{equation} \label{Lamb_irradiance_eq}
I(\x) = \vN(\x) \cdot \vomega. 
\end{equation} 
Under these assumptions, the orthographic SfS problem consists in determining the function $u:\overline{\Omega}\rightarrow\mathbb{R}$ that satisfies the equation \eqref{Lamb_irradiance_eq}. The unit vector $\vomega$ and the function $I(\x)$ are the only quantities known.\\
For Lambertian surfaces  \cite{HB86,HB89}, just considering an orthographic projection of the scene,  it is possible to model the SfS problem via a nonlinear PDE of the first order which describes the relation between the surface $u(\x)$ (our unknown) and the brightness function $I(\x)$. 
In fact, recalling the definition of the unit normal to a graph given in \eqref{normal_N}, we can write \eqref{Lamb_irradiance_eq} as 
\begin{equation}
\label{Lamb_HJE}
I(\x) \sqrt{1 + |\nabla u(\x)|^2} + \widetilde{\vomega}\cdot \nabla u(\x) - \omega_3 = 0, \hbox{ in }\Omega
\end{equation}
where $\widetilde{\vomega} = (\omega_1,\omega_2)$. This is an Hamilton-Jacobi type equation which does not admit in general a regular solution. It is known that the mathematical framework to describe its weak solutions is the theory of viscosity solutions as in \cite{LRT93}.\\

\begin{figure}[h!]
\centering{
\includegraphics[width=8cm]{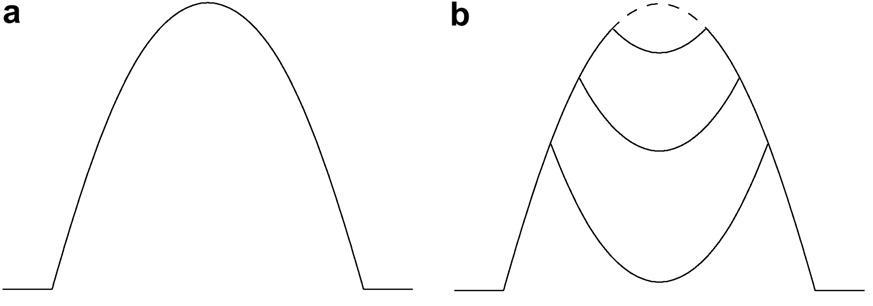}} 
\caption{Illustration of the concave/convex ambiguity: (a) maximal solution and (b) a.e. solutions giving the same image. Figure adapted from \cite{DFS08}.}
\label{fig:viscoSol}
\end{figure}

\noindent \textbf{The vertical light case.}\\ 
If we choose $ \grass{\omega} = (0,0,1)$, the equation \eqref{Lamb_HJE} becomes the so-called ``eikonal equation'':
\begin{equation}\label{Lamb_eiconale}
|\nabla u(\x)| = f(\x) \quad \textrm{for} \quad \x \in\Omega,
\end{equation}
where 
\begin{equation} \label{Lamb_f_eiconale}
f(\x) = \sqrt{\frac{1}{I(\x)^2} - 1}.
\end{equation}
The points $\x \in\Omega$ where $I(\x)$ assumes maximum value correspond to the case in which $\vomega$ and $\vN(\x)$ have the same direction: these points are usually called ``singular points''.\\
In order to make the problem well-posed, we need to add boundary conditions to the equations \eqref{Lamb_HJE} or \eqref{Lamb_eiconale}: they can require the value of the solution $u$ (Dirichlet boundary conditions type), or the value of its normal derivative (Neumann boundary conditions), or an equation that must be satisfied on the boundary (the so-called boundary conditions ``state constraint''). 
In this paper, we consider Dirichlet boundary conditions equal to zero assuming a surface on a flat background
\begin{equation}\label{DiriZERO}
u(\x) = 0, \quad \textrm{for} \quad \x \in\partial\Omega,
\end{equation}
but a second possibility of the same type occurs when it is known the value of $u$ on the boundary, which leads to the more general condition
\begin{equation}\label{Dirichlet}
u(\x) = g(\x), \quad \textrm{for} \quad \x \in\partial\Omega.
\end{equation}
Unfortunately, adding a boundary condition to the PDE that describes the SfS model is not enough to obtain a unique solution because of the concave/convex ambiguity. 
In fact, the Dirichlet problem \eqref{Lamb_HJE}-\eqref{Dirichlet} can have several weak solutions in the viscosity sense and also several classical solutions due to this ambiguity (see \cite{H75}). As an example, all the surfaces represented in Fig.  \ref{fig:viscoSol} are viscosity solutions of the same problem \eqref{Lamb_eiconale}-\eqref{DiriZERO} which is a particular case of \eqref{Lamb_HJE}-\eqref{Dirichlet} (in fact the equation is $|u'|=-2x$ with homogenous Dirichlet boundary condition). The solution represented in Fig.  \ref{fig:viscoSol}-a is the maximal solution and is smooth. All the non-smooth a.e. solutions, which can be obtained by a reflection with respect to a horizontal axis, are still admissible weak solutions (see Fig. \ref{fig:viscoSol}-b). In this example, the lack of uniqueness of the viscosity solution is due to the existence of a singular point where the right hand side of (\ref{Lamb_eiconale}) vanishes. An additional effort is then needed to define which is the preferable solution since the lack of uniqueness is also a big drawback when trying to compute a numerical solution. In order to circumvent these difficulties, the problem is usually solved by adding some information such as height at each singular point \cite{LRT93}.

For analytical and numerical reasons it is useful to introduce the exponential Kru\v zkov change of variable \cite{Kru75} $\mu v(\x) = 1- e^{-\mu u(\x)}$. 
In fact, setting the problem in the new variable $v$  we will have  values in $[0, 1/\mu]$ instead of $[0,\infty)$ as the original variable $u$ so an upper bound will be easy to find. 
Note that $\mu$ is a free positive parameter which does not have a specific physical meaning in the SfS problem. However,  it can play an  important role also in our convergence proof as we will see later (see the remark following the end of Theorem \ref{th:genprop}). 
Assuming that the surface is standing on a flat background and following \cite{FSS03}, we can write \eqref{Lamb_HJE} and \eqref{DiriZERO} in a fixed point form in the new variable $v$. To this end let us define 
$\vb^L:\Omega\times\partial B_3(0,1)\rightarrow \R^2$ and $f^L:\Omega\times\partial B_3(0,1)\times [0,1] \rightarrow \R$ as 
\begin{equation} \label{Lamb_vector_field}
\vb^L(\x,a) := \frac{1}{\omega_3} \left(I(\x) a_1 - \omega_1, I(\x) a_2 - \omega_2 \right), 
\end{equation}
\begin{equation} \label{Lamb_f}
f^L(\x, a, v(\x)) :=- \frac{I(\x) a_3}{\omega_3} (1-\mu v(\x)) +1
\end{equation}
and let $B_3$ denote the unit ball in $\R^3$. We obtain

\hspace{-0.6cm}
\begin{tabular}{|r|}\hline
\begin{minipage}{0.95\columnwidth}

\vspace{2mm}

{\it Lambertian Model} \\

\vspace{-6mm}

\begin{equation}
\label{Lamb_pde_v}
\left\{ \begin{array}{ll}
\mu v(\x) = T^L(\x,v(\x),\nabla v) & \hbox{for }\x \in \Omega,\\
v(\x) = 0 & \hbox{for }\x \in \partial \Omega,
\end{array} \right.
\end{equation}
\vspace{2mm}
\hbox{where}\\
$T^L(\cdot):=\min\limits_{a\in\partial B_3}  \{ \vb^L(\x,a) \cdot \nabla v(\x) +f^L(\x, a, v(\x))\}$.
\vspace{0.2cm}
\end{minipage}\\
\hline
\end{tabular}
\vspace{2mm}\\
It is important to note for the sequel that the structure of the above first order Hamilton-Jacobi equation is similar to that related to the dynamic programming approach in control theory where  $\vb$ is a vector field describing the dynamics of the system and $f$ is a running cost. In that framework the meaning of $v$ is that of a value function which allows to characterize the optimal trajectories (here they play the role of characteristic curves). The interested reader can find more details on this interpretation in \cite{FF14}.

\section{The Oren-Nayar model (ON--model)} \label{sec:ON_model}
\begin{figure}[h!]
\centering
\includegraphics[width=5.5cm]{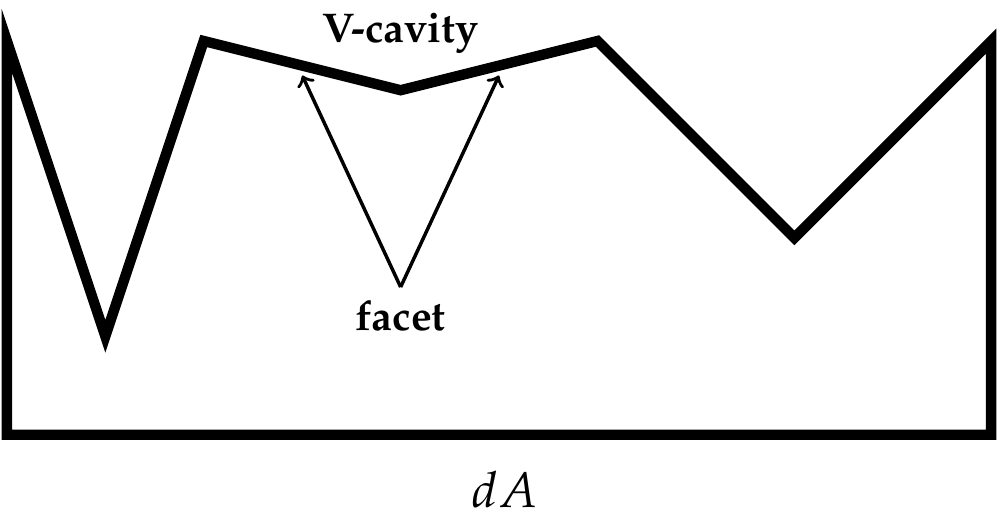}\\
\caption{Facet model for surface patch $dA$ consisting of many V-shaped Lambertian cavities. Figure adapted from \cite{JBBG12}.}
\label{fig:dA}
\end{figure}

\noindent The diffuse reflectance ON--model~\cite{ON_SIGGRAPH94,ON_ECCV94,ON_CVPR93,ON_IJCV95} is an extension of the previous L--model which
explicitly allows to handle {\em rough} surfaces. 
The idea of this model is to represent a rough surface as an aggregation of V-shaped cavities, each with Lambertian reflectance properties (see Fig. \ref{fig:dA}). \\
In \cite{ON_SIGGRAPH94} and, with more details, in \cite{ON_IJCV95}, Oren and Nayar derive a reflectance model for several type of surfaces with different slope-area distributions. 
In this paper we will refer to the model called by the authors the ``Qualitative Model'', a simpler version obtained by ignoring interreflections (see Section $4.4$ of \cite{ON_SIGGRAPH94} for more details).

\begin{figure}[h!]
\centering
\includegraphics[width=10.cm]{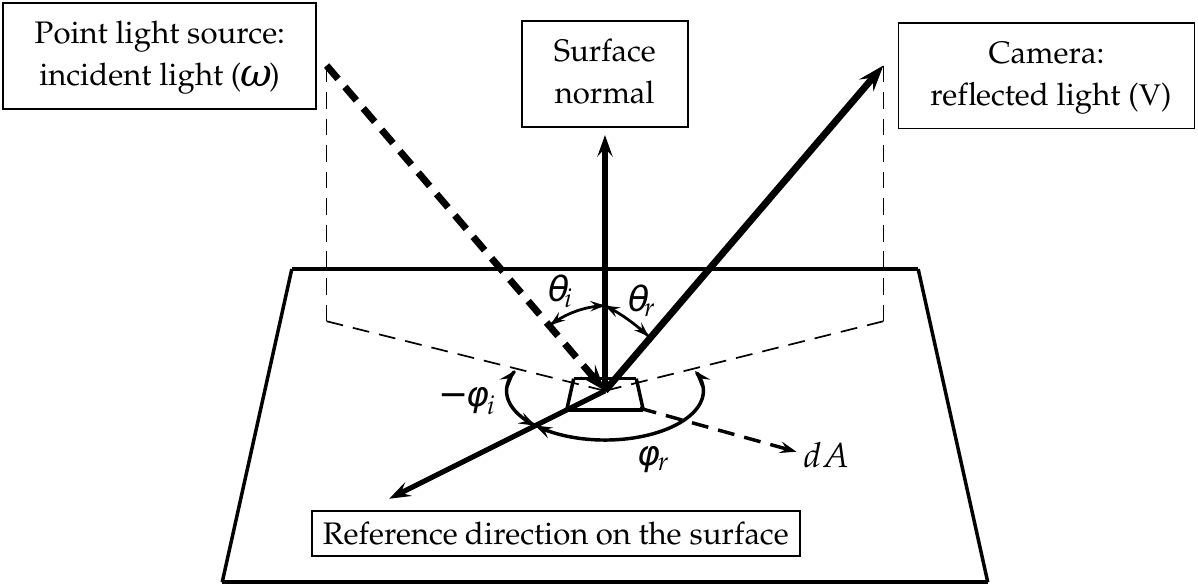}
\caption{Diffuse reflectance for the ON--model. Figure adapted from \cite{JBBG12}.}
\label{fig:angles_ON} 
\end{figure}

\noindent Assuming that there is a linear relation between the irradiance of the image and the image intensity,
the {\em $I_D$ brightness equation} for the ON--model  is given by 
\begin{equation}
\label{ON_reflectance_model}
I_D(\x)= \cos(\theta _i) (A + B \sin(\alpha) \tan (\beta )M(\varphi_i,\varphi_r)) 
\end{equation}
\begin{eqnarray} \label{eq:A_B}
&&\hbox{where  }A = 1 - 0.5 \, \sigma^2 (\sigma^2 + 0.33)^{-1} \\
&& \hspace{1cm} B = 0.45\sigma^2(\sigma^2 + 0.09)^{-1} \\
&& \hspace{1cm} M(\varphi_i,\varphi_r) = \max\{0,\cos(\varphi_r -  \varphi_i)\}.
\end{eqnarray}
Note that $A$ and $B$ are two non negative constants depending on the statistics of the cavities via the roughness parameter $\sigma$. We  set $\sigma \in [0,\pi/2)$, interpreting $\sigma$ as the slope of the cavities. In this model (see Fig. \ref{fig:angles_ON}), $\theta_{i}$ represents the angle between the unit normal to the surface $\vN(\x)$ and the light source direction $\vomega$,
$\theta_{r}$ stands for the angle between $\vN(\x)$ and the observer direction $\vV$,
$\varphi_{i}$ is the angle between the projection of the light source direction $\vomega$ and the $x_1$ 
axis onto the $(x_1, x_2)$-plane,  
$\varphi_{\mathrm r}$ denotes the angle between the projection of the observer direction $\vV$ and the $x_1$ axis
onto the $(x_1, x_2)$-plane  
and the two variables $\alpha$ and $\beta$ are given by
\begin{equation}\label{def:ab}
\alpha = \max \left\{ \theta_{i}, \theta_{r} \right\} \hbox{ and } \beta = \min \left\{ \theta_{i}, \theta_{r} \right\}.
\end{equation}
Since the vectors $\vomega$ and $\vV$ are fixed and given, their projection on the incident plane is obtained considering their first two components over three (see Eq. \eqref{formulas_trig_cos_phi}). In this way, the quantity $\max\{0,\cos(\varphi_r -  \varphi_i)\}$ is computed only once for a whole image. 

\noindent We define (see Fig. \ref{fig:angles_ON}):
\begin{eqnarray} \label{formulas_trig_cos_theta_i}
\cos(\theta_i) &=& \vN \cdot \vomega = \frac{-\widetilde{\vomega}\cdot \nabla u(\x) + \omega_3}{\sqrt{1 + |\nabla u(\x)|^2}} \vspace{0.2cm} \\
\cos(\theta_r) &=& \vN \cdot \vV = \frac{-\widetilde{\mathbf{v}}\cdot \nabla u(\x) + v_3}{\sqrt{1 + |\nabla u(\x)|^2}} \vspace{0.3cm} \label{formulas_trig_cos_theta_r} 
\\
\cos(\varphi_r &-&  \varphi_i) = (\omega_1,\omega_2) \cdot (v_1,v_2) =  \widetilde{\vomega}\cdot \widetilde{\mathbf{v}}   \vspace{0.2cm} \label{formulas_trig_cos_phi} \\
\sin(\theta_i) &=& \sqrt{1- (\cos(\theta_i))^2} =  \frac{g_{\omega}(\nabla u(\x))}{\sqrt{1 + |\nabla u(\x)|^2}} \vspace{0.2cm} \label{formulas_trig_sin_theta_i} \\
\sin(\theta_r) &=& \sqrt{1- (\cos(\theta_r))^2} =  \frac{g_{v}(\nabla u(\x))}{\sqrt{1 + |\nabla u(\x)|^2}} \label{formulas_trig_sin_theta_r} 
\end{eqnarray}
where
\begin{eqnarray}
g_{\omega}(\nabla u(\x))&:=& \sqrt{1 + |\nabla u(\x)|^2 - (-\widetilde{\vomega}\cdot \nabla u(\x) + \omega_3)^2} \vspace{0.4cm} \nonumber \\ 
g_{v}(\nabla u(\x))&:=& \sqrt{1 + |\nabla u(\x)|^2 - (-\widetilde{\mathbf{v}}\cdot \nabla u(\x) + v_3)^2}. \nonumber
\end{eqnarray}

For smooth surfaces, we have $\sigma=0$ and in this case the ON--model reduces to the L--model.  
In the particular case $\vomega = \vV = (0,0,1)$, or, more precisely, when $\cos(\varphi_r -  \varphi_i) \leq 0$ (e.g. the case when the unit vectors $\vomega$ and $\vV$ are perpendicular we get $\cos(\varphi_r -  \varphi_i) = -1$) the equation \eqref{ON_reflectance_model} simplifies and reduces to a L--model scaled by the coefficient $A$. This means that the model is more general and flexible than the L--model. 
This happens when only one of the two unit vectors is zero or, 
more in general, when the dot product between the normalized projections onto the $(x_1, x_2)$-plane of $\vomega$ and $\vV$ is equal to zero.

Also for this diffuse model we neglect the ambient component, setting $k_D = 1$. As a consequence, in the general equation \eqref{I_decomposition} the total light intensity $I(\x)$ is equal to the diffuse component $I_D(\x)$ (described by the equation \eqref{ON_reflectance_model}). This is why we write $I(\x)$ instead of $I_D(\x)$ in what follows.

To deal with this equation one has to compute the $\min$ and $\max$ operators which appear in \eqref{ON_reflectance_model} and \eqref{def:ab}. 
Hence, we must consider several cases described in detail in what follows. For each case we will derive a partial differential equation that is always a first order nonlinear HJ equation: \\ 
\\
{\bf Case 1:} $\theta_i \ge \theta_r$ and $(\varphi_r -  \varphi_i) \in [0, \frac{\pi}{2}) \cup (\frac{3}{2}\pi, 2 \pi]$ \vspace{0.2cm}  \\
The brightness equation \eqref{ON_reflectance_model} becomes
\begin{equation}
\label{ON_reflectance_model_case1}
I(\x)=\cos(\theta_i)\left(A \! + \! B \sin(\theta_i) \frac{\sin(\theta_r)}{\cos(\theta_r)} \cos(\varphi_r -  \varphi_i) \right) \hspace{-0.2cm}
\end{equation}
and by using the formulas \eqref{formulas_trig_cos_theta_i}-\eqref{formulas_trig_sin_theta_r} we arrive to the following HJ equation 
\begin{equation} 
\begin{array}{ll}
\label{ON_pde_u_case1}
& I(\x) (\sqrt{1 + |\nabla u(\x)|^2}) \vspace{0.2cm} \displaystyle + A (\widetilde{\vomega}\cdot \nabla u(\x) - \omega_3)\vspace{0.2cm}  \\
&\displaystyle - B \frac{(\widetilde{\vomega}\cdot \widetilde{\mathbf{v}}) \, g_{\omega}(\nabla u(\x)) \, g_{v}(\nabla u(\x)) (-\widetilde{\vomega}\cdot \nabla u(\x) + \omega_3)}{\sqrt{1 + |\nabla u(\x)|^2} \,\, (- \widetilde{\mathbf{v}}\cdot \nabla u(\x) +v_3)} = 0, 
\end{array}
\end{equation}
where $\widetilde{\vomega} := (\omega_1,\omega_2)$ and $\widetilde{\mathbf{v}} := (v_1, v_2)$.  \\
\\
\\
{\bf Case 2:} $\theta_i < \theta_r$ and $(\varphi_r -  \varphi_i) \in [0, \frac{\pi}{2}) \cup (\frac{3}{2}\pi, 2 \pi]$ \vspace{0.2cm} \\
In this case the brightness equation \eqref{ON_reflectance_model} becomes
\begin{equation}
\label{ON_reflectance_model_case2}
I(\x) =  \cos(\theta_i) \left(A \! + \! B \sin(\theta_r) \frac{\sin(\theta_i)}{\cos(\theta_i)} \cos(\varphi_r -  \varphi_i) \right) \hspace{-0.2cm}
\end{equation}
and by using again the formulas \eqref{formulas_trig_cos_theta_i}-\eqref{formulas_trig_sin_theta_r} we get 
\begin{equation} 
\label{ON_pde_u_case2}
\begin{array}{ll}
I(\x)& (1 + |\nabla u(\x)|^2) \vspace{0.2cm} \\
& + A (\widetilde{\vomega}\cdot \nabla u(\x) - \omega_3) \, \sqrt{1 + |\nabla u(\x)|^2} \vspace{0.2cm} \\
& - B (\widetilde{\vomega}\cdot \widetilde{\mathbf{v}}) \, g_{\omega}(\nabla u(\x)) \, g_{v}(\nabla u(\x)) = 0. \vspace{0.2cm}
\end{array}
\end{equation}
{\bf Case 3:} $\forall \, \theta_i, \theta_r$ and $(\varphi_r -  \varphi_i) \in [\frac{\pi}{2}, \frac{3}{2}\pi]$ \vspace{0.2cm} \\
In this case we have the implication $\max \{0, \cos(\varphi_r -  \varphi_i)\} = 0$.
The brightness equation \eqref{ON_reflectance_model} simplifies in
\begin{equation}
\label{ON_reflectance_model_case3} 
I(\x) = A \, \cos(\theta_i) 
\end{equation}
and the HJ equation associated to it becomes consequentially 
\begin{equation} 
\label{ON_pde_u_case3}
I(\x) (\sqrt{1 + |\nabla u(\x)|^2}) + A (\widetilde{\vomega}\cdot \nabla u(\x) - \omega_3) = 0,
\end{equation}
that is equal to the L--model scaled by the coefficient $A$.\\
\\
{\bf Case 4:} $\theta_i = \theta_r$ and $\varphi_r =  \varphi_i$ \vspace{0.2cm} \\
This is a particular case when the position of the light source $\vomega$ coincides with the observer direction $\vV$ but there are not on the vertical axis.
This choice implies $\max\{0,\cos(\varphi_i - \varphi_r)\} = 1$, then defining $\theta := \theta_i = \theta_r = \alpha = \beta$, the equation \eqref{ON_reflectance_model} simplifies to
\begin{equation}
\label{ON_reflectance_model_case4}
I(\x) =  \cos(\theta) \ \left(A \! + \! B \sin(\theta)^2\cos(\theta)^{-1} \right)
\end{equation}
and we arrive to the following HJ equation 
\begin{equation} 
\label{ON_pde_u_case_4}
\begin{array}{ll}
(I(\x)&- B) (\sqrt{1 + |\nabla u(\x)|^2}) + A (\widetilde{\vomega}\cdot \nabla u(\x) - \omega_3) \vspace{0.2cm}\\
&\displaystyle + B \frac{(-\widetilde{\vomega}\cdot \nabla u(\x) + \omega_3)^2}{\sqrt{1 + |\nabla u(\x)|^2}} = 0. \hspace{-0.2cm}
\end{array}
\end{equation}
Note that this four cases are exactly the same cases reported and analyzed in \cite{JTBBK13}. This is not surprising since the reflectance model used there is always the same one proposed by Oren and Nayar. However, here we get different HJ equations since we consider an orthographic camera projection and cartesian coordinates whereas in \cite{JTBBK13}  the HJ equations are derived in spherical coordinates under a generalized perspective camera projection. Another major difference is that in that paper the light source is close to the camera but is not located at the optical center of the camera.

\noindent \textbf{The vertical light case.} \\
If $\vomega = (0,0,1)$, independently of the position of $\vV$, the analysis is more simple. In fact, the first three cases considered above reduce to a single case corresponding to the following simplified PDE for the brightness equation \eqref{ON_reflectance_model}
\begin{equation} 
\label{ON_reflectance_model_vertical}
I(\x) = \displaystyle \frac{A}{\sqrt{1 + |\nabla u(\x)|^2}} .
\end{equation}
In this way we can put it in the following eikonal type equation, analogous to the Lambertian eikonal equation \eqref{Lamb_eiconale}:
\begin{equation}\label{ON_eiconale}
|\nabla u(\x)| = f(\x) \quad \textrm{for} \quad \x \in\Omega,
\end{equation}
where 
\begin{equation} \label{eq:def_f_ON_eiconale}
f(\x) = \displaystyle \sqrt{\frac{A^2}{I(\x)^2} - 1} .
 \end{equation}
 
Following \cite{TF14}, we write the surface as $S(\x,z) = z - u(\x)=0$, for $\x \in \Omega$, $z\in\R$, and $\nabla S(\x,z) = (-\nabla u(\x),1)$, so \eqref{ON_pde_u_case_4} becomes
\begin{equation}
\label{ON_pde_S}
\begin{array}{ll}
( I(\x) - B)& |\nabla S(\x,z)| \vspace{0.2cm}\\ 
& + A (- \nabla S(\x,z)\cdot \vomega)\vspace{0.2cm} \\
& + B \left(\frac{\nabla S(\x,z)}{|\nabla S(\x,z)|} \cdot \vomega \right)^2 |\nabla S(\x,z)| = 0.
\end{array}
\end{equation}
Defining 
\begin{equation}\label{def_d}
\vd(\x,z) := \nabla S(\x,z) / |\nabla S(\x,z)| 
\end{equation}
and 
\begin{equation}
c(\x,z) := I(\x) - B +B (\vd(\x,z) \cdot \vomega)^2 ,
\end{equation}
using the equivalence 
\begin{equation}
|\nabla S(\x,z)| \equiv \max_{a\in\partial B_3} \{a\cdot \nabla S(\x,z)\}
\end{equation}
 we get
\begin{equation}
\max_{a\in\partial B_3} \{ c(\x,z) \, a\cdot \nabla S(\x,z) -A\vomega\cdot \nabla S(\x,z)\} = 0.
\end{equation}
Let us define the vector field for the ON-model
\begin{equation}
\vb^{ON}(\x,a) := \frac{\left( c(\x,z) a_1 - A\omega_1, c(\x,z) a_2 - A\omega_2 \right)}{A\omega_3} , 
\end{equation}
and
\begin{equation}\label{f_ON}
f^{ON}(\x,z,a,v(\x)) := - \frac{c(\x,z) a_3}{A\omega_3} (1-\mu v(\x)) + 1.
\end{equation} 
Then, introducing the exponential Kru\v zkov change of variable $\mu v(\x) = 1- e^{-\mu u(\x)}$ as already done for the L--model, 
we can finally write the fixed point problem in the new variable $v$ obtaining the 
\vspace{1mm}

\noindent\begin{tabular}{|r|}\hline 
\begin{minipage}{0.97\columnwidth}

\vspace{2mm}

{\it Oren-Nayar Model} \\

\vspace{-6mm}

\begin{equation}
\label{ON_pde_v}
\left\{ \begin{array}{ll}
\mu v(\x) = T^{ON}(\x,v(\x),\nabla v),
& \textrm{$\x \in \Omega$}, \hspace{-0.cm}\\
v(\x) = 0, & \textrm{$\x \in \partial \Omega$.} \hspace{-0.cm} 
\end{array} \right.
\end{equation}
\vspace{1mm}
where
\[ T^{ON}(\cdot):=\hspace{-1.1mm}\min\limits_{a\in\partial B_3}  \{ \vb^{ON}(\x,a)  \cdot \nabla v(\x) 
+ f^{ON}(\x,z,a,v(\x)) \} \]
\vspace{1mm}
\end{minipage}\\
\hline
\end{tabular}
\vspace{1mm}\\
Note that the simple homogeneous Dirichlet boundary condition is due to the flat background behind the object but a condition like $u(\x) = g(\x)$ can also be considered if necessary. Moreover, the structure is similar to the previous Lambertian model although the definition of the vector field and of the cost are different.
\\
In the particular case when $\cos(\varphi_r -  \varphi_i) = 0$, 
the equation \eqref{ON_reflectance_model} simply reduces to
\begin{equation}
\label{ON_reflectance_model_case_remark1}
I(\x) =  A \cos(\theta)
\end{equation}
and, as a consequence, the Dirichlet problem in the variable $v$ is equal to \eqref{ON_pde_v} with $c(\x,z) = I(\x)$.

\section{The Phong model (PH--model)}  \label{sec:PH_model}
The PH--model is an empirical model that was developed by Phong \cite{Phong75} in 1975. This model introduces a specular component to the brightness function $I(\x)$, representing the diffuse component $I_D(\x)$ in \eqref{I_decomposition} as the Lambertian reflectance model.\\
A simple specular model is obtained putting the incidence angle equal to the reflection one and $\vomega$, $\vN(\x)$ and $\vR(\x)$ belong to the same plane.\\
This model describes the specular light component $I_S(\x)$ as a power of the cosine of the angle between the unit vectors $\vV$ and $\vR(\x)$ (it is the vector representing the reflection of the light $\vomega$ on the surface). 
Hence, the brightness equation for the PH--model is
\begin{equation}
\label{eq:phong_model}
I(\x) = k_D (\vN(\x) \cdot \vomega) + k_S (\vR(\x)\cdot \vV)^{\alpha},
\end{equation}
where the parameter $\alpha\in[1,10]$ is used to express the specular reflection properties of a material and $k_D$ and $k_S$ indicate the percentages of diffuse and specular components, respectively. 
Note that the contribution of the specular part decreases as  the value of $\alpha$ increases. 

\begin{figure}[h!]
\centering{
\includegraphics[width=7.5cm]{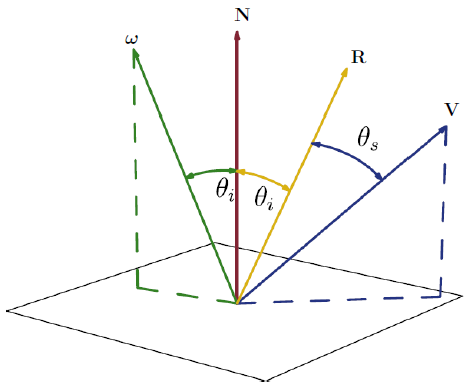}}
\caption{Geometry of the Phong reflection model.}
\label{fig:Phong_geometry}
\end{figure}

Starting to see in details the PH--model in the case of oblique light source $\vomega$ and oblique observer $\vV$. \\
Assuming that $\vN(\x)$ is the bisector of the angle between $\vomega$ and $\vR(\x)$ (see Fig. \ref{fig:Phong_geometry}), we obtain
\begin{equation}
\vN(\x) = \frac{\vomega + \vR(\x)}{||\vomega + \vR(\x)||}  
\end{equation}
which implies
\begin{equation}
\vR(\x) = ||\vomega + \vR(\x)|| \vN(\x) - \vomega.
\end{equation}
From the parallelogram law, taking into account that $\vomega, \vR(\x)$ and $\vN(\x)$ are unit vectors,
we can write $||\vomega + \vR(\x)|| = 2(\vN(\x) \cdot \vomega)$, then we can derive the unit vector $\vR(\x)$ as follow:
\begin{eqnarray} \label{vector_R}
&&\vR(\x) = 2(\vN(\x) \cdot \vomega) \vN(\x) - \vomega \nonumber\\
&&= \displaystyle 2 \left(\frac{-\widetilde{\vomega}\cdot \nabla u(\x) + \omega_3}{\sqrt{1 + |\nabla u(\x)|^2}} \right) \vN(\x) - (\omega_1,\omega_2, \omega_3) \\
&&= \displaystyle \left(\frac{-2\widetilde{\vomega}\cdot \nabla u(\x) + 2\omega_3}{1 + |\nabla u(\x)|^2} \right) \left(-\nabla u(\x), 1\right) - (\omega_1,\omega_2, \omega_3).\nonumber
\end{eqnarray}
\\
With this definition of the unit vector $\vR(\x)$ we have
\begin{eqnarray}
\label{R_cdot_V_omega_e_V_oblique}
&&\vR(\x)\cdot\vV = \\
&& \displaystyle \left( \frac{-2\widetilde{\vomega}\cdot \nabla u(\x) + 2\omega_3}{1 + |\nabla u(\x)|^2} \right) \left(-\nabla u(\x)\cdot \widetilde{\mathbf{v}}+ v_3 \right) - \vomega \cdot \vV. \nonumber
\end{eqnarray}
Then, setting $\alpha=1$, that represents the worst case, equation \eqref{eq:phong_model} becomes
\begin{eqnarray}
\label{eq:PH_pde_u_omega_e_V_oblique} 
&&I(\x)(1 + |\nabla u(\x)|^2) \nonumber\\
&&- k_D (-\nabla u(\x)\cdot \widetilde{\vomega} + \omega_3) (\sqrt{1 + |\nabla u(\x)|^2})  \\
&&- 2 k_S ( -\nabla u(\x)\cdot \widetilde{\vomega} + \omega_3) \left( -\widetilde{\mathbf{v}}\cdot \nabla u(\x) + v_3 \right) \nonumber\\
&& + k_S ( \vomega \cdot \vV) (1 + |\nabla u(\x)|^2)= 0,\nonumber
\end{eqnarray}
to which we add a Dirichlet boundary condition equal to zero, assuming that the surface is standind on a flat blackground.
Note that the cosine in the specular term is usually replaced by zero if $\vR(\x)\cdot \vV < 0$ (and in that case we get back to the L--model).\\
As we have done for the previous models, we write the surface as $S(\x,z) = z - u(\x)=0$, for $\x \in \Omega$, $z\in\R$, and $\nabla S(\x,z) = (-\nabla u(\x),1)$, so \eqref{eq:PH_pde_u_omega_e_V_oblique} will be written as 
\begin{eqnarray}\label{PH_pde_S}
&& ( I(\x) + k_S (\vomega\cdot\vV) ) |\nabla S(\x,z)|^2 \nonumber\\
&& - k_D (\nabla S(\x,z)\cdot \vomega) (|\nabla S(\x,z)|) \\
&& - 2 k_S (\nabla S(\x,z)\cdot \vomega) (\nabla S(\x,z)\cdot \vV) = 0.\nonumber 
\end{eqnarray}
Dividing by $|\nabla S(\x,z)|$, defining $\vd(\x,z)$ as in \eqref{def_d} and
$c(\x) :=  I(\x) + k_S (\vomega\cdot\vV)$, 
 we get
\begin{eqnarray}\label{PH_pde_S_con_c_e_d}
&& c(\x) |\nabla S(\x,z)| - k_D (\nabla S(\x,z)\cdot \vomega) \\
 &&- 2 k_S (\nabla S(\x,z)\cdot \vomega) (\vd(\x,z)\cdot \vV) = 0. \nonumber 
 \end{eqnarray}
By the equivalence
$ \displaystyle |\nabla S(\x,z)| \equiv \max_{a\in\partial B_3} \{a\cdot \nabla S(\x,z)\}$ we obtain
\begin{eqnarray}
&&\max_{a\in\partial B_3} \{ c(\x) \, a\cdot \nabla S(\x,z) - k_D (\nabla S(\x,z)\cdot \vomega ) \\
&&- 2 k_S (\nabla S(\x,z)\cdot \vomega) (\vd(\x,z)\cdot \vV) \} = 0.\nonumber
\end{eqnarray}
Let us define the vector field
\begin{equation}\label{b^PH}
\vb^{PH}(\x,a) :=  \frac{1}{Q^{PH}(\x,z)}  \vM^{PH}(\x,z)
\end{equation}
where
\begin{equation}
Q^{PH}(\x,z) := 2 k_S \omega_3 (\vd(\x,z)\cdot \vV) + k_D \omega_3,
\end{equation}
and 
\begin{eqnarray}
&&\vM_i^{PH}(\x,z) := \left( c(\x) a_i - k_D \omega_i \right . \nonumber\\
&&\left. \hspace{1.7cm} -2 k_S \, \omega_i (\vd(\x,z)\cdot \vV) \right), i=1,2. 
\end{eqnarray}

\noindent Let us also define 
\begin{equation}
f^{PH}(\x,z,a,v(\x)) := - \frac{c(\x) a_3}{Q^{PH}(\x,z)} (1-\mu v(\x)) +1.
\end{equation}
Again, using the exponential Kru\v zkov change of variable $\mu v(\x) = 1- e^{-\mu u(\x)}$ as done for the previous models,
we can finally write the nonlinear fixed point problem 
\vspace{1mm}

\noindent\begin{tabular}{|r|}\hline 
\begin{minipage}{0.97\columnwidth}

\vspace{2mm}

{\it Phong Model} \\

\vspace{-6mm}

\begin{equation}
\label{PH_pde_v}
\left\{ \begin{array}{ll}
\mu v(\x) = T^{PH}(\x,v(\x),\nabla v),
& \textrm{$\x \in \Omega$}, \hspace{-0.cm}\\
v(\x) = 0, & \textrm{$\x \in \partial \Omega$.} \hspace{-0.cm} 
\end{array} \right.
\end{equation}
\vspace{1mm}
where
\[ T^{PH}(\cdot):=\hspace{-1.1mm}\min\limits_{a\in\partial B_3}  \{ \vb^{PH}(\x,a)  \cdot \nabla v(\x) 
+ f^{PH}(\x,z,a,v(\x)) \} \]
\vspace{1mm}
\end{minipage}\\
\hline
\end{tabular}
\vspace{1mm}\\
%
\noindent \textbf{A. Oblique light source and vertical position of the observer.} \\
In the case of oblique light source $\vomega$ and vertical observer $\vV = (0,0,1)$, the dot product $\vR(\x) \cdot \vV$ becomes
\begin{equation}
\label{R_cdot_V_vertical_V}
\begin{array}{ll}
\vR(\x)\cdot\vV &= \displaystyle \frac{-2\widetilde{\vomega}\cdot \nabla u(\x) + 2\omega_3}{1 + |\nabla u(\x)|^2} - \omega_3 \vspace{1.0mm}\\
& \displaystyle = \frac{-2\widetilde{\vomega}\cdot \nabla u(\x) + \omega_3 (1 - |\nabla u(\x)|^2)}{1 + |\nabla u(\x)|^2}.
\end{array}
\end{equation}
The fixed point  problem in $v$ will be equal to \eqref{PH_pde_v} with the following choices
\begin{eqnarray}
\begin{array}{lll}
c(\x) &:=& I(\x) + \omega_3 k_S, \vspace{2mm} \\
Q^{PH}(\x,z) &:=& 2 k_S (\vd(\x,z)\cdot \vomega) + k_D \omega_3, \vspace{2mm} \\
 \vb^{PH}(\x,a) &:=& \frac{\left( c(\x) a_1 - k_D \omega_1, c(\x) a_2 - k_D \omega_2 \right)}{Q^{PH}(\x,z)}
\end{array}
\end{eqnarray}
\noindent \textbf{B. Vertical light source and oblique position of the observer.} \\
When $\vomega = (0,0,1)$ the definition of the vector $\vR(\x)$ reported in \eqref{vector_R} becomes
\begin{equation}
\label{vector_R_with_vertical_omega}
\vR(\x) = \displaystyle \left(\frac{-2 \nabla u(\x)}{1 + |\nabla u(\x)|^2}, \frac{2}{1 + |\nabla u(\x)|^2} - 1 \right)
\end{equation}
and, as a consequence, the dot product $\vR(\x) \cdot \vV$ with general $\vV$ is 
\begin{equation}
\label{R_cdot_V_vertical_omega}
\vR(\x)\cdot\vV = \frac{-2\widetilde{\mathbf{v}}\cdot \nabla u(\x) + v_3 (1 - |\nabla u(\x)|^2)}{1 + |\nabla u(\x)|^2}.
\end{equation}
Hence, the fixed point  problem in $v$ is equal to \eqref{PH_pde_v} with
\begin{eqnarray}
\begin{array}{lll}
c(\x) &:=& I(\x) + v_3 k_S, \vspace{2mm} \\
Q^{PH}(\x,z) &:=& 2 k_S (\vd(\x,z)\cdot \vV) + k_D, \vspace{2mm} \\
 \vb^{PH}(\x,a) &:=& \frac{1}{Q^{PH}(\x,z)} \left( c(\x) a_1, c(\x) a_2 \right).
\end{array}
\end{eqnarray} 
\\
\noindent \textbf{C. Vertical light source and vertical position of the observer.} \\
If we choose $\vomega \equiv \vV = (0,0,1)$ the equation \eqref{eq:PH_pde_u_omega_e_V_oblique} simplifies in
\begin{eqnarray}\label{eq:PH_pde_u_vertical}
&&I(\x) (1 + |\nabla u(\x)|^2) - k_S (1 - |\nabla u(\x)|^2) \\
&&\quad - k_D (\sqrt{1 + |\nabla u(\x)|^2}) = 0.\nonumber
\end{eqnarray}
Working on this equation one can put it in the following eikonal type form, which is analogous to the Lambertian eikonal equation \eqref{Lamb_eiconale}:
\begin{equation}\label{PH_eiconale}
|\nabla u(\x)| = f(\x) \quad \textrm{for} \quad \x \in\Omega,
\end{equation}
where now
\begin{equation} \label{eq:def_f_PH_eiconale}
f(\x) = \displaystyle \sqrt{\frac{k_D^2-2I_+(\x)I_-(\x) + k_D^2 \sqrt{Q(\x)}}{2\biggl(I(\x)+k_S\biggr)^2}}, 
\end{equation}
with
\begin{eqnarray}
I_+(\x) &:=& I(\x) + k_S, \label{eq:def_Ip_PH_eiconale} \\
I_-(\x) &:=& I(\x) - k_S, \label{eq:def_Im_PH_eiconale} \\
Q(\x) &:=& k_D^2 + 8 k_S^2 + 8 \, I(\x) \, k_S. \label{eq:def_Q_PH_eiconale}
\end{eqnarray}

{\em Remark on the control interpretation.} 
The above analysis has shown that all the cases corresponding to the models proposed by Oren-Nayar 
and by Phong lead to a stationary Hamilton-Jacobi equation of the same form, namely
\[
\mu v({\bf x})=\min\limits_{a\in \partial B_3} \{ {\bf b}({\bf x},a)\cdot \nabla v({\bf x})+ f({\bf x},z,a,v({\bf x}))\}
\]
where the vector field ${\bf b}$ and the cost $f$ can vary according to the model and to the case. This gives to these models a control theoretical interpretation which can be seen as a generalization of the control interpretation for the original Lambertian model (which was related to the minimum time problem). 
In this framework, $v$ is the value function of a rescaled (by the Kru\v zkov change of variable) control problem in which one wants to drive the controlled system governed by 
\[ \dot y(t)=  {\bf b}({\bf y}(t),a(t))\]
($a(\cdot)$ here is the control function taking values in $\partial B_3$) from the initial position ${\bf x}$ to the target (the silhouette of the object) minimizing the cost associated to the trajectory.
The running cost associated to the position and the choice of the control will be given by $f$. More informations on this interpretation, which is not crucial to understand 
the application to the SfS problem presented in this paper, can be found in \cite{FF14}. 
\section{Semi-Lagrangian Approximation} 	\label{sec:general_convergence_theorem}
Now, let us state a general convergence theorem suitable for the class of differential operators appearing in the models described in the previous sections.
As we noticed, the unified approach presented in this paper has the big advantage to give a unique formulation for the three models in the form of a fixed point problem 
\begin{equation}
\mu v(\x) = T^M(\x,v,\nabla v), \quad \hbox{for } \x\in \Omega,
\end{equation}
where $M$ indicates the model, i.e. $M=L, ON, PH$. 

We will see  that the discrete operators of the ON--model and the PH--model described in the previous sections satisfy the properties listed here.
In order to obtain the fully discrete approximation we will adopt the semi-Lagrangian approach described in the book by Falcone and Ferretti \cite{FF14}. The reader can also refer to \cite{CFF13} for a similar approach to various image processing problems (including nonlinear diffusions models and segmentation).

Let $W_i=w(x_i)$ so that $W$ will be the vector solution giving the approximation of the height $u$ at every node $x_i$ of the grid. Note that in one dimension the index $i$ is an integer number, in two dimensions $i$ denotes a multi-index, $i=(i_1,i_2)$. 
We consider a semi-Lagrangian scheme written in a fixed point form, so we will write the fully discrete scheme as
\begin{equation} \label{eq:generic_operator}
W_i=\widehat T_i(W).
\end{equation}
Denoting by $G$ the global number of nodes in the grid, the operator corresponding to the oblique light source is $\widehat T:\R^G\rightarrow \R^G$ that is defined componentwise by
\begin{eqnarray} \label{eq:generic_T}
 \widehat T_i(W) \hspace{-0.5mm} := 
 \min\limits_{a\in\partial B_3} \{ e^{-\mu h} I[W](x_i^+) - \tau F(x_i,z,a) \} + \tau 
\end{eqnarray}
where $I[W]$ represents an interpolation operator based on the values at the grid nodes and 
\begin{eqnarray}
&& x_i^+:= x_i + h b(x_i,a) \\
&&\tau := (1-e^{-\mu \, h})/{\mu} \\
&&F(x_i,z,a) := P(x_i,z) a_3 (1-\mu W_i) \\
&&P: \Omega\times\R \rightarrow \R \hbox{ is continuous and nonnegative}.
\end{eqnarray}
Since  $w(x_i + h b(x_i,a))$ is approximated via $I[W]$ by interpolation on $W$ (which is defined on the grid $G$), it is important to use a monotone interpolation in order to preserve the properties of the continuous operator $T$  in the discretization. To this end, the typical choice is to apply a piecewise linear (or bilinear) interpolation operator $I_1 [W]: \Omega\rightarrow \R$   which allows to define a function defined for every $x\in \Omega$ (and not only on the nodes)
\begin{equation}\label{int1}
w(x)= I_1[W](x)=\sum_j \lambda_{ij} (a) W_j
\end{equation}
where
\begin{equation}\label{int2}
\sum_j \lambda_{ij} (a)=1 \quad \hbox{for} \quad x=\sum_j \lambda_{ij}(a) x_j.
\end{equation}
A simple explanation for \eqref{int1}-\eqref{int2} is that the coefficients $\lambda_{ij}(a)$ represent the local coordinates of the point $x$ with respect to the grid nodes (see \cite{FF14} for more details and other choices of interpolation operators). Clearly, in \eqref{eq:generic_T} we will apply the interpolation operator to the point $x^+_i=x_i + h b(x_i,a)$ and we will denote by $w$ the function defined by $I_1[W]$.

\noindent Comparing \eqref{eq:generic_T} with its analogue for the vertical light case we can immediately note that the former has the additional term $\tau F(x_i,z,a)$ which requires analysis.

\begin{theorem} \label{th:genprop}
Let $\widehat T_i(W)$ the $i$-th component of the operator defined as in \eqref{eq:generic_T}. Then, the following properties hold true:
\begin{enumerate}
\item Let \\
$\overline{a}_3 \equiv arg \min\limits_{a\in\partial B_3} \{ e^{-\mu h} w(x_i + h b(x_i,a)) -  \tau F(x_i,z,a) \}$ and assume
\begin{equation}\label{eq:hp_condition1}
P(x_i,z) \overline{a}_3 \leq 1.
\end{equation}
Then $0\leq W\leq \frac{1}{\mu}$ implies $0\leq \widehat T(W)\leq \frac{1}{\mu}$.\\
\item $\widehat T$ is a monotone operator, i.e., $W\leq \overline W$ implies $\widehat T(W)\leq \widehat T(\overline W)$. \\
\item $\widehat T$ is a contraction mapping in $L^{\infty}([0,1/\mu)^G)$ if 

$P(x_i,z) \,\overline{a}_3 < \mu$.
\end{enumerate}
\end{theorem}

{\em Proof.}
\begin{enumerate}
\item To prove that $W\leq \frac{1}{\mu}$ implies $T(W)\leq \frac{1}{\mu}$ we just note that
\begin{equation}
\widehat T(W) \leq \frac{e^{-\mu h}}{\mu} + \tau = \frac{1}{\mu}.
\end{equation}
Let $W \geq 0$; then
\begin{equation}
\begin{aligned}
\widehat T(W) & \geq - \tau P(x_i,z) \,\overline{a}_3 (1-\mu W_i) + \tau \\
& = \tau \left( 1 - P(x_i,z) \,\overline{a}_3 (1-\mu W_i) \right).
\end{aligned}
\end{equation}
This implies that $\widehat T(W) \geq 0$ if $ P(x_i,z) \,\overline{a}_3 \leq 1$ since $0 \leq 1-\mu W_i\leq 1$.\\
\item In order to prove that $\widehat T$ is monotone, let us observe first that for each couple of functions $w_1$ and $w_2$ such that $w_1(x)\leq w_2(x)$ for every $x\in \Omega$  implies
\begin{eqnarray}
&&e^{-\mu h} [ w_1(x + h b(x,a^*)) - w_2(x + h b(x,\overline{a})) ]  \vspace{0.1cm} \\
&&\quad - \tau P(x,z) ( a^*_3 (1-\mu w_1(x)) - \overline{a}_3 (1-\mu w_2(x)) ) \nonumber \\
&&\leq e^{-\mu h} [ w_1(x + h b(x,\overline{a})) - w_2(x + h b(x,\overline{a})) ] \nonumber \\
&&\quad  + \tau P(x,z) \overline{a}_3 (w_1(x) - w_2(x))\nonumber
\end{eqnarray}
where $a^*$ and $\overline a$ are the two arguments corresponding to the minimum on $a$ of the expression
\begin{equation}
e^{-\mu h}  w(x + h b(x,a)) - \tau P(x,z)  a_3 (1-\mu w(x)) 
\end{equation}
respectively for $w=w_1,w_2$. Hence, if we take now two vectors $W$ and $\overline W$ such that $W\le \overline W$ and we denote respectively by $w$ and $\overline w$ the corresponding functions defined by interpolation, we will have $w(x)=I_1[W](x)\le I_1[\overline W](x)=\overline w(x)$, for every $x\in \Omega$,   because  the linear interpolation operator $I_1$ is monotone. Then, by \eqref{eq:generic_T}, setting $x_i^+=x_i + h b(x_i,\overline{a})$ we get 
\begin{eqnarray}\label{eq:cond2}
&&\widehat T_i(W) - \widehat T_i(\overline W) \leq e^{-\mu h} [ I_1[W](x_i^+)) - I_1[\overline W](x_i^+)) \nonumber]\\
&&  + \tau P(x_i,z) \overline{a}_3 (W_i - \overline W_i)\le 0
\end{eqnarray}
 where the last inequality follows from $P\ge 0$. So we can conclude that  $\widehat T(W)-\widehat T(\overline W)\le 0$.
Note that this property does not require condition \eqref{eq:hp_condition1} to be satisfied.\\
\item Let us consider now two vectors $W$ and $\overline W$ (dropping the condition $W\le \overline W$) and assuming 
\begin{equation}\label{eq:hp_condition3}
P(x_i,z) \overline{a}_3 < \mu. 
\end{equation}
To prove that $T$ is a contraction mapping note that following the same argument used to prove the second statement we can obtain  \eqref{eq:cond2}. Then, by  applying the definition of $I_1$, we get 
\begin{equation*}\label{eq:cond3}
\widehat T(W) -\widehat T(\overline W) \leq \left( e^{-\mu h} + \tau P(x_i,z) \overline{a}_3 \right) ||W- \overline W ||_{\infty}.
\end{equation*}
Reversing the role of $W$ and $\overline W$, one can also obtain
\begin{equation*}\label{eq:cond4}
\widehat T(\overline W) - \widehat T(W) \leq \left( e^{-\mu h} + \tau P(x_i,z) \overline{a}_3 \right) ||W- \overline W||_{\infty}
\end{equation*}
and conclude then that  $\widehat T$ is a contraction mapping in $L^\infty$ if and only if
\begin{equation}
e^{-\mu h} + \tau P(x_i,z) \overline{a}_3 < 1
\end{equation}
and this holds true if the bound \eqref{eq:hp_condition3} is satisfied.
\end{enumerate}
\begin{flushright}
    $\Box$
\end{flushright}

{\em Remark on the role of $\mu$.} The parameter $\mu$ can be tuned to satisfy the inequality which guarantees the contraction map property for $\widehat T$. 
This parameter adds a degree of freedom in the Kru\v zkov change of variable and modifies the slope of $v$. However, in the practical applications we have done in our tests, this parameter has been always set to 1 so in our experience this parameter does not seem to require a fine tuning. \\

\vspace{0.2cm}
{\em Remark on the choice of the interpolation operator.} Although $ I_1$ can be replaced by a high-order interpolation operator (e.g. a cubic local Lagrange interpolation), the monotonicity of the interpolation operator plays a crucial role in the proof because we have to guarantee that the extrema of interpolation polynomial stay bounded by the minimum and maximum of the values at the nodes. This property is not satisfied by quadratic or cubic local interpolation operators and the result is that this choice introduces spurious oscillations in the numerical approximation. A cure could be to adopt Essentially Non Oscillatory (ENO) interpolations. A detailed discussion on this point is contained in \cite{FF14}.

\vspace{0.2cm}
{\em Remark on the minimization.} In the definition of the fixed point operator $\widehat T$ there is a minimization over $a\in\partial B_3$.  A simple way to solve it is to build a discretization of   $\partial B_3$ based on a finite number of points and get the minimum by comparison. One way to do it is to discretize the unit sphere by spherical coordinates, even a small number of nodes will be sufficient to get convergence. A detailed discussion on other methods to solve the minimization problem  is contained in \cite{FF14}.  

Let us consider now the algorithm based on the fixed-point iteration
\begin{equation}\label{def:sucpf}
 \left\{ \begin{array}{ll}
W^n &=\widehat T(W^{n-1}),\\
W^0 & \textrm{given}.
\end{array} \right.
\end{equation}
We can state the following convergence result
\begin{theorem} \label{th:genconv}
Let $W^k$ be the sequence generated by \eqref{def:sucpf}. Then the following results hold:\\
1.  Let $W^0\in \mathcal{S}=\{W\in\R^G: W\ge \widehat T(W)\}$, then the $W^k$ converges monotonically decreasing to a fixed point $W^*$ of the $\widehat T$ operator;\\
2. Let us choose $\mu>0$ so that the condition $P(x_i,z) \overline{a}_3 \leq \mu$ is satisfied. Then, $W^k$ converges to the unique fixed point $W^*$ of the $\widehat T$. Moreover, if $W^0\in \mathcal{S}$ the convergence is monotone decreasing.
\end{theorem}
{\em Proof.}\\
\begin{enumerate}
\item  Starting from a point in the set of super-solutions $\mathcal{S}$, the sequence is non increasing and lives in $\mathcal{S}$ which is a closed set bounded from below (by 0). Then, $W^k$ converges and the limit is necessarily a fixed point for $\widehat T$. 
\item The assumptions guaranteed by Theorem \ref{th:genprop} are satisfied and $\widehat T$ is a contraction mapping in $[0,1/\mu]$, so the fixed point is unique. The monotonicity of $\widehat T$ implies that starting from $W^0\in \mathcal{S}$ the convergence is monotone decreasing.
\end{enumerate}
\begin{flushright}
    $\Box$
\end{flushright}
It is important to note that the change of variable allows for an easy choice of the initial guess $W^0\in \mathcal{S}$ for which we have the natural choice $W^0=(1/\mu, 1/\mu, \dots, 1/\mu)$ and monotonicity can be rather useful to accelerate convergence as shown  in \cite{Falcone97}. 
A different way to improve convergence is to apply Fast Sweeping or Fast Marching methods as illustrated in \cite{Tozza14,ToFal_Dagstuhl}. A crucial role is played by boundary conditions on the boundary of $\Omega$, where usually we impose the homogeneous Dirichlet boundary condition, $v=0$. This condition implies that the shadows must not cross the boundary of $\Omega$, so the choice $\omega_3 = 0$ corresponding to an infinite shadow behind the surface is not admissible. However, other choices are possible: to impose the height of the surface on $\partial \Omega$ we can set $v=g$ or to use a more neutral boundary condition we can impose $v=1$ (state constraint boundary condition). More informations on the use of boundary conditions for these type of problems can be found in \cite{FF14}.   

\subsection{Properties of the discrete operators $\widehat T^{ON}$ and $\widehat T^{PH}$} \label{sec:propertie}
We consider a semi-Lagrangian (SL) discretization of \eqref{ON_pde_v} written in a fixed point form, so we will write the {\em SL fully discrete scheme for the ON--model} as
\begin{equation} \label{eq:generic_operator_ON}
W_i=\widehat T^{ON}_i(W),
\end{equation}
where $ON$ is the acronym identifying the ON--model. 
Using the same notations of the previous section, the operator corresponding to the oblique light source is $\widehat T^{ON}:\R^G\rightarrow \R^G$ that with linear interpolation can be written as
\begin{equation*} \label{eq:T_ON}
\widehat T^{ON}_i(W) := \min\limits_{a\in\partial B_3} \{ e^{-\mu h} I_1[W](x_i^+) - \tau F^{ON}(x_i,z,a) \} + \tau,
\end{equation*}
where 
\begin{eqnarray*} \label{components_T_ON}
&&\tau :=  \displaystyle \frac{1-e^{-\mu \, h}}{\mu} \vspace{0.2cm} \\
&&b^{ON}(x_{i},a) :=  \displaystyle  \frac{1}{A\omega_3} \left( c(x_{i},z) a_1 - A\omega_1, c(x_{i},z) a_2 - A\omega_2 \right) \vspace{0.2cm} \\
&&c(x_{i},z) := \displaystyle I(x_{i}) - B +B (d(x_{i},z) \cdot \vomega)^2  \vspace{0.2cm} \\
&&d(x_{i}, z) :=  \displaystyle  \nabla S(x_{i},z) / |\nabla S(x_{i},z)|  \vspace{0.2cm} \\ 
&&F^{ON}(x_i, z,a) := P^{ON}(x_{i},z) a_3 (1-\mu W_i))\vspace{0.2cm} \\
&&P^{ON}(x_{i},z) :=  \displaystyle \frac{ c(x_{i},z)}{A\omega_3}. 
\end{eqnarray*}
Note that, in general, $P^{ON}$ will not be positive but that condition can be obtained tuning the parameter $\sigma$ since the coefficients $A$ and $B$ depend on $\sigma$. This explains why in some tests we will not be able to get convergence for every value of $\sigma\in [0, \pi/2)$. 
Once the non negativity of $P^{ON}$ is guaranteed, we can follow the same arguments of Theorem \ref{th:genprop}  to check that the  discrete operator $\widehat T^{ON}$ satisfies the three properties which are necessary to guarantee convergence as in Theorem \ref{th:genconv} provided we set $P=P^{ON}$ in that statement.

For the Phong model, the semi-Lagrangian discretization of \eqref{PH_pde_v} written in a fixed point form gives 
\begin{equation} \label{eq:generic_operator_PH}
W_i=\widehat  T^{PH}_i(W),
\end{equation}
where $\widehat T^{PH}:\R^G\rightarrow \R^G$, that is defined componentwise by
\begin{equation*} \label{eq:T_PH}
\widehat T^{PH}_i(W) := \min\limits_{a\in\partial B_3} \{ e^{-\mu h} I_1[W](x_i^+)  - \tau F^{PH}(x_i,z,a) \} + \tau,
\end{equation*}
where, in the case of oblique light source and vertical position of the observer,
\begin{eqnarray*} \label{components_T_PH}
&&\tau :=  \displaystyle \frac{1-e^{-\mu \, h}}{\mu} \vspace{0.2cm} \\
&&b^{PH}(x_{i},a) :=  \displaystyle  \frac{\left( c(x_{i}) a_1 - k_D \omega_1, c(x_{i}) a_2 - k_D \omega_2 \right)}{Q^{PH}(x_{i},z)}  \vspace{0.2cm} \\
&&c(x_{i}) := \displaystyle I(x_{i}) + \omega_3 k_S \vspace{0.2cm} \\
&&d(x_{i}, z) :=  \displaystyle  \nabla S(x_{i},z) / |\nabla S(x_{i},z)| \vspace{0.2cm} \\ 
&&Q^{PH}(x_{i},z) := 2 k_S (d(x_{i},z) \cdot \vomega) + k_D \omega_3 \vspace{0.2cm} \\ 
&&F^{PH}(x_i,z,a) := P^{PH}(x_{i},z) a_3 (1-\mu W_i))\vspace{0.2cm} \\
&&P^{PH}(x_{i},z) :=  \displaystyle \frac{ c(x_{i})}{Q^{PH}(x_{i}, z)}. 
\end{eqnarray*}
Here the model has less parameters and $P^{PH}$ will always be nonnegative. Again, following the same arguments of Theorem \ref{th:genprop}, we can check that the  discrete operator $\widehat T^{PH}$ satisfies the three properties which are necessary to guarantee convergence as in Theorem \ref{th:genconv} provided we set $P=P^{PH}$ in that statement.

\section{Numerical Simulations}\label{sec:numerical_tests}
In this section we show some numerical experiments on synthetic and real images in order to analyze the behavior of the parameters involved in the ON--model and the PH--model and to compare the performances of these models with respect to the classical L--model and with other numerical methods too. 
All the numerical tests in this section have been implemented in language C++. The computer used for the simulations is a MacBook Pro 13" Intel Core 2 Duo with speed of 2.66 GHz and 4GB of RAM (so the CPU times in the tables refer to this specific architecture).

We denote by $\mathcal{G}$ the discrete grid in the plane getting back to the double index notation $x_{ij}$, $G:=card(\mathcal{G})=n\times m$. We define $G_{in} := \{x_{ij} : x_{ij} \in \Omega \}$ as the set of grid points inside $\Omega$; $G_{out} := G\setminus G_{in}$. The boundary $\partial \Omega$ will be then approximated by the nodes such that at least one of the neighboring points belongs to $G_{in}$. 
For each image we define a map, called {\em  mask}, representing the pixels $x_{ij} \in G_{in}$ in white and the pixels
$x_{ij} \in G_{out}$ in black. In this way it is easy to distinguish the nodes that we have to use for the reconstruction (the nodes inside $\Omega$) and the nodes on the boundary $\partial \Omega$ (see e.g. Fig. \ref{fig:mask_sphere}).  

Regarding the minimization over $a\in\partial B_3$ that appears in the definition of the fixed point operators associated to the models, in all the tests we discretize the unit sphere by spherical coordinates, considering 12 steps in $\theta$ and 8 in $\phi$, where $\theta$ is the zenith angle and $\phi$ is the azimuth angle. 

\subsection{Synthetic tests}
If not otherwise specified, all the synthetic images are defined on the same rectangular domain containing the support of the image, $\Omega  \equiv [-1,1]\times[-1,1]$.
We can easily modify the number of the pixels choosing different values for the steps in space $\Delta\,x$ and $\Delta\,y$. 
The size used for the synthetic images is  $256\times 256$ pixels,  unless otherwise specified. $X$ and $Y$ represent the real size (e.g. for $\Omega  \equiv [-1,1]\times[-1,1]$, $X=2, Y=2$).
For all the synthetic tests, since we know the algebraic expression of the surfaces, the input image rendering in gray-levels is obtained using the corresponding reflectance model. This means that for each model and each value of parameter involved in it, the reconstruction will start from a different input image. Clearly, this is not possible for real images, so for these tests the input image will be always the same for all the models, independently of the values of the parameters.
Moreover, we fix $\mu = 1$ and we choose the value of the tolerance $\eta$ for the iterative process equal to $10^{-8}$ for the tests on synthetic images, using as stopping rule $|W^{k+1} - W^k|_{max} \leq \eta$, where $k+1$ denotes the current iteration. We will see that dealing with real images, sometimes  we will need to increase $\eta$. 

\paragraph{Test 1: Sphere.}
For this first test we will use the semisphere  defined as 
\begin{equation}
		\begin{cases}
			u(x,y) = \sqrt{r^2-x^2-y^2} & (x,y) \in G_{in},\\
			u(x,y) = 0 &  (x,y) \in G_{out},
		\end{cases}
\end{equation}
where 
\begin{eqnarray}
&& G_{in} := \{(x,y): x^2+y^2\le r^2\},\\
&& r := \frac{\min\{X,Y\}}{2}+{2\,\tilde{\delta}},
\end{eqnarray}
and 
$\tilde{\delta} := \max\lbrace \Delta{}x,\Delta{}y\rbrace$.

\begin{figure}[h!]
\centering
 \subfigure[Sphere Input]
   {\includegraphics[width=0.23\textwidth]{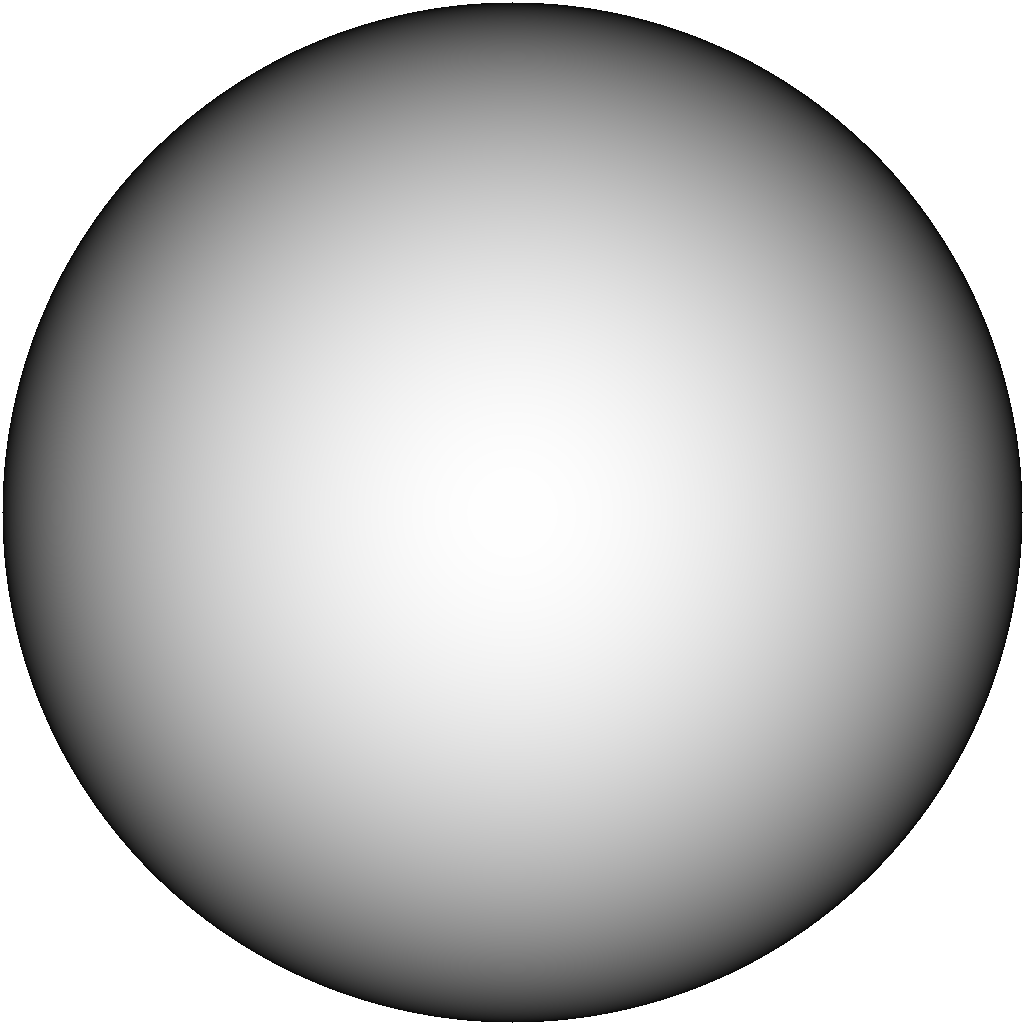} \label{fig:input_sphere}  }
 \hspace{2mm}
  \subfigure[Sphere Mask]
  {\includegraphics[width=0.23\textwidth]{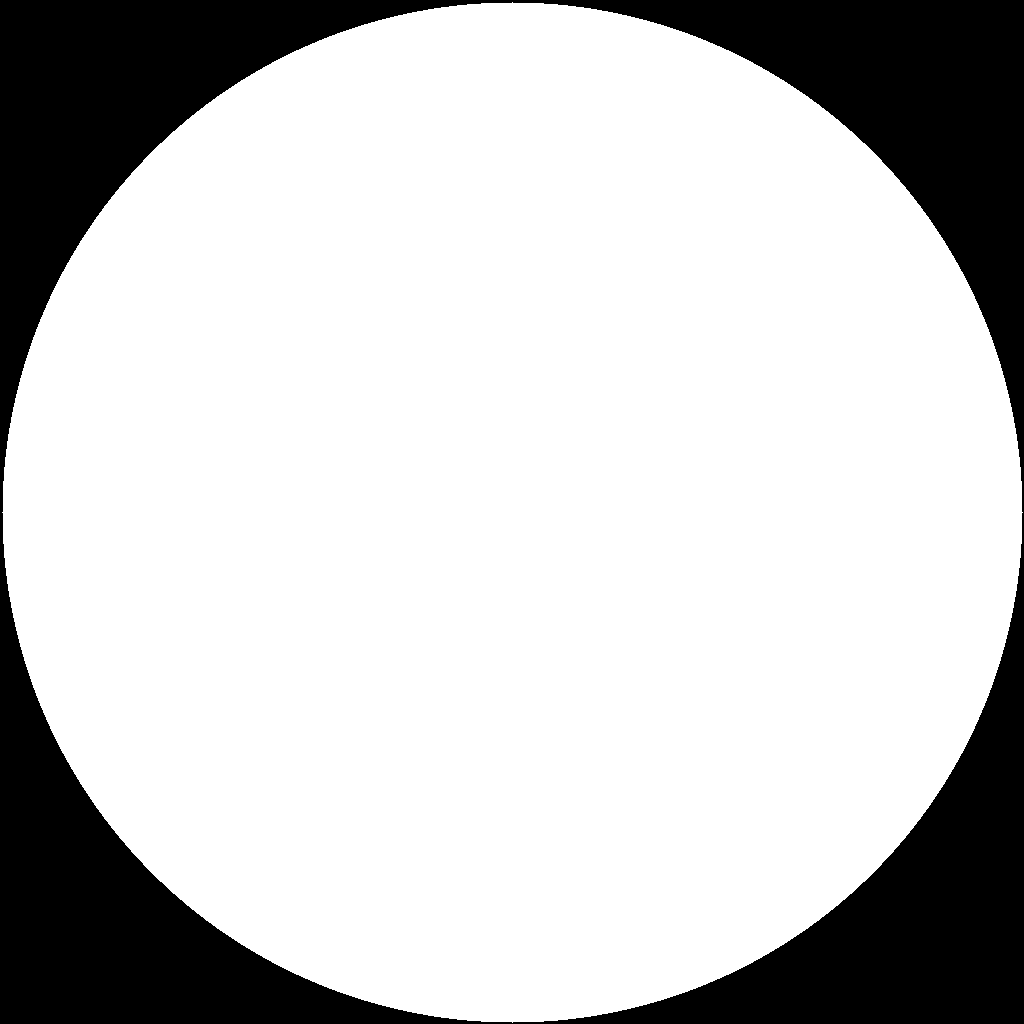} \label{fig:mask_sphere}  }
 \hspace{2mm}
 \subfigure[Sphere surface]
   {\includegraphics[width=0.25\textwidth]{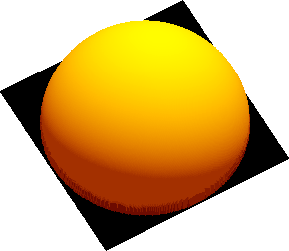}\label{fig:surface_sphere}} 
 \caption{Sphere via the L--model:
            (a)~Input image; 
            (b)~Mask;
            (c)~surface.}
 \label{fig:sphere_images}
 \end{figure}

As example, we can see in Fig. \ref{fig:sphere_images} the input image, the corresponding mask and the surface reconstructed by the L--model. 
\noindent The values of the parameters used in the simulations are indicated in Table \ref{tab:sphere_parameter}. Note  (in Table \ref{tab:sphere_errors}) that  when the specular component is zero for the PH--model, we just have the contribution of the diffuse component so we have exactly the same error values of the L--model, as expected. 
By increasing the value of the coefficient $k_S$ and, as a consequence, decreasing the value of $k_D$, in the PH--model the $L^2(I)$ and $L^\infty(I)$ errors on the image grow albeit slightly and still remain of the same order of magnitude, whereas the errors on the surface decrease.
For the ON--model, the same phenomenon appears when we set the roughness parameter $\sigma$ to zero: we bring back to the L--model and, hence, we obtain the same errors on the image and the surface.  Note that the errors in $L^2(I)$ and $L^\infty(I)$ norm for the image and the surface, decrease by increasing the value of $\sigma$. 
This seem to imply that the model and the approximation work better for increasing roughness values.

We point out that the errors computed on the images $I$ in the different norms are over integer between the input image and the image computed a posteriori using the value of $u$ just obtained by the methods. This is why small errors becomes bigger since can jump from an integer to the other one, as for the Phong case.

\begin{table}[h!]
 \caption{Sphere: parameter values used in the models.}  \label{tab:sphere_parameter} 
   \begin{tabular}{p{1.3cm}p{1.1cm}p{1.1cm}p{1.1cm}p{1.1cm}}
\hline\noalign{\smallskip}
Model & \hspace{0.cm}$\sigma$ & \hspace{0.cm}$k_D$ & \hspace{0.cm}$k_S$ & \hspace{0.cm}$\alpha$   \\
\hline 
LAM & 	\hspace{0.cm} & & &  \\
ON-00 &		\hspace{0.cm}0 &\hspace{0.cm} &\hspace{0.cm} & \hspace{0.cm} \\
ON-04 & 	 \hspace{0.cm}0.4 & \hspace{0.cm} & \hspace{0.cm} & \hspace{0.cm}   \\
ON-08 & 	 \hspace{0.cm}0.8 &  \hspace{0.cm} & \hspace{0.cm} &\hspace{0.cm}   \\
PH-s00 &		\hspace{0.cm} &\hspace{0.cm} 1&\hspace{0.cm} 0& \hspace{0.cm} 1 \\
PH-s04 & 	 \hspace{0.cm} &  \hspace{0.cm} 0.6 & \hspace{0.cm} 0.4 &\hspace{0.cm} 1   \\
PH-s08 & 	\hspace{0.cm} &\hspace{0.cm} 0.2 &\hspace{0.cm} 0.8 & \hspace{0.cm} 1 \\
  \hline 
\end{tabular}
\end{table}

In Table \ref{tab:sphere_iter&CPU} we reported the number of iterations and the CPU time (in seconds) referred to the three models with the parameter indicated  in Table \ref{tab:sphere_parameter}. For all the models, also varying the parameters involved, the number of iterations is always about $2000$ e the CPU time slightly greater than $2$ seconds, so the computation is really fast.

\begin{table} [h!]
     \caption{Synthetic sphere: iterations and CPU time in seconds for the models with  vertical light source $\vomega = (0,0,1)$.}   \label{tab:sphere_iter&CPU} 
  \begin{tabular}{p{2.5cm}p{1.5cm}p{1.5cm}}
\hline\noalign{\smallskip}
  SL--Schemes & Iter.  & $[sec.]$ \\
  \hline 
LAM  & 2001 &  2.14 \\ 
ON-00  & 2001 &  2.06 \\  
ON-04  & 2020 &  2.24 \\  
ON-08  & 2016 &  2.24 \\  
PH-s00 & 2001 &  2.11 \\  
PH-s04  & 2008 &  2.45 \\  
PH-s08  & 2056 &  2.27 \\  
  \hline 
\end{tabular}
\end{table}

\begin{table} [h!]
     \caption{Synthetic sphere:  $L^2$, $L^{\infty}$ errors with vertical light source $\vomega = (0,0,1)$.}   \label{tab:sphere_errors} 
\begin{tabular}{p{2.1cm}p{1cm}p{1cm}p{1cm}p{1cm}}
\hline\noalign{\smallskip}
  SL--Schemes & $L^2(I)$ & $L^\infty(I)$ & $L^2(S)$ & $L^\infty(S)$  \\
  \hline 
LAM  &  0.0046 &  0.0431 &  0.0529 &  0.0910  \\
ON-00  &  0.0046 &  0.0431 &  0.0529 &  0.0910  \\
ON-04  &   0.0039 &  0.0353 &  0.0513 &  0.0882  \\
ON-08  &  0.0035 &  0.0314 &  0.0506 &  0.0881   \\
PH-s00 &  0.0046 &  0.0431 &  0.0529 &  0.0910   \\
PH-s04  &  0.0064 &  0.0471 &  0.0511 &  0.0896 \\
PH-s08  &  0.0090 &  0.0706 &  0.0386 &  0.0752   \\
  \hline 
\end{tabular}
\end{table}

\begin{figure}[h!]
\centering
 \subfigure[Tent Input]
   {\includegraphics[width=0.31\textwidth]{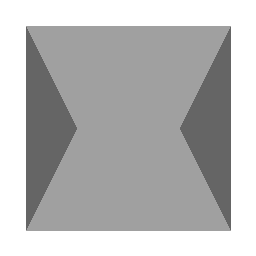} \label{fig:input_tent}  }
 \hspace{2mm}
 \subfigure[Tent surface]
   {\includegraphics[width=0.33\textwidth]{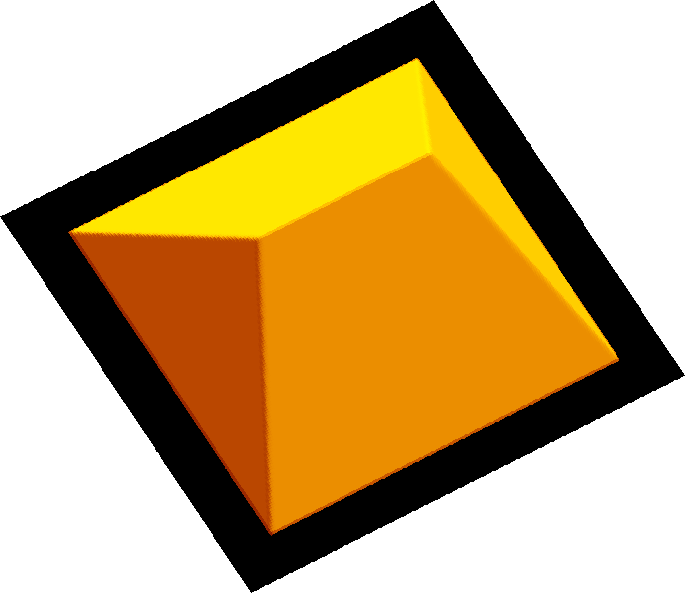}\label{fig:surface_tent}} 
 \caption{Tent via the ON--model with $\sigma = 0.3$:
            (a)~Input image; 
            (b)~3D reconstruction.}
 \label{fig:tent_images}
 \end{figure}

\paragraph{Test 2: Ridge tent.}
In tests on synthetic images, the relevance of the choice of a model depends on which model was used to compute the images. 
In the previous test, the parameters that are used for the 3D reconstruction are identical to those  used to compute the synthetic sphere input images, so there is a perfect match. However,  for real applications, it is relevant to examine the influence of an error in the  parameter values. To this end  we can produce an input image with the Oren-Nayar model using $\sigma = 0.1$ and  then process this image with the same model using a different value of $\sigma$ to see how the results are affected by this error. This is what we are going to do for the ridge tent. 
Let us consider the ridge tent defined by the following equation
	\begin{equation}\label{eqn:tendacanadese}
	\begin{cases}
		u(x,y)=\min\left\lbrace 
			-2\, |x|+\frac{4}{5}\,X, 
			- |y|+\frac{2}{5}\,Y \right \rbrace &\hspace{-0.2cm}\\
			&	\hspace{-2cm}	(x,y) \in G_{in}, \\
		u(x,y)=0 &\hspace{-2cm}
				(x,y) \in G_{out},
	\end{cases}
	\end{equation}
	where
	\[
	G_{in}:=\left\{(x,y):  \frac{x}{X},\frac{y}{Y}  
            <  \frac{2}{5}  \right\}.
	\]
In Fig. \ref{fig:tent_images} we can see an example of reconstruction obtained by using the ON--model with $\sigma = 0.3$, under a vertical light source $\vomega = (0,0,1)$.
A first remark is that the surface reconstruction is good even if in this case the surface is not differentiable. Moreover, note that there are no oscillations near the kinks 
where there is jump in the gradient direction.
Let us examine the stability with respect to the parameters. We have produced seven input images for the ridge tent, all of size $256\times{}256$, with the following combinations of models and parameters:
        \begin{description}
           \item{LAM} Lambertian model;
           \item{ON1} Oren-Nayar model with $\sigma=0.1$;
           \item{ON3} Oren-Nayar model with $\sigma=0.3$;
           \item{ON5} Oren-Nayar model with $\sigma=0.5$;
           \item{PH1} Phong model with $\alpha=1$ and $k_S=0.1$;
           \item{PH3} Phong model with $\alpha=1$ and $k_S=0.3$;
           \item{PH5} Phong model with $\alpha=1$ and $k_S=0.5$.
        \end{description}
Then we have computed the surfaces corresponding to all the parameter choices (i.e. matching and not matching the first choice).
The results obtained in this way have been compared in terms of $L^2$ and $L^\infty$ norm errors with respect to the original surface. The errors obtained by the ON--model are shown in Table~\ref{tab:tent_ON} and in Table \ref{tab:tent_PH} for the PH--model.
  
\begin{table}[h!]
\caption{Ridge tent, ON--model: $L^2,L^\infty$ errors for the surface. In each column the model used to produce the input image, in the row the model used for the $3D$ reconstruction. }
\label{tab:tent_ON}
\begin{tabular}{lllll}
\hline\noalign{\smallskip}
   $L^2$ & LAM     &ON1      &ON3     &ON5     \\
\noalign{\smallskip}\hline\noalign{\smallskip}
  LAM & 0.0067 & 0.0172	& 0.0933	& 0.1920	\\
 ON1 & 0.0082	& 0.0068	& 0.0821	& 0.1801 \\
 ON3 & 0.0832 	& 0.0700	& 0.0086	&  0.1033	\\
 ON5 & 0.1946	& 0.1784	& 0.0923 	& 0.0067	\\
\noalign{\smallskip}\hline
\hline\noalign{\smallskip}
   $L^\infty$ & LAM     &ON1      &ON3     &ON5     \\
\noalign{\smallskip}\hline\noalign{\smallskip}
 LAM	  	& 0.0094	& 0.0315 & 0.1942   & 0.4060	\\
 ON1 	& 0.0199  & 0.0093 & 0.1701   & 0.3805	\\
 ON3 	& 0.1784 & 0.1507  & 0.0118   &  0.2156 \\
 ON5 	& 0.4104 &  0.3769 &  0.1976  &  0.0094	\\
\noalign{\smallskip}\hline
\end{tabular}
\end{table}
\begin{table}[h!]
\caption{Ridge tent, PH-model: $L^2,L^\infty$ errors for the surface. In each column the model used to produce the input image, on the row the model used for the $3D$ reconstruction. }
\label{tab:tent_PH}
\begin{tabular}{lllll}
\hline\noalign{\smallskip}
   $L^2$ & LAM  &PH1     &PH3      & PH5    \\ 
\noalign{\smallskip}\hline\noalign{\smallskip}
 LAM &  0.0067 & 0.0841 &  0.2867  &  0.6146 \\
 PH1 &  0.0586 &  0.0067 &  0.1996 &  0.5031 \\
 PH3 &  0.1403 &  0.0955 &  0.0073 &  0.2687 \\
 PH5 &  0.1915 &  0.1664 &  0.0976 &  0.0060 \\
\noalign{\smallskip}\hline
\hline\noalign{\smallskip}
   $L^\infty$ & LAM   &PH1     &PH3      & PH5   \\ 
\noalign{\smallskip}\hline\noalign{\smallskip}
 LAM &  0.0094 &  0.1740 &  0.6068 &  1.3123 \\
 PH1 &  0.1243  &  0.0108 &  0.4202 &  1.0741 \\
 PH3 &  0.3167 &  0.2245 &  0.0149 &  0.5718   \\
 PH5 &  0.4503  &  0.4241 &  0.2907  &  0.0093  \\
\noalign{\smallskip}\hline
\end{tabular}
\end{table}
 Analyzing the errors in Table \ref{tab:tent_ON} and Table \ref{tab:tent_PH}, we can observe that using the same model to generate the input image and to reconstruct the surface is clearly the optimal choice. The errors on the surface grow more as we consider a parameter $\sigma$ other than the one used to generate the input image as data for the 3D reconstruction. For the ON--model we loose one or two order of magnitude, depending on the ``distance" of the parameter from the source model. For the PH--model we can observe that the $L^2$ and $L^\infty$ errors grow more as we consider a different $k_S$ for the generation of the image and for the reconstruction, loosing one or two order of magnitude. However,  the two models seems to be rather stable with respect to a variation of the parameters since the errors do not increase dramatically varying the parameters.

\paragraph{Test 3: Concave/convex ambiguity for the ON--model.}
We consider this test in order to show that the ON--model is not able to overcome the concave/convex ambiguity typical of the SfS problem although it is a model more realistic than the classical L--model.
Let consider the following function
\begin{equation}\label{example_conc_conv_ON}
u(x,y) =
\left\{ \begin{array}{ll}
  -(1-(x^2-y^2))^2+1, & \\
  & \hspace{-1.2cm} \hbox{if } (x^2 + y^2)<2, \\
 0 &\hspace{-1.2cm} \hbox{ otherwise}.
\end{array}  \right.
\end{equation}

\begin{figure}[h!]
\centering
 \subfigure[]
   {\includegraphics[width=0.28\textwidth]{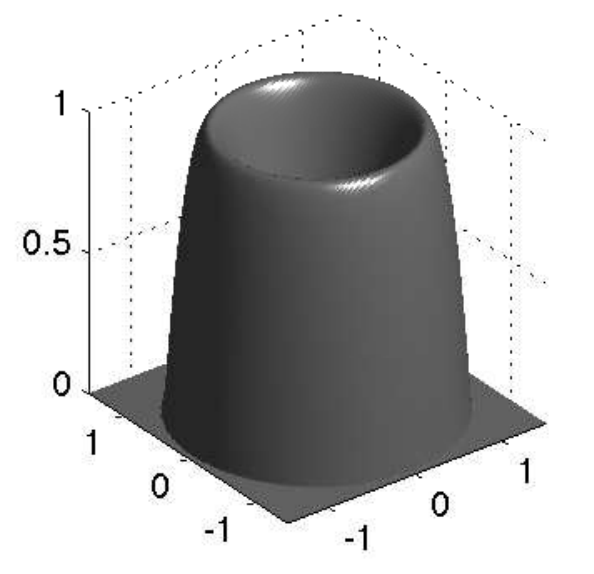} \label{fig:original_surface}  }
 \hspace{1mm}
  \subfigure[]
  {\includegraphics[width=0.3\textwidth]{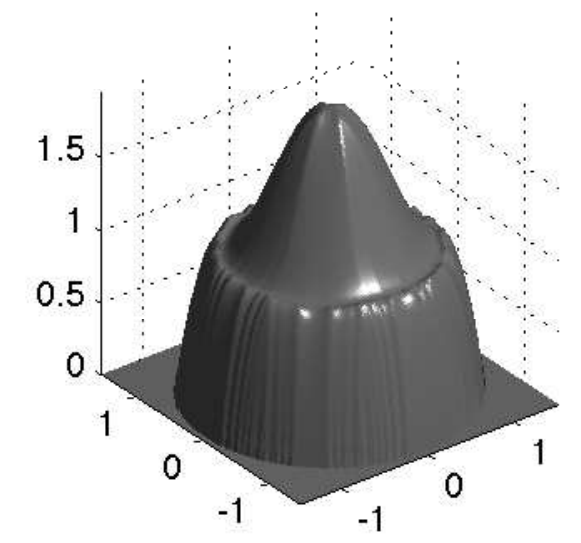} \label{fig:maximal_sol}  }
 \hspace{1mm}
 \subfigure[]
   {\includegraphics[width=0.3\textwidth]{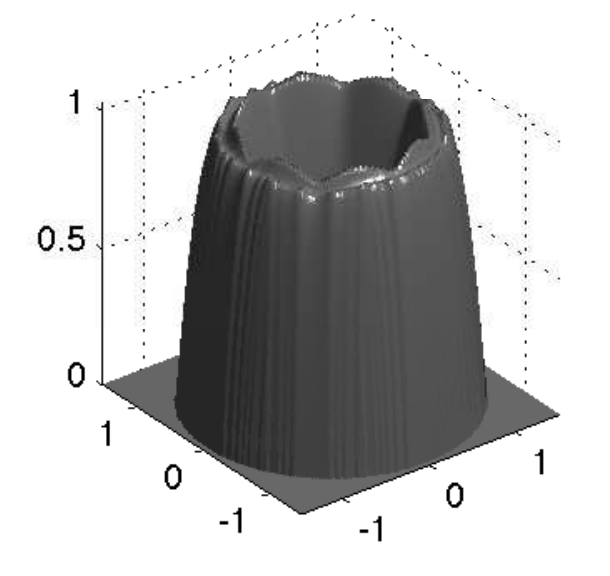}\label{fig:approx_sol}} 
 \caption{Example of concave/convex ambiguity for the ON--model with $\sigma = 0.5$ and $\vomega = (0,0,1)$:
            (a)~Original Surface; 
            (b)~Maximal Solution;
            (c)~Approximated surface with value in the origin equal to zero.}
 \label{fig:conc_conv_ON_images}
 \end{figure}
 \noindent We discretize the domain $\Omega = [-1.5,1.5]\times[-1.5,1.5]$ with $151\times151$ nodes. The fixed point has been computed with an accuracy of $\eta = 10^{-4}$ and the stopping rule for the algorithm based on the fixed point iteration defined in  \eqref{def:sucpf} is $|W^{k+1} - W^k|_{max} \leq \eta$, where $k+1$ denotes the current iteration. The iterative process starts with $W^0 = 0$ on the boundary and $W^0=1$ inside in order to proceed from the boundary to the internal constructing a monotone sequence
 (see \cite{BH85,FSS03} for details on the approximation of maximal solutions).

\noindent Looking at Fig. \ref{fig:conc_conv_ON_images} we can note that the scheme chooses the maximal viscosity solution stopping after 105 iterations, which does not coincide with the original surface. In order to obtain a reconstruction closer to the original surface, we fix the value in the origin at zero. In this way we forced the scheme to converge to a solution different from the maximal solution (see Fig. \ref{fig:approx_sol}). We obtained this different solution shown in Fig. \ref{fig:approx_sol} after 82 iterations. 

\paragraph{Test 4: Concave/convex ambiguity for the PH--model.}
The fourth synthetic numerical experiment is related to the sinusoidal function defined as follow
\begin{equation}
\begin{cases}
u(x,y)=\displaystyle 0.5+0.5 \, \sin(\frac{\pi\,x}{\Delta x}) \, \sin(\frac{\pi\,y}{\Delta y}), & \\
&\hspace{-1.2cm} (x,y) \in G_{in},\\
u(x,y)=0, &\hspace{-1.2cm} (x,y) \in G_{out}.
\end{cases}
\label{eq:sinusoidal}
\end{equation}
With this test we want to show that also the PH--model is not able to overcome the concave/convex ambiguity typical of the SfS problem. 

\begin{figure}[h!]
\centering
   \begin{tabular}{p{1.8cm}p{1.8cm}p{2.9cm}}
   \hline\noalign{\smallskip}
  \hspace{0.6cm} \textbf{in}  &  \hspace{0.4cm} \textbf{out} &  \hspace{0.cm} \textbf{reconstruction} \\  
  \hline \\
\includegraphics[width=1.8cm]{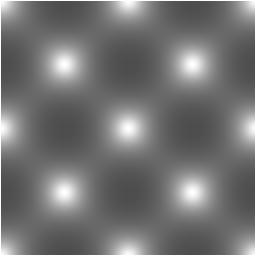} &\hspace{-0.2cm}
\includegraphics[width=1.8cm]{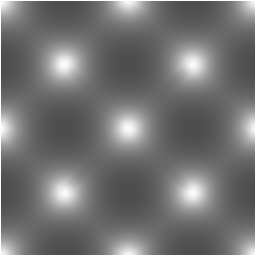} & \hspace{0.3cm}
\includegraphics[width=2.2cm]{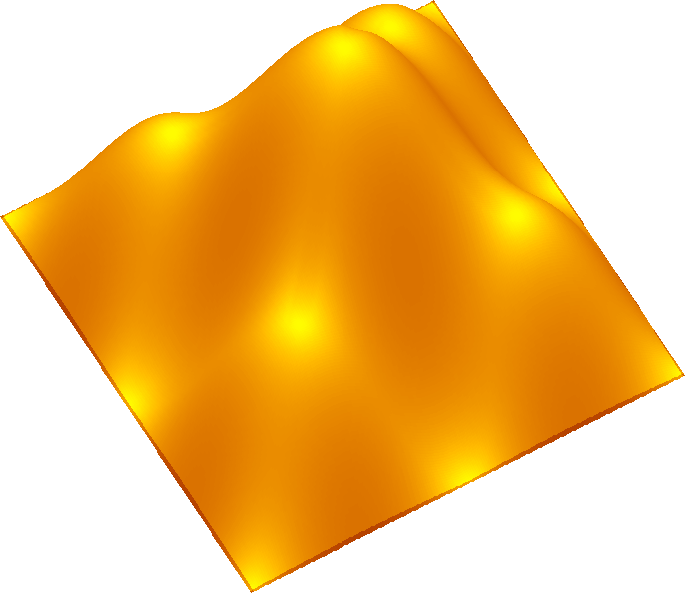} \\ 
\includegraphics[width=1.8cm]{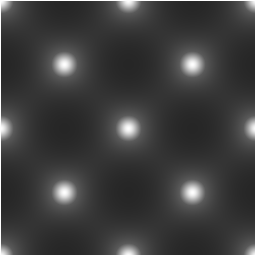} &\hspace{-0.2cm}
\includegraphics[width=1.8cm]{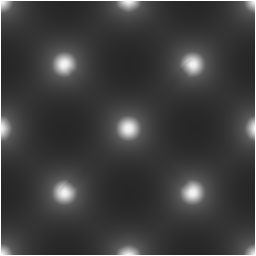} & \hspace{0.3cm}
\includegraphics[width=2.2cm]{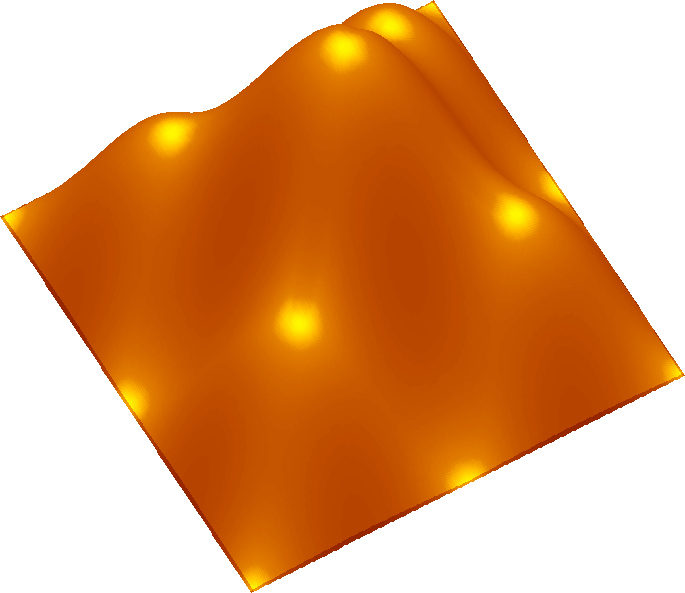} \\ 
\includegraphics[width=1.8cm]{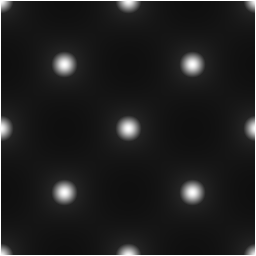} &\hspace{-0.2cm}
\includegraphics[width=1.8cm]{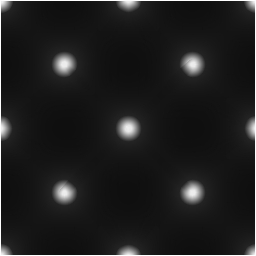} & \hspace{0.3cm}
\includegraphics[width=2.2cm]{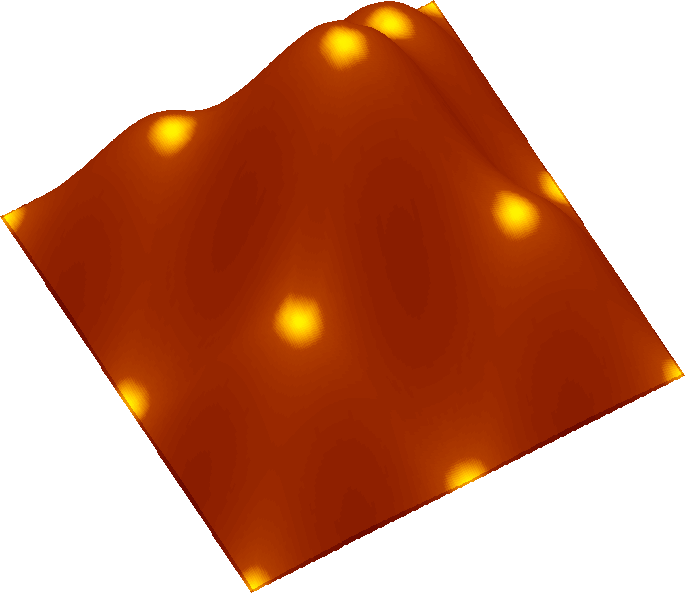} \\ 
  \hline 
\end{tabular}
 \caption{Synthetic sinusoidal function: example of concave/convex ambiguity for the PH--model with $k_S = 0,0.5,0.8$ for each row from the top to the bottom, respectively.}  \label{fig:sinusoidal_figures}
\end{figure}
Fig. \ref{fig:sinusoidal_figures} shows the results related to the PH--model with $k_S = 0, 0.5, 0.8$. 
In the first column one can see the input images generated by the PH--model using the values of the parameter before mentioned. In the second column we can see the output images computed a posteriori using the depth just computed and approximating the gradient via finite difference solver. What we can note is that even if the reconstructed a posteriori images match with the corresponding input images, the SL method always chooses the maximal solution even varying the parameters $k_D$ and $k_S$. By adding some informations as shown in the previous test it is possible to achieve a better result, but these additional informations are not available for real images.

\begin{figure}[h!]
\centering
   \begin{tabular}{p{3.4cm}p{3.4cm}p{3.4cm}}
   \hline\noalign{\smallskip}
  \hspace{0.cm} \textbf{L--model}  &  \hspace{0.cm} \textbf{ON--model} &  \hspace{0.cm} \textbf{PH--model} \\ 
  \hline \\
\includegraphics[width=3.45cm]{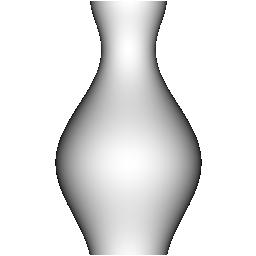} &\hspace{-0.2cm}
\includegraphics[width=3.45cm]{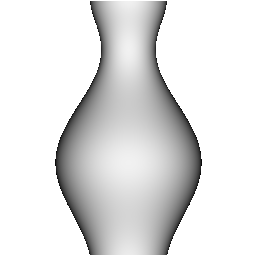} & \hspace{-0.2cm}
\includegraphics[width=3.45cm]{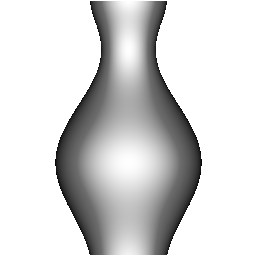} \vspace{0.mm}\\ 
\includegraphics[width=3.45cm]{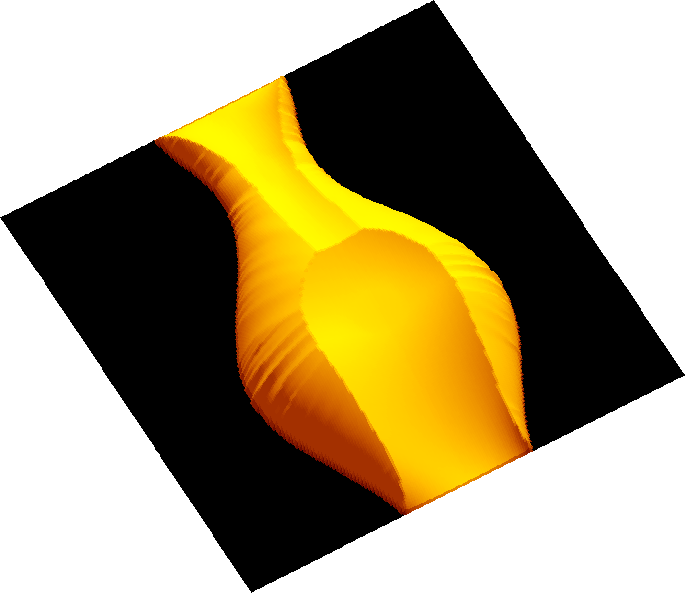} &\hspace{-0.2cm}
\includegraphics[width=3.45cm]{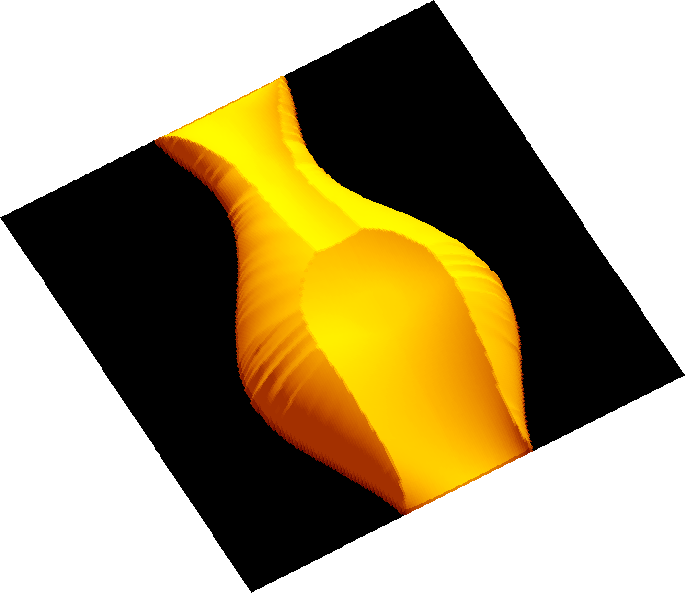} & \hspace{-0.2cm}
\includegraphics[width=3.45cm]{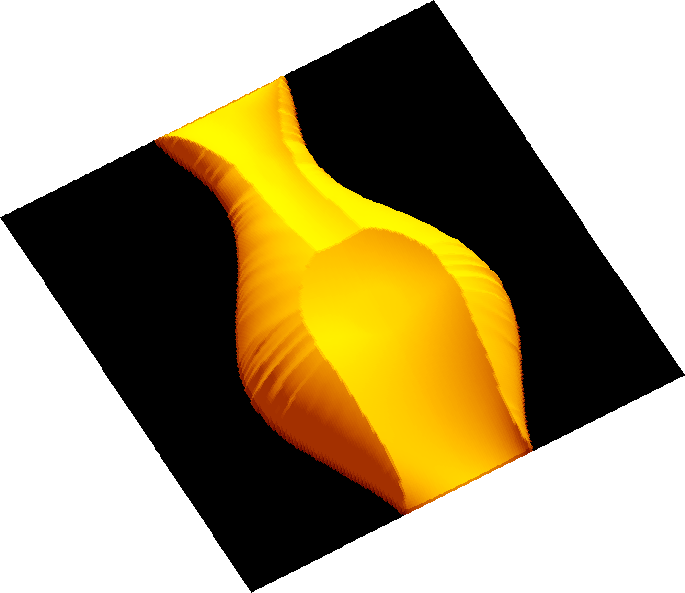} \\ 
\includegraphics[width=3.45cm]{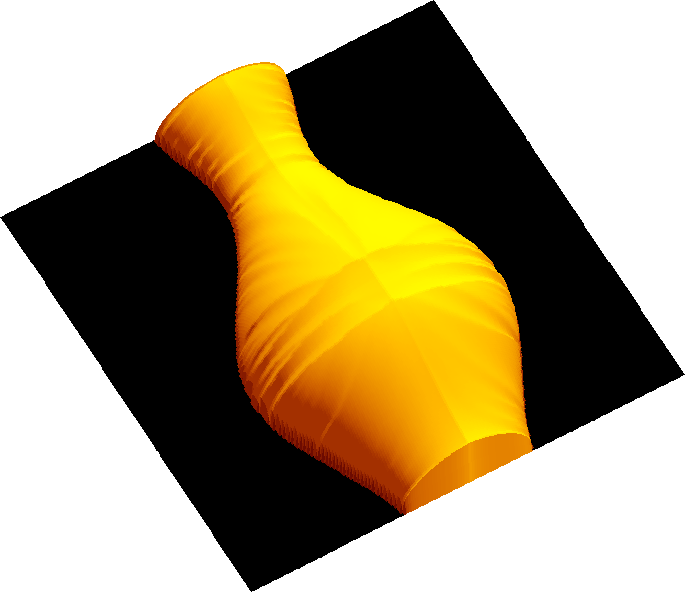} &\hspace{-0.2cm}
\includegraphics[width=3.45cm]{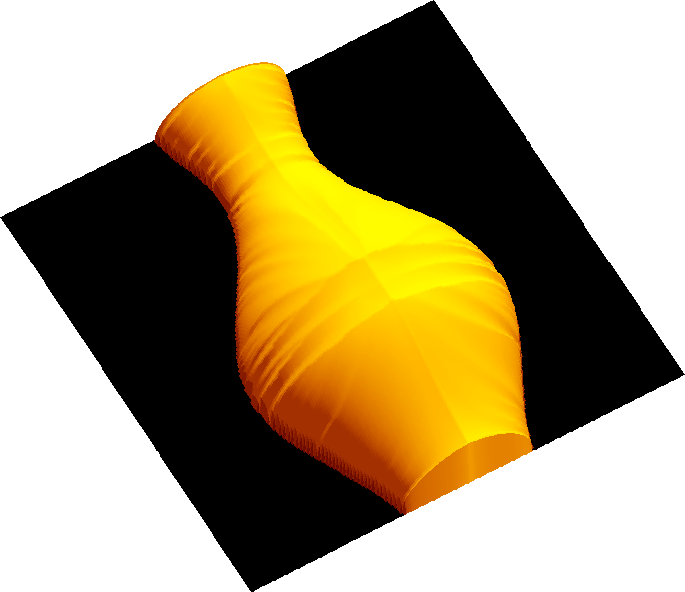} & \hspace{-0.2cm}
\includegraphics[width=3.45cm]{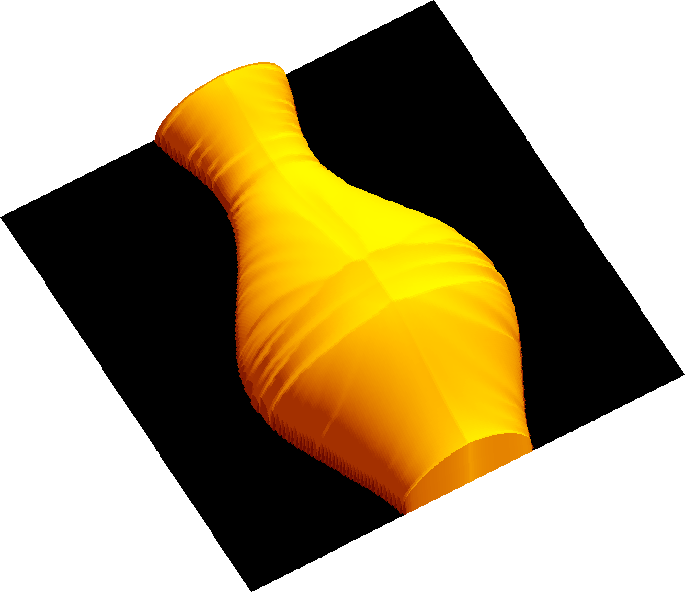} \\ 
  \hline 
\end{tabular}
 \caption{Synthetic vase under vertical light source ($\vomega=(0,0,1)$): example of concave/convex ambiguity solved by using correct Dirichlet BC. On the first row from left to right: input images generated by L--model, ON--model with $\sigma = 0.2$ and PH--model with $k_S = 0.4$. On the second row: 3D reconstruction with homogeneous Dirichlet BC. On the last row: 3D reconstruction with Dirichlet BC $u(x,y)= g(x,y)$.}
  \label{fig:vase_roleBC_figures}
  \end{figure}

\paragraph{Test 5: Role of the boundary conditions.}
With this fifth test we want to point out the role of the boundary condition (BC), showing how good BC can significantly improve the results on the 3D reconstruction. 
We will use the synthetic vase defined as follow
	\begin{equation}
		\begin{cases}
			u(x,y)=\sqrt{P(\bar{y})^2-x^2} &
				(x,y) \in G_{in}\\
			u(x,y)=g(x,y) &
				(x,y) \in G_{out},
		\end{cases}
		\label{eqn:vasosintetico}
	\end{equation}
where 
$\bar{y} := y/Y$,
\begin{eqnarray}
&P(\bar{y}) := &(-10.8\,\bar{y}^6+7.2\,\bar{y}^5+6.6\,\bar{y}^4-3.8\,\bar{y}^3\\
	&&-1.375\,\bar{y}^2+0.5\,\bar{y}+0.25)\,X  \nonumber
\end{eqnarray}		
and
	\[
		G_{in} := \{(x,y)| P(\bar{y})^2 > x^2\}.  
	\]
In Fig. \ref{fig:vase_roleBC_figures} one can see on the first row the input images (size $256\times256$) generated by the L--model, ON--model and the PH--model, from left to right respectively. On the second row we reported the 3D reconstruction with homogeneous Dirichlet BC ($g(x,y)=0$). As we can see, there is a concave/convex ambiguity in the reconstruction of the surface. If we consider the correct boundary condition, that is the height of the surface at the boundary of the silhouette that we can easily derive in this case being the object a solid of rotation, what we obtain is visible in the third row of the same Fig. \ref{fig:vase_roleBC_figures}.

  In Tables \ref{tab:SV_errors} and \ref{tab:SV_errors_BC} we can see the number of iterations, the CPU time and the error measures in  $L^2$ and $L^{\infty}$ norm for the method used with homogeneous and not homogeneous Dirichlet BC respectively. Looking at these errors we can note that in each table the values are almost the same for the different models. Comparing the values of the two tables one can see that we earn an order of magnitude using good BC, that confirms what we noted looking at Fig. \ref{fig:vase_roleBC_figures}.
\begin{table} [h!]
     \caption{Synthetic Vase:  Number of iterations, CPU time in seconds and $L^2$, $L^{\infty}$ errors on the surface with vertical light source $\vomega = (0,0,1)$ using homogeneous Drichlet BC.} 
       \label{tab:SV_errors} 
\begin{tabular}{p{2.1cm}p{1cm}p{1cm}p{1cm}p{1cm}}
\hline\noalign{\smallskip}
  SL--Schemes &  Iter. & [sec.]  & $L^2(S)$ & $L^\infty(S)$  \\
  \hline 
LAM  &  1010 &  0.54  &  0.1614 &  0.3015  \\
ON-02  &  1011 &  0.53 &   0.1613 &  0.3015  \\
ON-04  &  1008 & 0.54 & 0.1610 & 0.3015  \\
PH-s02 & 1007 &  0.54  &  0.1619 &  0.3015    \\
PH-s04  & 1009 &  0.55 &  0.1621 &  0.3015   \\
  \hline 
\end{tabular}
\end{table}

\begin{table} [h!]
     \caption{Synthetic Vase: Number of iterations, CPU time in seconds and $L^2$, $L^{\infty}$ errors on the surface with vertical light source $\vomega = (0,0,1)$ using not homogeneous Drichlet BC.}  
     \label{tab:SV_errors_BC} 
\begin{tabular}{p{2.1cm}p{1cm}p{1cm}p{1cm}p{1cm}}
\hline\noalign{\smallskip}
  SL--Schemes & Iter. & [sec.] & $L^2(S)$ & $L^\infty(S)$  \\
  \hline 
LAM  &  1337 &  0.70 &  0.0286 &  0.0569   \\
ON-02  &  1335 &  0.70 &  0.0284 &  0.0558   \\
ON-04  & 1341  & 0.70  & 0.0282 & 0.0562  \\
PH-s02 &  1330 &  0.70 &  0.0284 &  0.0560   \\
PH-s04  & 1330 &  0.70 &  0.0280 &  0.0558  \\
  \hline 
\end{tabular}
\end{table}

\paragraph{Test 6: Comparison with other numerical approximations.}
In this sixth and last test of this section dedicated to synthetic tests we will compare the performance of our semi-Lagrangian approach with other methods used in the literature.
For this reason, we will use a very common image used in the literature, that is the vase, defined in the previous test through the \eqref{eqn:vasosintetico}. 
More in details, we will compare the performance of our semi-Lagrangian method with the Lax-Friedrichs Sweeping (LFS) scheme adopted by Ahmed and Farag \cite{AF07} under vertical and oblique light source. Also these authors derive some similar HJ equations for the L--model in \cite{AF06}, and generalize this approach for various image conditions in \cite{AF07}, comparing their results on the only L--model with the results shown in \cite{Samaras_TPAMI2003} and the algorithms reported in \cite{ZTCS99}. 
    Unfortunately, we cannot compare our semi-Lagrangian approximation for the PH--model with no other schemes since, 
    to our knowledge, there are no table of errors for the PH--model under orthographic projection in the literature.
  In order to do the comparison, we will consider the vase image of size $128\times128$ as used in the other papers. 
  We start to remind the error estimations used: given a vector $\vA$ representing the reference depth map  
  on the grid and a vector $\widetilde \vA$ representing its approximation, we define the  error vector as $\ve=\vA-\widetilde \vA$ and 
\begin{eqnarray*}\label{eqn:error_norm_definition}
            & err_1:= ||\ve||_{L^1} & =  \frac{1}{N}\sum_i |e_i|\\
            & err_2 := ||\ve||_{L^2} & =  \left\lbrace \frac{1}{N}\sum_i |e_i|^2  \right\rbrace^{1/2} 
\end{eqnarray*}
where $N$ is the total number of grid points used for the computation, i.e. the grid points belonging to $G_{in}$. These estimators are called mean and standard deviation of the absolute error.
In Table \ref{tab:vase_comparison_errors} we compare the error measures for the different SfS algorithms under the L--model with $\vomega = (0,0,1)$. What we can note is that
our semi-Lagrangian method is better than the other ones also if we consider Dirichlet boundary condition equal to zero (as used by Ahmed and Farag in their work \cite{AF07}), but the better result is with correct BC, shown in bold in the last row.
\begin{table} [h!]
     \caption{Synthetic vase: error measures related to the different methods for the L--model with vertical light source $\vomega = (0,0,1)$. In bold the best performances.}   \label{tab:vase_comparison_errors} 
     \begin{tabular}{p{4.7cm}p{1.2cm}p{1.2cm}}
\hline\noalign{\smallskip}
  Methods & $err_1$ & $err_2$  \\
  \hline 
Best \cite{ZTCS99} & 8.1 & 11.1 \\
\cite{Samaras_TPAMI2003} & 2.8 & 2.0 \\
\cite{AF07} & 0.22 & 0.4 \\
Our proposed ($BC=0$) & 0.1570 &  0.1717 \\
{\bf \boldmath Our proposed ($BC \not= 0$)} & {\bf 0.0349} &  {\bf 0.0385} \\
  \hline 
\end{tabular}
\end{table}
In Table \ref{tab:vase_comparison_obl_errors} we can see the same methods applied to the L--model but under a different  light source,  that is $(1,0,1)$. Also in this case, our approach obtains always the smallest errors and the best is with not homogeneous Dirichlet BC, as noted before.
\begin{table} [h!]
     \caption{Synthetic vase: error measures related to the different methods for the L--model with oblique light source $(1,0,1)$. In bold the best performances.}   \label{tab:vase_comparison_obl_errors} 
     \begin{tabular}{p{4.7cm}p{1.2cm}p{1.2cm}}
\hline\noalign{\smallskip}
  Methods & $err_1$ & $err_2$  \\
  \hline 
Best \cite{ZTCS99} & 7.9 & 13.9 \\
\cite{Samaras_TPAMI2003} & 4.1 & 2.6 \\
\cite{AF07} & 1.2 & 2.2  \\
Our proposed ($BC=0$) & 0.0683 & 0.1061   \\
{\bf \boldmath Our proposed ($BC \not= 0$)} & {\bf 0.0218} &  {\bf 0.0242} \\
  \hline 
\end{tabular}
\end{table}

\begin{figure}[h!]
\centering
   \begin{tabular}{p{3.4cm}p{3.4cm}p{3.4cm}}
   \hline\noalign{\smallskip}
  \hspace{0.3cm} \textbf{L--model}  &  \hspace{-0.1cm} \textbf{L--model} &  \hspace{-0.35cm} \textbf{ON--model} \\ 
   \hspace{0.5cm} vertical  &  \hspace{0.1cm} oblique &  \hspace{0.cm} vertical \\ 
  \hline \\
\hspace{-0.3cm}\includegraphics[width=0.24\textwidth]{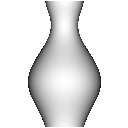} &\hspace{-0.4cm}
\hspace{-0.2cm}\includegraphics[width=0.24\textwidth]{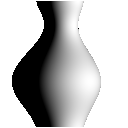} &\hspace{-0.4cm}
\hspace{-0.1cm}\includegraphics[width=0.24\textwidth]{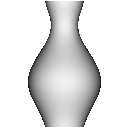} \\ 
\includegraphics[width=0.25\textwidth]{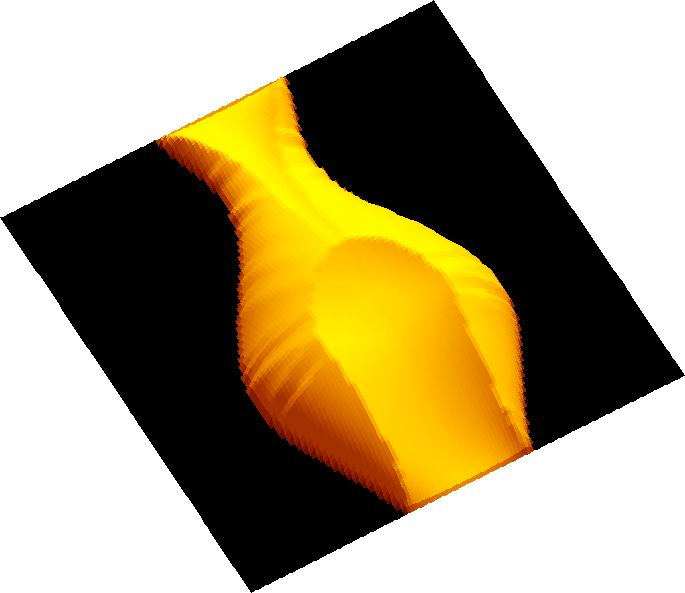} &\hspace{-0.4cm}
\includegraphics[width=0.25\textwidth]{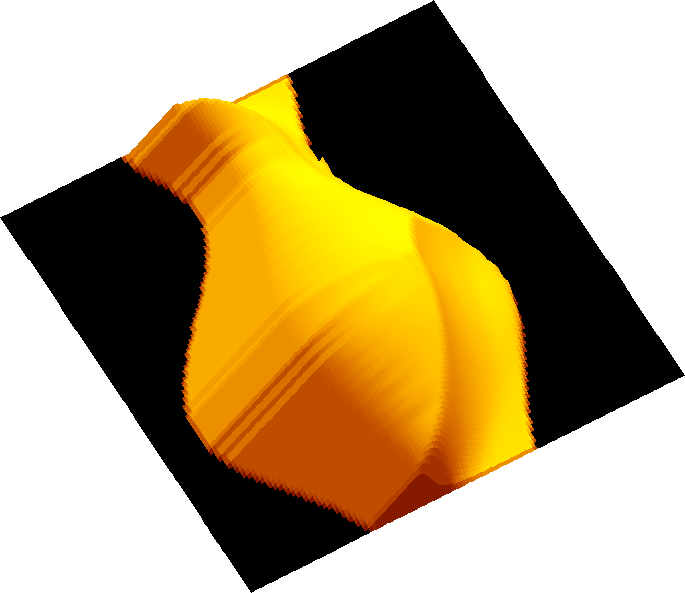} & \hspace{-0.4cm}
\includegraphics[width=0.25\textwidth]{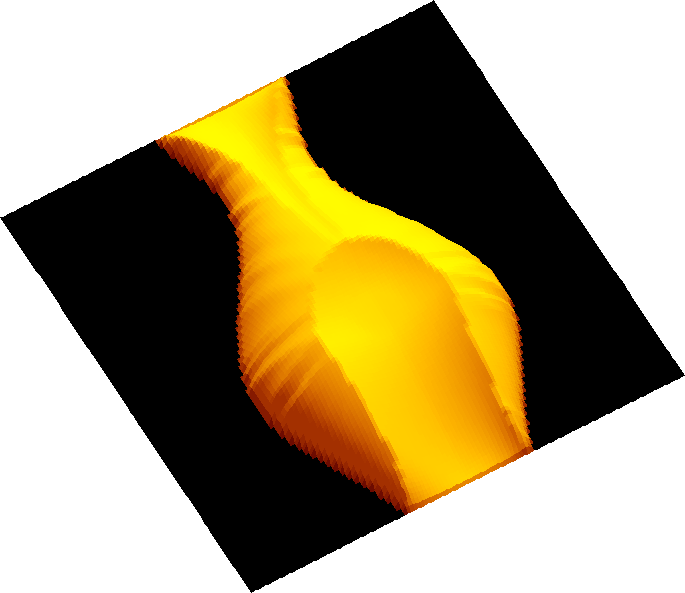} \\ 
\includegraphics[width=0.25\textwidth]{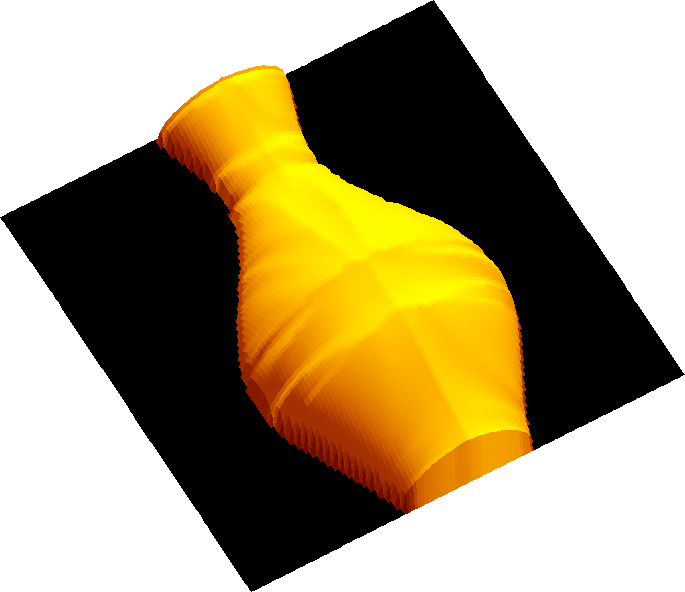} & \hspace{-0.4cm}
\includegraphics[width=0.25\textwidth]{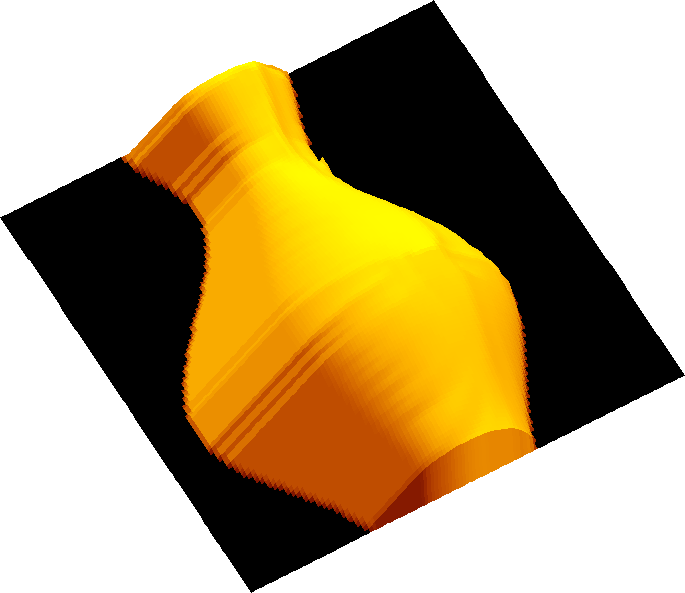} & \hspace{-0.4cm}
\includegraphics[width=0.25\textwidth]{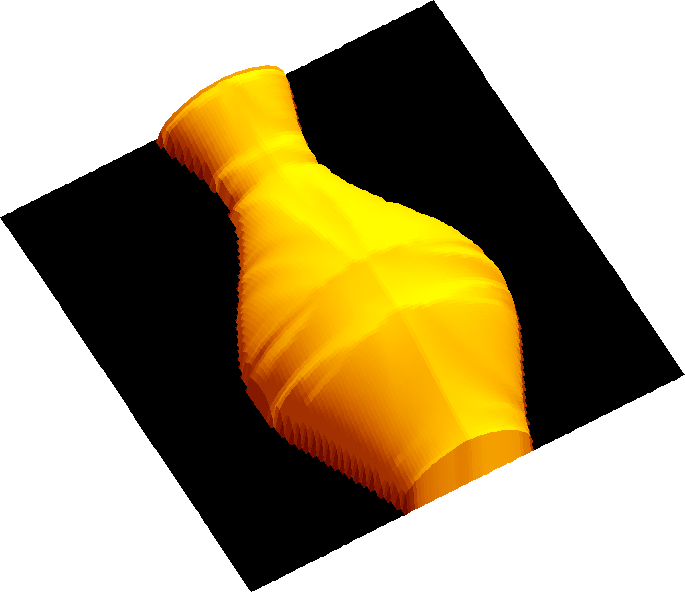} \\ 
\includegraphics[width=0.25\textwidth]{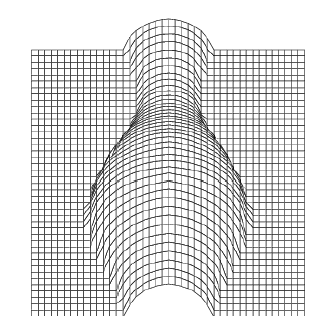} &\hspace{-0.4cm}
\includegraphics[width=0.25\textwidth]{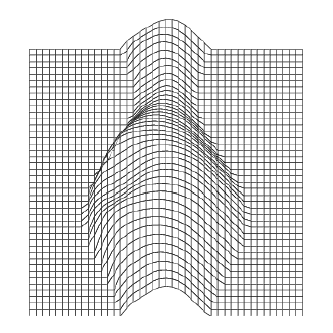} &\hspace{-0.4cm}
\includegraphics[width=0.25\textwidth]{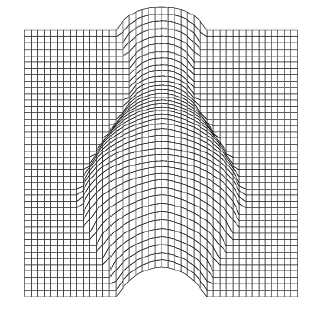} \\
  \hline 
\end{tabular}
 \caption{
 Synthetic vase: From top to bottom, Input images, recovered shapes by our approach with homogeneous Dirichlet BC and with not homogeneous BC, recovered shape by \cite{AF07}. 
First column: L--model with vertical light source (0,0,1). Second column: L--model with oblique light source (1,0,1). Third column: ON--model with $\sigma = 0.2$, $\vomega = (0,0,1)$, $\vV = (0,0,1)$. Input images size: $128\times128$.}
 \label{fig:vase_corresp_to_tables} 
\end{figure}

Only with respect to the LFS used by Ahmed and Farag, we can compare the performance of our semi-Lagrangian scheme under the ON--model as well, since the other authors only consider the L--model. In this context, we show in Table \ref{tab:vase_comparison_ON_errors} the error measures for the two SfS algorithms 
with $\sigma = 0.2$, under vertical position of light source and viewer ($\vomega = (0,0,1), \vV = (0,0,1)$). 
\begin{table} [h!]
     \caption{Synthetic vase: error measures related to the SL and  LFS methods for the ON--model with $\sigma = 0.2$ under vertical light source ($\vomega = (0,0,1)$) and vertical viewer $\vV = (0,0,1)$. In bold the best performances.} \label{tab:vase_comparison_ON_errors} 
     \begin{tabular}{p{4.7cm}p{1.2cm}p{1.2cm}}
\hline\noalign{\smallskip}
  Methods & $err_1$ & $err_2$  \\
  \hline 
\cite{AF07} & 0.6 & 1.0  \\
Our proposed ($BC=0$) & 0.1568 &  0.1715   \\
{\bf \boldmath Our proposed ($BC \not= 0$)} & {\bf 0.0348} &  {\bf 0.0384} \\ 
  \hline 
\end{tabular}
\end{table}
As before, the best result is obtained using the semi-Lagrangian scheme with not homogeneous Dirichlet BC. 
The reconstructions corresponding to the error measures shown in the three last Tables \ref{tab:vase_comparison_errors}, \ref{tab:vase_comparison_obl_errors}, \ref{tab:vase_comparison_ON_errors} obtained applying our scheme and compared to \cite{AF07} are shown in Fig. \ref{fig:vase_corresp_to_tables}. The reconstruction obtained by the two methods are comparable. 
In particular, in the second column regarding the oblique light source case we can note that our scheme reconstructs a surface that incorporates the black shadow part (see \cite{FSS03} for more details on this technique), avoiding the effects of ``dent" present in the reconstruction obtained by Ahmed and Farag visible in the last row, second column.

Finally, in Table \ref{tab:vase_iter&CPU} we reported the number of iterations and the CPU time in seconds with the comparison with respect to \cite{AF07}.
\begin{table} [h!]
     \caption{Synthetic vase: iterations and CPU time in seconds for the L--model and the ON--model with vertical light source $\vomega = (0,0,1)$. Image size: $128\times128$.}   
     \label{tab:vase_iter&CPU}
  \begin{tabular}{p{4.3cm}p{2.2cm}p{1.1cm}p{1.1cm}}
\hline\noalign{\smallskip}
\hspace{1.2cm}  Schemes & Model & Iter.  & $[sec.]$ \\
  \hline 
\hspace{1.5cm}\cite{AF07}  & L--model  & \hspace{0.2cm} - & 0.5   \\ 
Our approach ($BC = 0$) & L--model   &  611 &  {\bf 0.09}   \\ 
Our approach ($BC \not= 0$) & L--model & 792 &  0.11 \\
\hspace{1.5cm}\cite{AF07}  & ON--model  & \hspace{0.2cm} - & 1.5   \\ 
Our approach ($BC = 0$) & ON--model   &  612 &  {\bf 0.09}  \\ 
Our approach ($BC \not= 0$) & ON--model & 791 &  0.12 \\
  \hline 
\end{tabular}
\end{table}
This shows that the SL-scheme is competitive also in terms of CPU time. Of course, in the case of oblique light source the number of iterations, and hence the CPU time needed is much more bigger. Just think that for the reconstruction under the L--model with oblique light source $(1,0,1)$ and $BC \not= 0$ visible in the second column of the third row of Fig. \ref{fig:vase_corresp_to_tables}, we need $15513$ iterations (CPU time: 147.40 seconds). If we double the size of the input image, considering the vase $256\times256$ visible in Fig. \ref{fig:vase_roleBC_figures}, for the reconstruction visible in the third row, first column of the same Fig. \ref{fig:vase_roleBC_figures}, we need 30458 iterations to get convergence, that we obtain in 1163 seconds.

\begin{figure}[h!]
\centering
   \begin{tabular}{p{3.4cm}p{3.4cm}}
   \hline\noalign{\smallskip}
  \hspace{0.3cm} \textbf{Input}  &  \hspace{-0.35cm} \textbf{reconstruction} \\ 
  \hline \\
\hspace{-0.1cm}\includegraphics[width=0.25\textwidth]{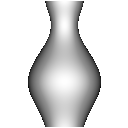} &\hspace{-0.cm}
\hspace{-0.1cm}\includegraphics[width=0.255\textwidth]{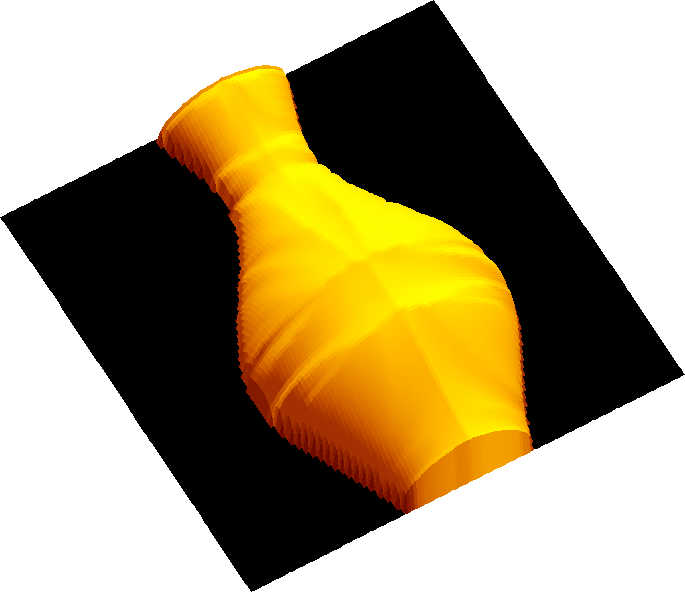} \\ 
\hspace{-0.1cm}\includegraphics[width=0.25\textwidth]{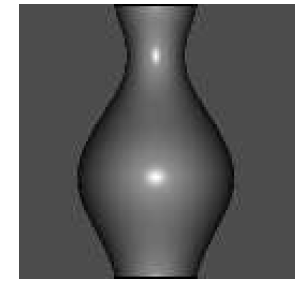} &\hspace{-0.cm}
\vspace{-0.2cm} \hspace{-0.cm}\includegraphics[width=0.25\textwidth]{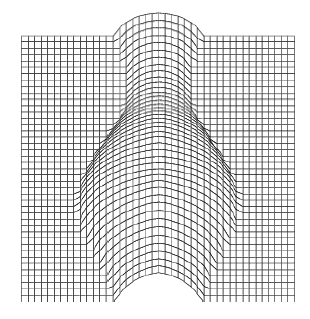} \\ 
  \hline 
\end{tabular}
 \caption{
 Synthetic vase: The first row shows the input image and the recovered shape by our approach with not homogeneous BC based on the PH--model with $k_S = 0.2$, $\vomega = (0,0,1)$ and $\vV = (0,0,1)$. 
The second row shows the input image and the recovered shape by the approximate Ward’s method illustrated in \cite{AF07} with 
$\sigma = 0.2$, $\rho_d = 0.67$, $\rho_s = 0.075$, $\vomega = (0,0,1)$, $\vV = (0,0,1)$. Input images size: $128\times128$.}
 \label{fig:vase_Ward&phong} 
\end{figure}

Since it is difficult to compare the performance of our approach based on the PH--model with other schemes based on the same reflectance model under orthographic projection for the consideration made before, we report in Fig. \ref{fig:vase_Ward&phong} the performances obtained by our method based on the PH--model compared to the approximate Ward’s method on the vase test (Cf. Fig. 7 in \cite{AF07}). What we can see is that both the two methods reconstruct the surface in a quite good way, without particular distinctions in the goodness of the reconstructions. In order to analyze the results not only in a qualitative way but also in a quantitative one, we report the mean and the standard deviation of the absolute errors in Table \ref{tab:vase_Ward&phong}.  Looking at this table, we can note that our approach seems to be superior, obtaining the reconstruction with errors smaller of one order with respect to the other model.

\begin{table} [h!]
     \caption{
     Synthetic vase: error measures related to the cases shown in Fig. \ref{fig:vase_Ward&phong}. In bold the best performances.}   
     \label{tab:vase_Ward&phong} 
     \begin{tabular}{p{4.7cm}p{1.2cm}p{1.2cm}}
\hline\noalign{\smallskip}
  Methods & $err_1$ & $err_2$  \\
  \hline 
Ward in \cite{AF07} & 0.8 & 1.3  \\
{\bf \boldmath PH--model} & {\bf 0.03} &  {\bf 0.04} \\  
  \hline 
\end{tabular}
\end{table}

\subsection{Real tests}
\begin{figure}[h!]
\centering
 \subfigure[Agamemnon Input]
   {\includegraphics[width=0.32\textwidth]{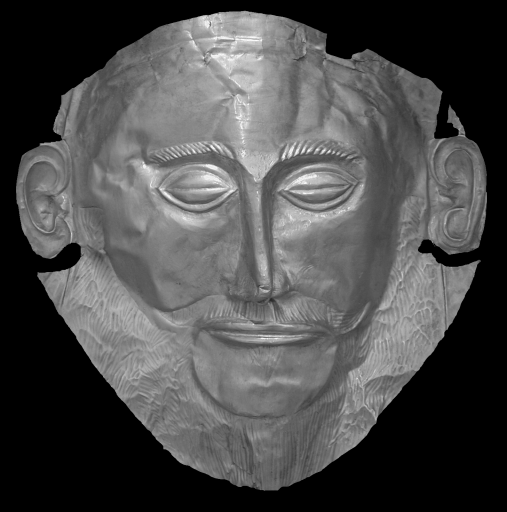} \label{fig:input_Agamemnon}  }
 \hspace{2mm}
  \subfigure[Agamemnon Mask]
  {\includegraphics[width=0.32\textwidth]{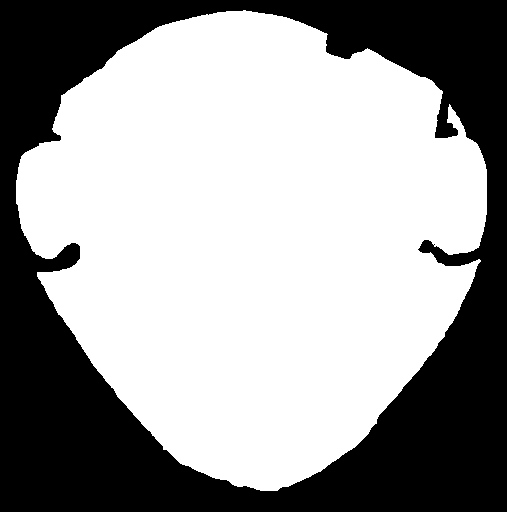} \label{fig:mask_Agamemnon}  }
 \caption{Agamemnon images (size $507\times512$):
            (a)~Input image; 
            (b)~Mask.}
 \label{fig:Agamemnon_images}
 \end{figure}

In this subsection we consider real input images. We start considering the golden mask of Agamemnon taken from \cite{wikipedia_Agamemnon} and then modified in order to get a picture in gray tones. The size of the modified image really used is $507\times512$. The input image is visible in Fig. \ref{fig:input_Agamemnon}, the associated mask used for the 3D reconstruction in Fig. \ref{fig:mask_Agamemnon}. 
The second real test is concerning the real vase (RV in the following) visible in Fig. \ref{fig:RV_images}, taken from \cite{DFS08}. The size of the input image shown in Fig. \ref{fig:input_RV} is $256\times256$. The reconstruction domain $\Omega_{RV}$, shown in Fig. \ref{fig:mask_RV}, is constituted by the pixels situated on the vase. 
For the real cases, the input image is the same for all the models and we can compute only errors on the images since we do not know the height of the original surface.  
For the real tests we will use the same stopping criterion for the iterative method before defined for the synthetic tests, i.e. $|W^{k+1} - W^k|_{max} \leq \eta$. 

\paragraph{Test 7: Agamennon mask.}
 
For this test we will compare the results regarding 3D reconstruction of the surface obtained with a vertical light source $\vomega_{vert} = (0,0,1)$ and an oblique light source $\vomega_{obl} = (0, 0.0995, 0.9950)$. 
\begin{table}[h!]
 \caption{Real Agamemnon mask: parameter values used in the models with vertical light source $\vomega = (0,0,1)$.}  \label{tab:Agam_parameter}
   \begin{tabular}{p{1.3cm}p{1.1cm}p{1.1cm}p{1.1cm}p{1.1cm}}
\hline\noalign{\smallskip}
Model & \hspace{0.cm}$\sigma$ & \hspace{0.cm}$k_D$ & \hspace{0.cm}$k_S$ & \hspace{0.cm}$\alpha$   \\
\hline 
LAM & 	\hspace{0.cm} & & &  \\
ON-04 & 	 \hspace{0.cm}0.4 & \hspace{0.cm} & \hspace{0.cm} & \hspace{0.cm}   \\
ON-08 & 	 \hspace{0.cm}0.8 &  \hspace{0.cm} & \hspace{0.cm} &\hspace{0.cm}   \\
ON-10 & 	 \hspace{0.cm}1 & \hspace{0.cm} & \hspace{0.cm} & \hspace{0.cm} \\
PH-s04 & 	 \hspace{0.cm} &  \hspace{0.cm} 0.6 & \hspace{0.cm} 0.4 &\hspace{0.cm} 1   \\
PH-s08 & 	\hspace{0.cm} &\hspace{0.cm} 0.2 &\hspace{0.cm} 0.8 & \hspace{0.cm} 1 \\
PH-s10 &	 	\hspace{0.cm} &  \hspace{0.cm} 0 & \hspace{0.cm} 1 & \hspace{0.cm} 1   \\
  \hline 
\end{tabular}
\end{table}
The values of the parameters used in this test are reported in Table \ref{tab:Agam_parameter}.
For a vertical light source, we refer to Table \ref{tab:Agam_iter&CPU} for the number of iterations and the CPU time (in seconds) and to Table \ref{tab:Agam_errors} for the errors obtained with a tolerance $\eta=10^{-8}$ for the stopping rule of the iterative process. 
Clearly,  the number of iteration and the errors of the two non-Lambertian models are the same of the classical Lambertian model when $\sigma$ for the ON--model and $k_S$ for the PH--model are equal to zero (and for this reason we do not report them in the tables).
In all the other cases, the non-Lambertian models are faster in terms of CPU time and need a lower number of iterations with respect to the L--model. 

\begin{table} [h!]
     \caption{Real Agamemnon mask: iterations and CPU time in seconds for the models with vertical light source $\vomega = (0,0,1)$.}   \label{tab:Agam_iter&CPU} 
  \begin{tabular}{p{2.5cm}p{1.5cm}p{1.5cm}}
\hline\noalign{\smallskip}
  SL--Schemes & Iter.  & $[sec.]$ \\
  \hline 
LAM  & 3921    &   24.48 \\ 
ON-04  & 2751 &  12.48 \\  
ON-08  & 1943 &  11.41 \\  
ON-10  & 1818 &  8.79 \\  
PH-s04  & 2127 & 9.89  \\  
PH-s08  & 1476 &  6.90 \\  
PH-s10  & 1325 & 6.33 \\ 
  \hline 
\end{tabular}
\end{table}

In Table \ref{tab:Agam_errors}  we can observe that the $L^2$ errors produced by the ON--model increase by increasing the value of $\sigma$. However, the $L^{\infty}$ errors are lower than the error obtained with the Lambertian model.
With respect to the PH--model, all the errors increase by increasing the value of the parameter $k_S$, as observed for synthetic images.

\begin{figure}[h!]
\centering
   \begin{tabular}{p{3.4cm}p{3.4cm}p{3.4cm}}
   \hline\noalign{\smallskip}
  \hspace{0.5cm} \textbf{LAM}  &  \hspace{0.5cm} \textbf{ON4} &  \hspace{0.5cm} \textbf{PH4} \\ 
  \hline \\
\includegraphics[width=0.25\textwidth]{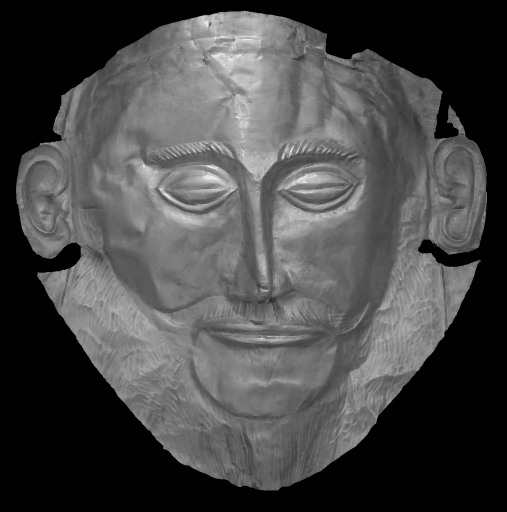} &
\includegraphics[width=0.25\textwidth]{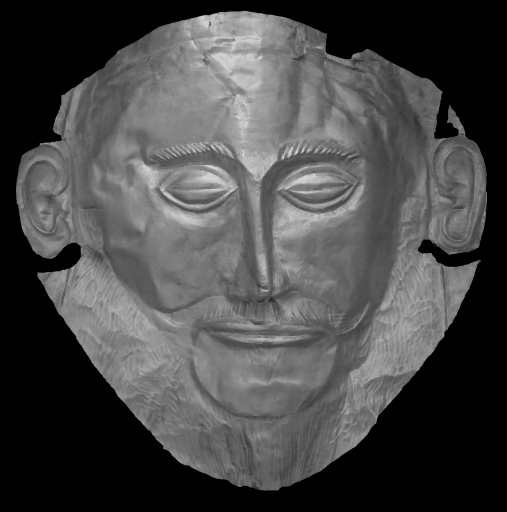} &
\includegraphics[width=0.25\textwidth]{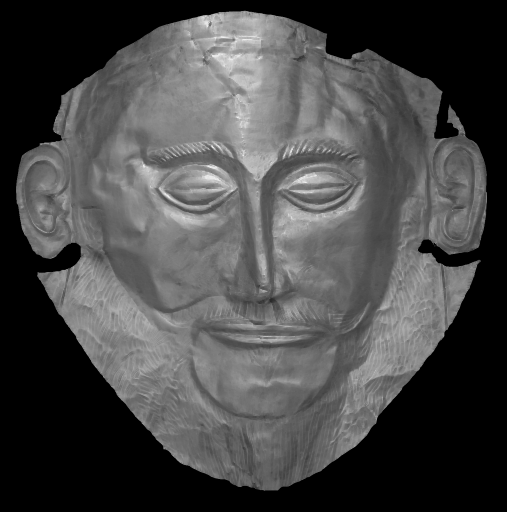} \\ 
\includegraphics[width=0.255\textwidth]{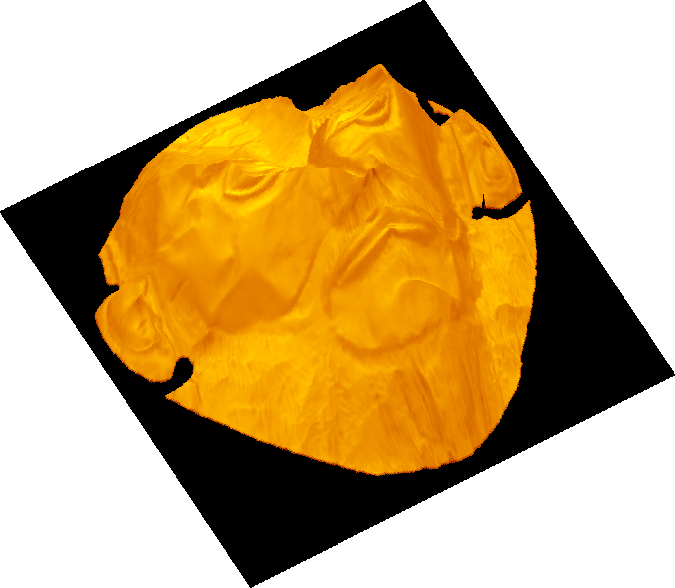} &
\includegraphics[width=0.255\textwidth]{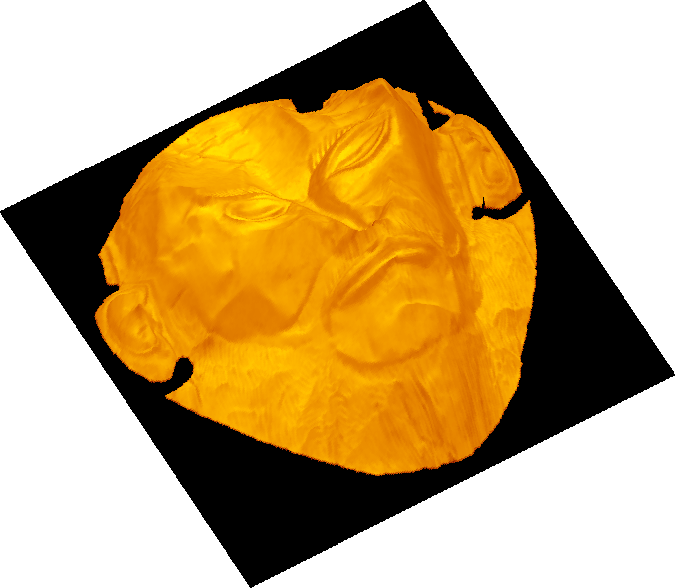} &
\includegraphics[width=0.255\textwidth]{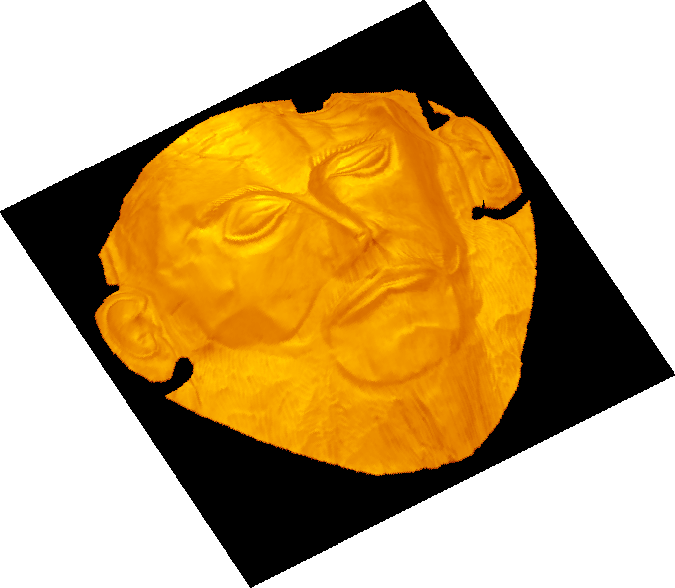} \\ 
\includegraphics[width=0.255\textwidth]{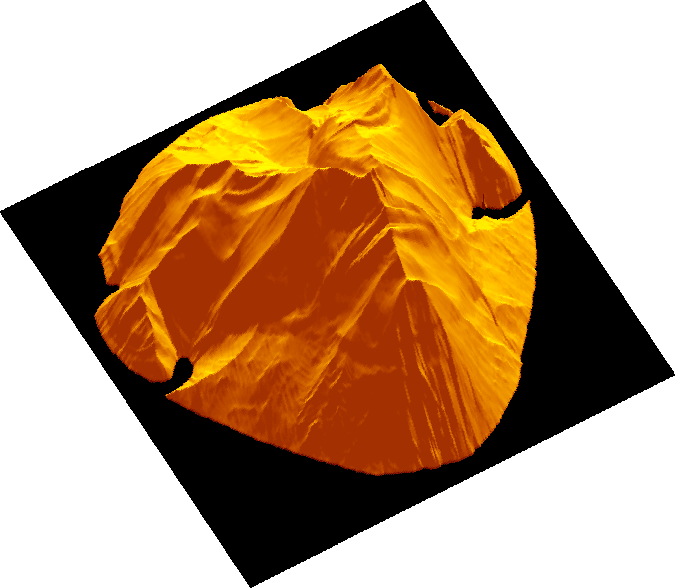} &
\includegraphics[width=0.255\textwidth]{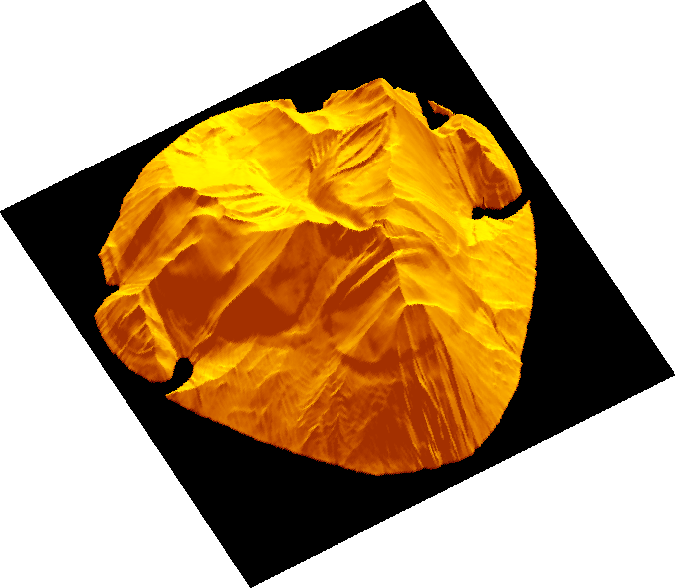} &
\includegraphics[width=0.255\textwidth]{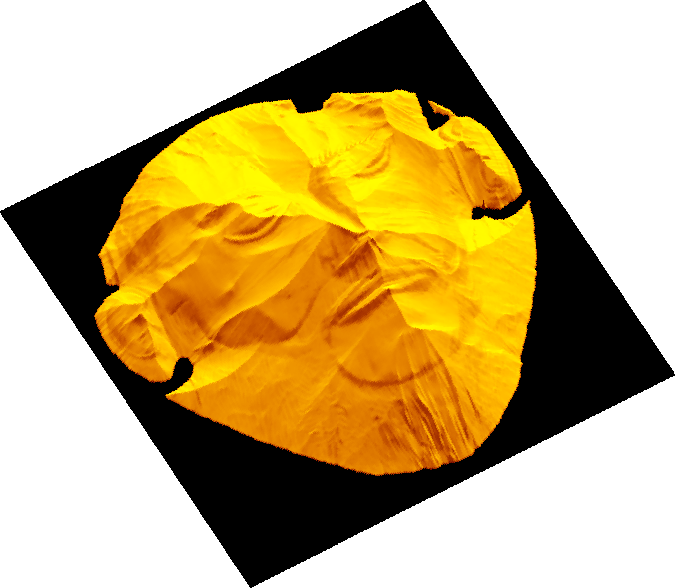} \\ 
  \hline 
\end{tabular}
 \caption{Agamennon mask: results with vertical light source. On the first row the output images, on the second row the 3D reconstruction with vertical view, on the third row the 3D reconstruction with oblique view.}  \label{fig:Agamennon_vert_results}
\end{figure}

\begin{table} [h!]
     \caption{Real Agamemnon mask:  $L^2$, $L^{\infty}$ errors with vertical light source $\vomega = (0,0,1)$.}   \label{tab:Agam_errors} 
\begin{tabular}{p{2.1cm}p{2cm}p{2cm}}
\hline\noalign{\smallskip}
  SL--Schemes & $L^2(I)$ & $L^\infty(I)$  \\
  \hline 
LAM     &  0.0371 &  0.4745   \\
ON-04  &  0.0375 &  0.4627  \\
ON-08  &  0.0440 &  0.4627  \\
ON-10   &  0.0501 &  0.4627   \\
PH-s04  &  0.0383 &  0.4824 \\
PH-s08  &  0.0391 &  0.4941   \\
PH-s10  &  0.0393 &  0.5098 \\
  \hline 
\end{tabular}
\end{table}
In Fig. \ref{fig:Agamennon_vert_results} we can see the output image and the 3D reconstruction in a single case for each models. What we can note is that no big improvements we can obtain visually.

For the oblique light case, we consider the values for the parameters reported in Table \ref{tab:Agam_parameter_obl}.
\begin{table}[h!]
 \caption{Real Agamemnon mask: parameter values used in the models with an oblique light source $\vomega_{obl} = (0, 0.0995, 0.9950)$.}  \label{tab:Agam_parameter_obl}
   \begin{tabular}{p{1.3cm}p{1.1cm}p{1.1cm}p{1.1cm}p{1.1cm}}
\hline\noalign{\smallskip}
Model & \hspace{0.cm}$\sigma$ & \hspace{0.cm}$k_D$ & \hspace{0.cm}$k_S$ & \hspace{0.cm}$\alpha$   \\
\hline 
LAM & 	\hspace{0.cm} & & &  \\
ON-01 & 	 \hspace{0.cm}0.1 & \hspace{0.cm} & \hspace{0.cm} & \hspace{0.cm}   \\
ON-02 & 	 \hspace{0.cm}0.2 &  \hspace{0.cm} & \hspace{0.cm} &\hspace{0.cm}   \\
ON-03 & 	 \hspace{0.cm}0.3 & \hspace{0.cm} & \hspace{0.cm} & \hspace{0.cm} \\
PH-s02 & 	\hspace{0.cm} &\hspace{0.cm} 0.8 &\hspace{0.cm} 0.2 & \hspace{0.cm} 1 \\
PH-s03 &	 	\hspace{0.cm} &  \hspace{0.cm} 0.7 & \hspace{0.cm} 0.3 & \hspace{0.cm} 1   \\
PH-s04 & 	 \hspace{0.cm} &  \hspace{0.cm} 0.6 & \hspace{0.cm} 0.4 &\hspace{0.cm} 1   \\
  \hline 
\end{tabular}
\end{table}

Looking at Table \ref{tab:Agam_iter&CPU_obl} we can note that the oblique cases require higher CPU time with respect to the 
 vertical cases due to the fact that the equations are more complex because of additional terms involved.
Because of these additional terms involved in the oblique case, in Table \ref{tab:Agam_errors_obl} we have reported the results obtained using the parameters shown in Table \ref{tab:Agam_parameter_obl} with a value of the tolerance $\eta$ for the stopping rule of the iterative method equal to $10^{-3}$. This is the maximum accuracy achieved by the non-Lambertian models since roundoff errors coming from several terms occur and limit the accuracy.

\begin{table} [h!]
     \caption{Real Agamemnon mask: number of iterations and CPU time in seconds for the different models with  oblique light source $\vomega_{obl} = (0, 0.0995, 0.9950)$.}   \label{tab:Agam_iter&CPU_obl} 
  \begin{tabular}{p{2.5cm}p{1.5cm}p{1.5cm}}
\hline\noalign{\smallskip}
  SL--Schemes & Iter.  & $[sec.]$ \\
  \hline 
LAM      &  321  & 117.9   \\ 
ON-01  & 315 & 246.0  \\  
ON-02  & 361 & 281.5   \\  
ON-03  & 396 & 264.6   \\   
PH-s02  & 427 & 285.2  \\  
PH-s03  & 564 & 373.6  \\ 
PH-s04  & 680 & 484.1  \\  
  \hline 
\end{tabular}
\end{table}

\begin{figure}[h!]
\centering
   \begin{tabular}{p{3.4cm}p{3.4cm}p{3.4cm}}
   \hline\noalign{\smallskip}
  \hspace{0.5cm} \textbf{LAM}  &  \hspace{0.5cm} \textbf{ON4} &  \hspace{0.5cm} \textbf{PH4} \\ 
  \hline \\
\includegraphics[width=0.25\textwidth]{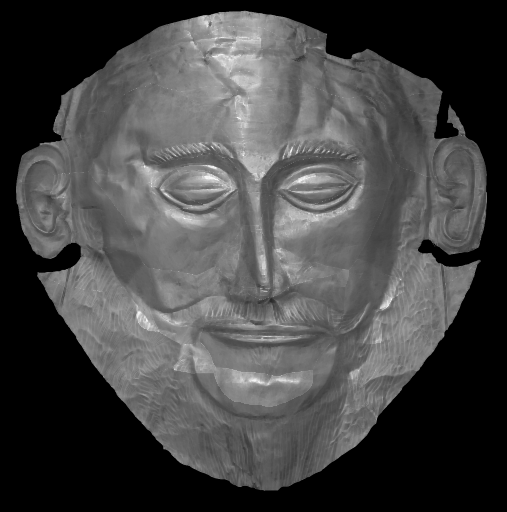} &
\includegraphics[width=0.25\textwidth]{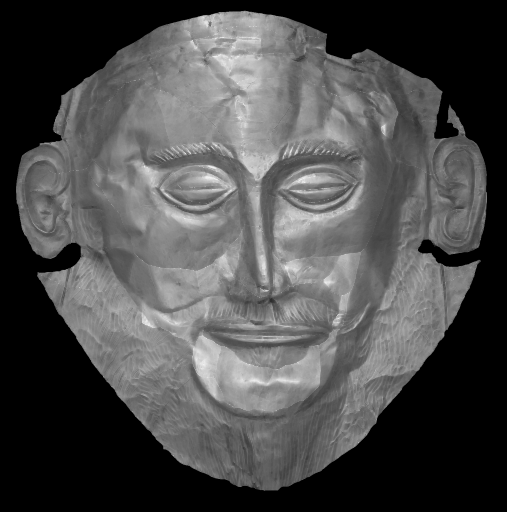} &
\includegraphics[width=0.25\textwidth]{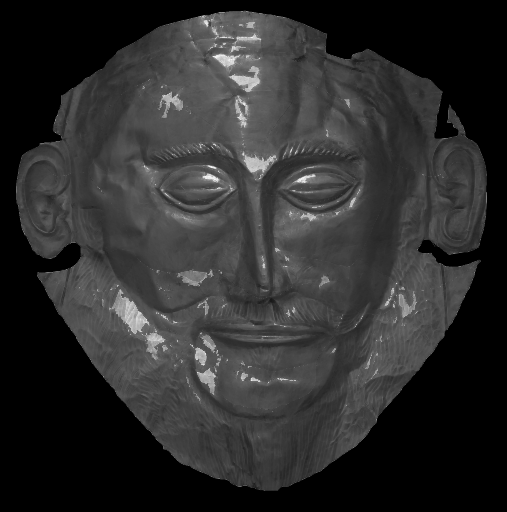} \\ 
\includegraphics[width=0.255\textwidth]{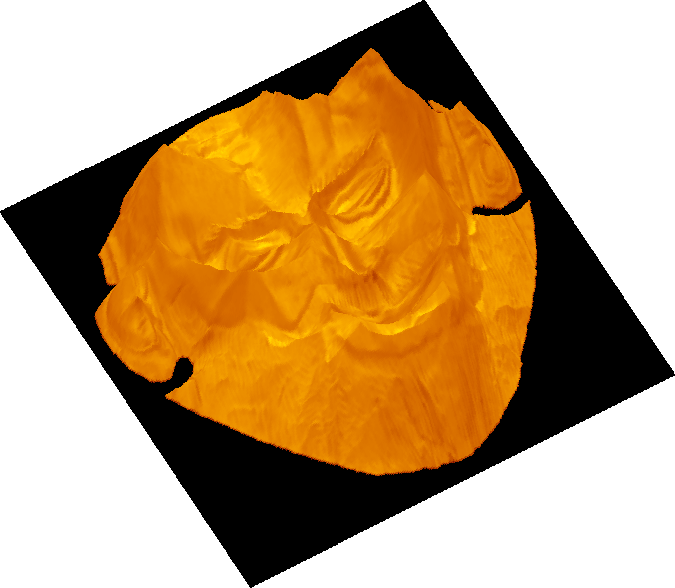} &
\includegraphics[width=0.255\textwidth]{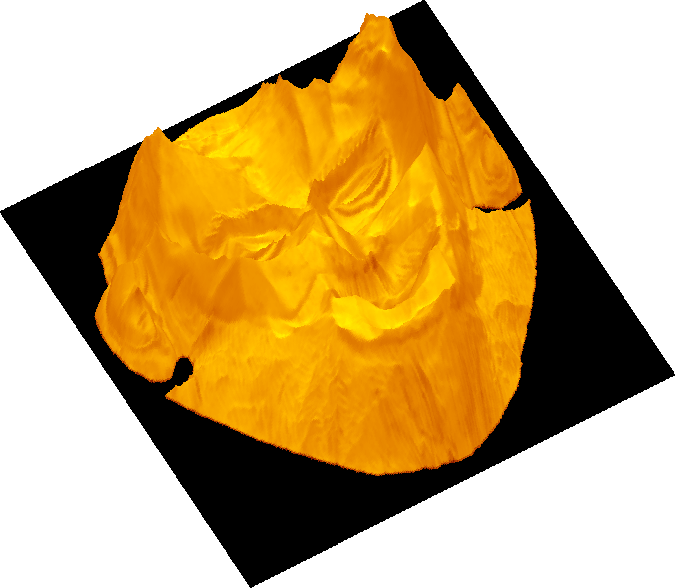} &
\includegraphics[width=0.255\textwidth]{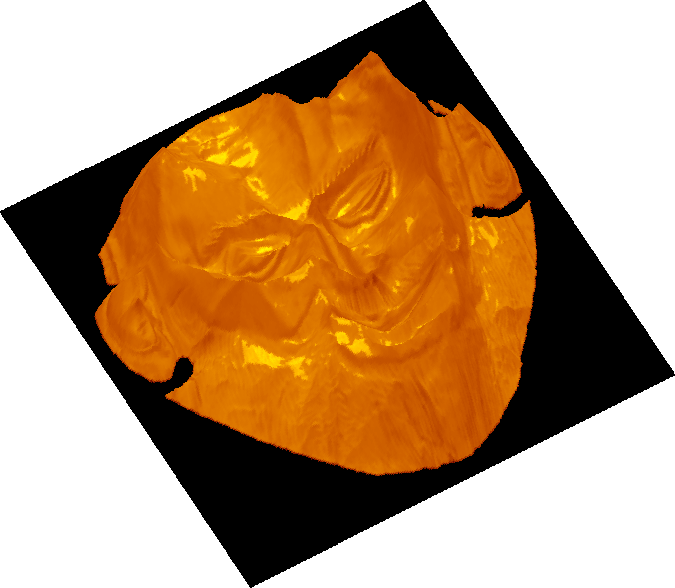} \\ 
\includegraphics[width=0.255\textwidth]{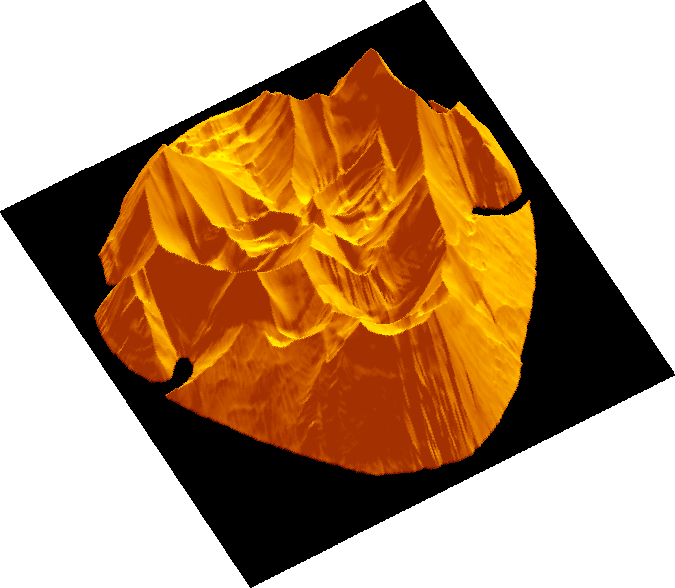} &
\includegraphics[width=0.255\textwidth]{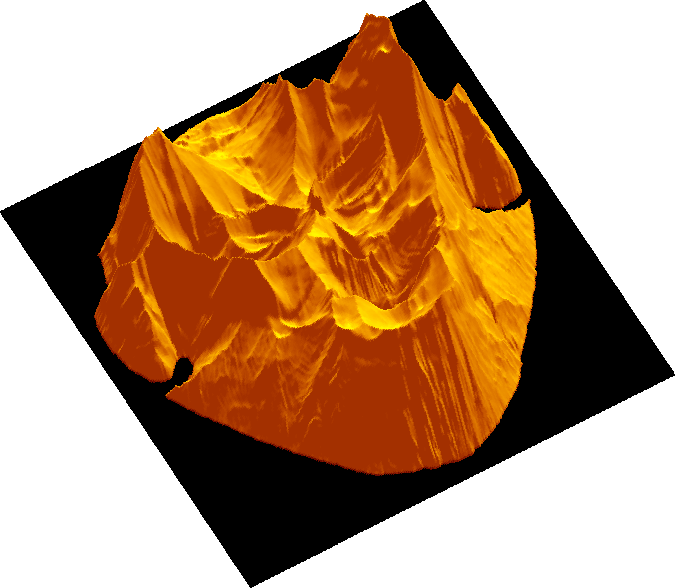} &
\includegraphics[width=0.255\textwidth]{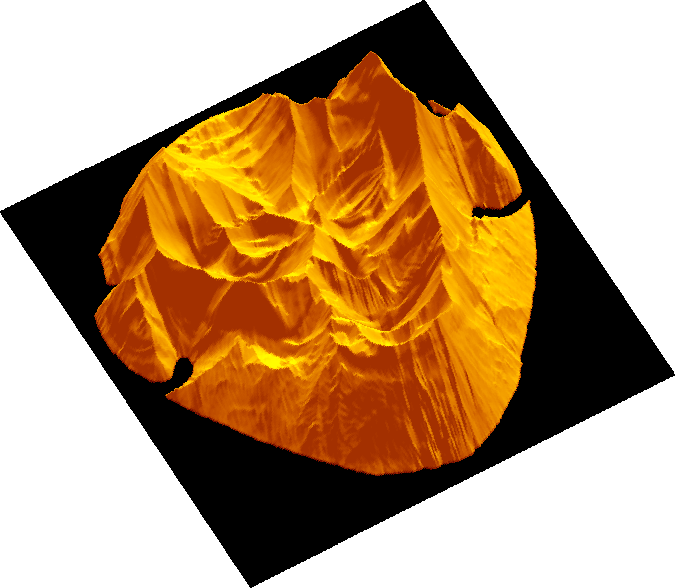} \\ 
  \hline 
\end{tabular}
 \caption{Agamennon mask: results with oblique light source $\vomega_{obl} = (0, 0.0995, 0.9950)$. On the first row the output images, on the second row the 3D reconstruction  with vertical view, on the third row the 3D reconstruction with oblique view.}  \label{fig:Agamennon_obl_results}
\end{figure}

\begin{figure}[h!]
\centering
 \subfigure[Real Vase Input]
   {\includegraphics[width=0.32\textwidth]{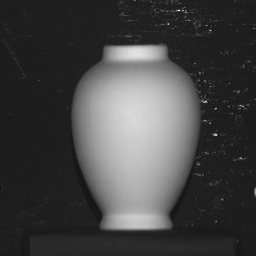} \label{fig:input_RV}  }
 \hspace{2mm}
  \subfigure[Real Vase Mask]
  {\includegraphics[width=0.32\textwidth]{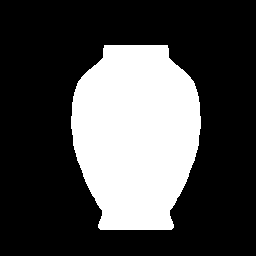} \label{fig:mask_RV}  }
 \caption{Real Vase images (size $256\times256$):
            (a)~Input image; 
            (b)~Mask.}
            
 \label{fig:RV_images}
 \end{figure}

\begin{table} [h!]
     \caption{Real Agamemnon mask:  $L^2$, $L^{\infty}$ errors with oblique light source $\vomega_{obl} = (0, 0.0995, 0.9950)$.}   \label{tab:Agam_errors_obl} 
\begin{tabular}{p{2.1cm}p{2cm}p{2cm}}
\hline\noalign{\smallskip}
  SL--Schemes & $L^2(I)$ & $L^\infty(I)$  \\
  \hline 
LAM     & 0.0585 & 0.4863  \\
ON-01  & 0.0663 & 0.4588  \\
ON-02  & 0.0670 & 0.4471  \\
ON-03  & 0.0708 & 0.5451   \\
PH-s02 & 0.1141 & 0.5725   \\
PH-s03  & 0.1580 & 0.6196 \\
PH-s04   & 0.2063 & 0.6706 \\
  \hline 
\end{tabular}
\end{table}
In Fig. \ref{fig:Agamennon_obl_results} we can see the output image and the 3D reconstruction in a single case for each models. What we can note is that also using more realistic illumination models as the two non-Lambertian considered, we do not obtain a so better approximation of the ground truth solution. This is due to the missing important informations (e.g. the correct oblique light source direction, the values of the parameter involved that are known for the models but not available in the real cases) and also due to the fact that we are using a first order method of approximation.

\begin{figure}[h!]
\centering
   \begin{tabular}{p{3.6cm}p{3.8cm}}
   \hline\noalign{\smallskip}
  \hspace{1.2cm} \textbf{Out}  &  \hspace{0.cm} \textbf{3D reconstruction} \\ 
  \hline \\
\includegraphics[width=0.25\textwidth]{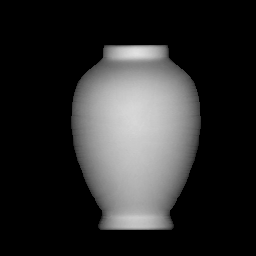} &
\includegraphics[width=0.255\textwidth]{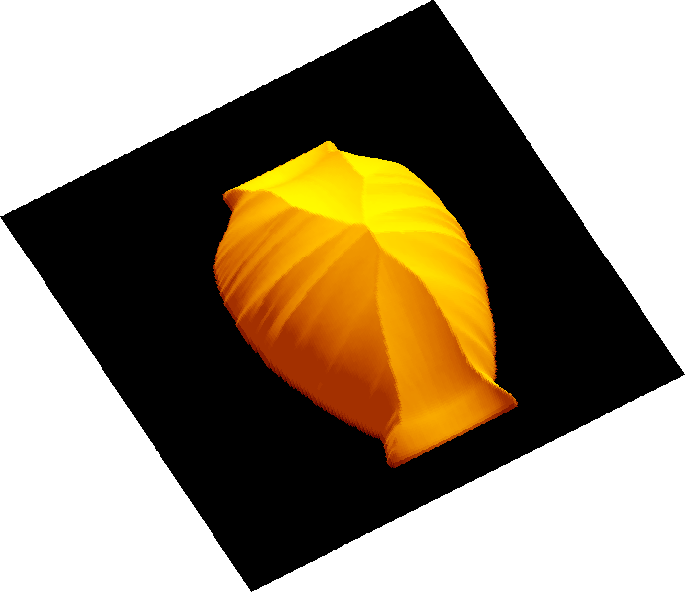} \\
\includegraphics[width=0.25\textwidth]{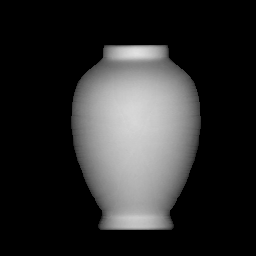} &
\includegraphics[width=0.255\textwidth]{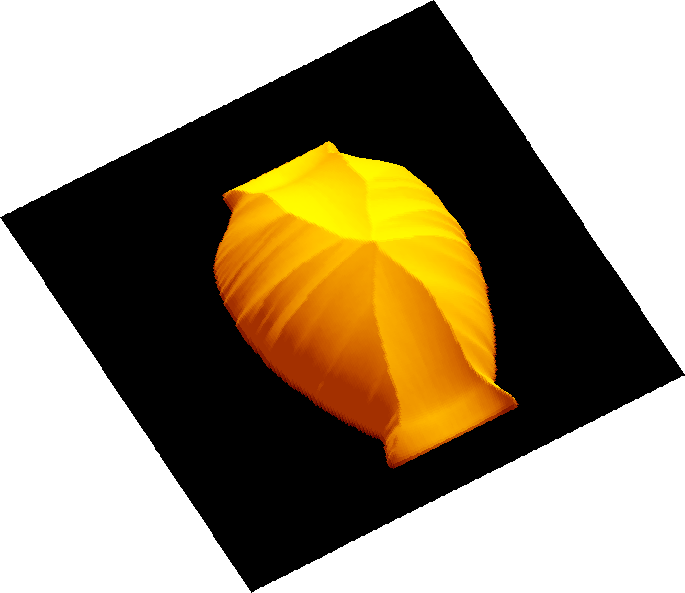} \\
\includegraphics[width=0.25\textwidth]{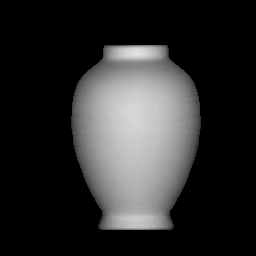} &
\includegraphics[width=0.255\textwidth]{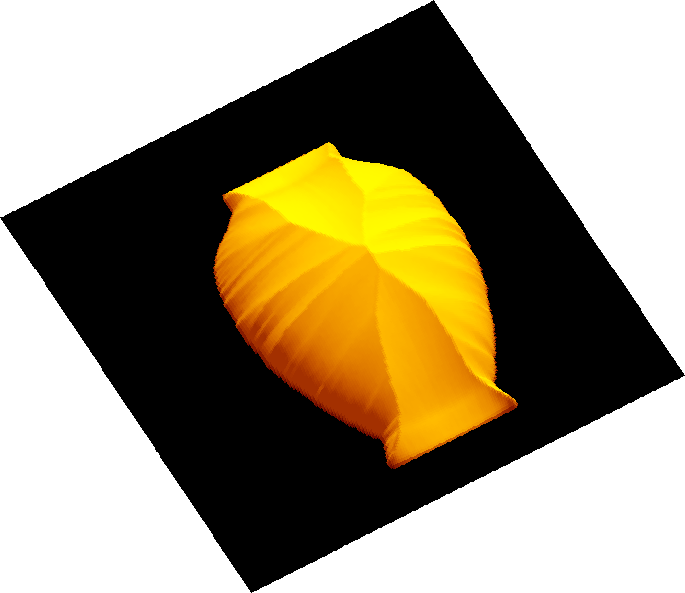} \\
  \hline 
\end{tabular}
 \caption{Real Vase: Output images and 3D reconstructions. On the first row the L--model, on the second row the ON--model with $\sigma = 0.2$, on the third row the PH--model with $k_S = 0.2$.}
   \label{fig:RV_reconstruction}
\end{figure}

\paragraph{Test 8: Real Vase.}
With this test we want to investigate the stability of our method with respect to the presence of noise. In fact, looking at the Fig. \ref{fig:input_RV} we can consider that RV is a noisy version of the synthetic vase used in the Tests 5 and 6.
 The test was performed using a vertical light source $\omega = (0,0,1)$. The output images, computed a posteriori by using the gradient of $u$ approximated via centered finite differences starting from the values of $u$ just computed by the numerical scheme, are visible in Fig. \ref{fig:RV_reconstruction}, first column. The reconstruction obtained with the three models are visible in the same Fig. \ref{fig:RV_reconstruction}, second column. What we can note is that all the reconstructions, obtained using homogeneous Dirichlet BC, suffer for a concave/convex ambiguity, as already noted for the synthetic vase (SV in the following).
Since it is visible looking at the SV and the RV tests that we obtain results with errors of the same order of magnitude around $10^{-2}$, 
considering that RV is a noisy version of SV, this shows that the method is stable in the presence of noise in the image.
\begin{table} [h!]
     \caption{Real Vase:  $L^2$, $L^{\infty}$ errors with vertical light source $\vomega = (0, 0, 1)$.}
        \label{tab:RV_errors} 
     \begin{tabular}{p{2.1cm}p{2cm}p{2cm}}
\hline\noalign{\smallskip}
  SL--Schemes & $L^2(I)$ & $L^\infty(I)$  \\
  \hline 
LAM     & 0.0093 &  0.0784   \\
ON-02  & 0.0094 &  0.0784  \\
ON-04  &  0.0111 &  0.0824  \\
PH-s02 &  0.0095 &  0.0824 \\
PH-s04   & 0.0098 &  0.0824 \\
  \hline 
\end{tabular}
\end{table}

\begin{figure}[h!]
\centering
   \begin{tabular}{p{3.4cm}p{3.4cm}p{3.4cm}}
   \hline\noalign{\smallskip}
  \hspace{0.3cm} \textbf{ON $\sigma = 0.2$}  &  \hspace{0.3cm} \textbf{ON $\sigma = 0.4$} &  \hspace{0.3cm} \textbf{ON $\sigma = 0.6$} \\ 
  \hline \\
\includegraphics[width=0.25\textwidth]{SHONAsigma02-SL-ortografico-RV.png} &
\includegraphics[width=0.25\textwidth]{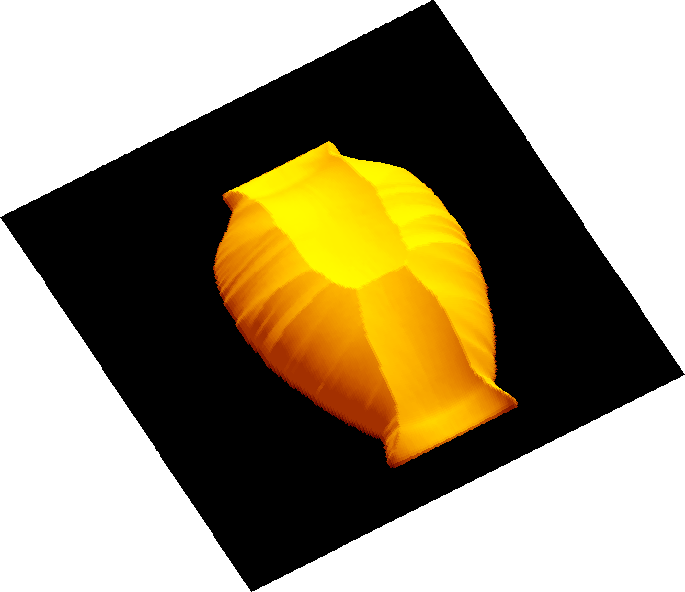} &
\includegraphics[width=0.25\textwidth]{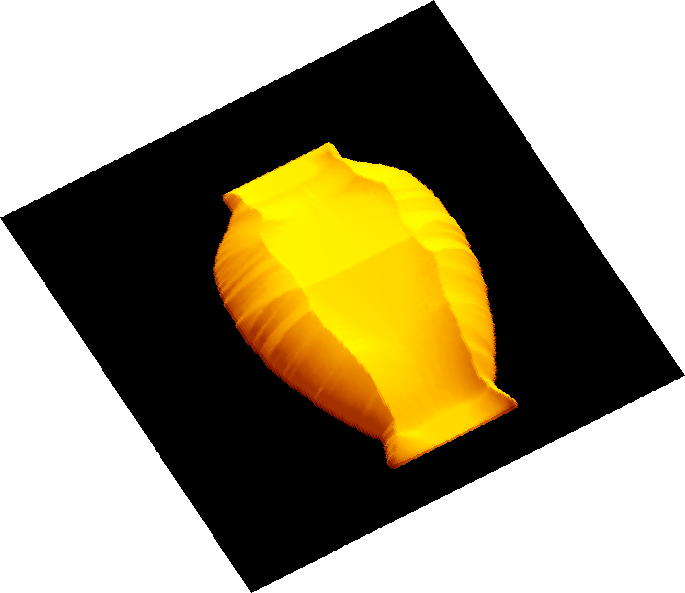} \\
  \hline 
\end{tabular}
 \caption{Real Vase: 3D reconstructions related to the ON--model, varying $\sigma$. From left to right $\sigma = 0.2, 0.4, 0.6$.}
   \label{fig:RV_recon_ON}
\end{figure}

Finally, in Fig. \ref{fig:RV_recon_ON} one can note the behavior of the ON--model by varying the value of the parameter $\sigma$. 
Since in real situations we do not know it (as other parameters like the light source direction) we can only vary it in order to see the one which gives the best fit with the image.  
Looking at Fig. \ref{fig:RV_recon_ON} what we can observe is that increasing the value of $\sigma$ the reconstruction shows a wider concave/convex ambiguity, which affects more pixels. But this holds only in this specific case, not in general for all the real images as  a rule.

\begin{figure}[h!]
\centering
\includegraphics[width=0.52\textwidth]{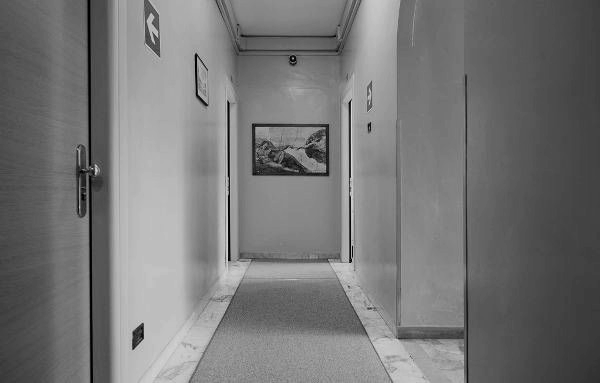}
 \caption{Image of a real scene  (size $600\times383$).}  \label{fig:corridor_input}
\end{figure}

\begin{figure}[h!]
\centering
   \begin{tabular}{p{4.6cm}p{4.8cm}}
   \hline\noalign{\smallskip}
  \hspace{1.2cm} \textbf{Out}  &  \hspace{0.cm} \textbf{3D reconstruction} \\ 
  \hline \\
\includegraphics[width=0.32\textwidth]{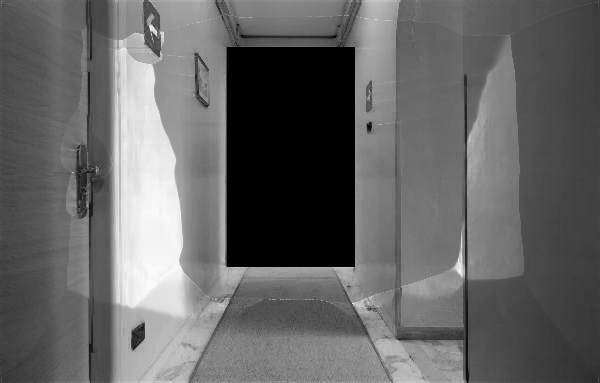} &
\includegraphics[width=0.32\textwidth]{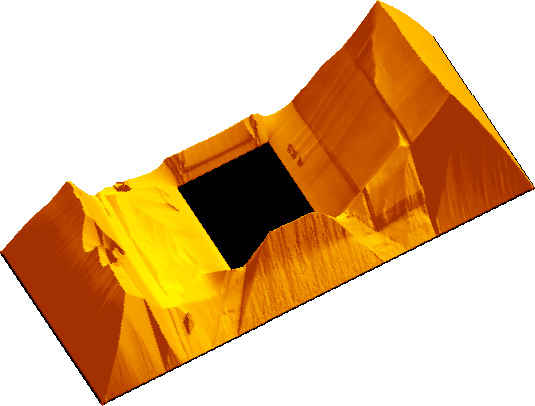} \\
\includegraphics[width=0.32\textwidth]{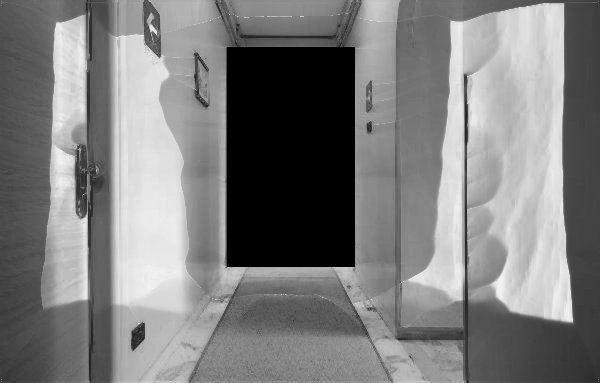} &
\includegraphics[width=0.32\textwidth]{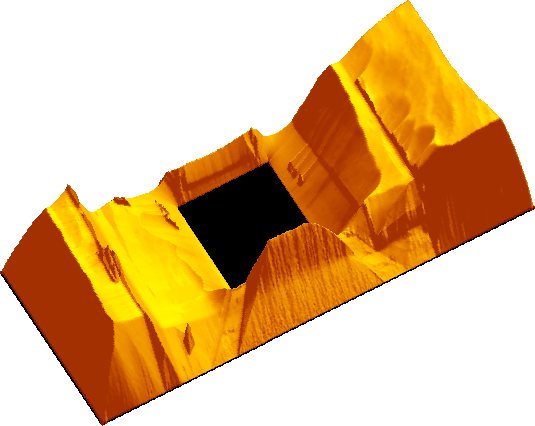} \\
\includegraphics[width=0.32\textwidth]{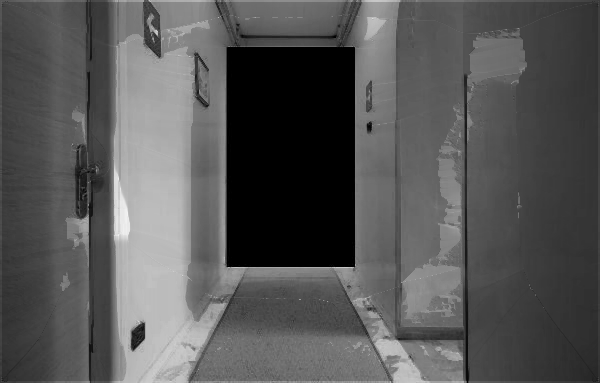} & 
\includegraphics[width=0.32\textwidth]{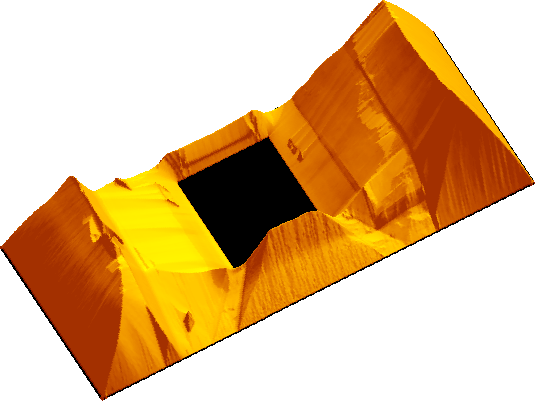} \\ 
  \hline 
\end{tabular}
 \caption{Output images and 3D reconstructions of a scene for a robot path planning application. 
 On the first row the L--model, on the second row the ON--model with $\sigma = 0.1$, on the third row the PH--model with $k_S = 0.2$.}  \label{fig:corridor}
\end{figure}

\paragraph{Test 9: Corridor.}
As illustrative example, let us consider a real image of a corridor (see Fig. \ref{fig:corridor_input}) as a typical example of a scene which can be useful for a robot navigation problem. 
The test has been added as illustrative example to show that even for a real scene which does not satisfies all the assumptions and for which several informations are missing
(e.g. boundary conditions) the method is able to compute a reasonable accurate reconstruction in the central part of the corridor (clearly the boundaries are wrong due to a lack of information). 
The size of this image is $600\times383$. Note that for this picture we do not know the parameters and the light direction in the scene. It seems that there is a diffused light and more than one light source. So this picture does not satisfy many assumptions we used in the theoretical part. 
In order to apply our numerical scheme we considered a Dirichlet boundary condition equal to zero at the wall located at the bottom of the corridor. In this way we have a better perception of the depth of the scene. 
In Fig. \ref{fig:corridor} we can see the output images (on the first column) and the 3D reconstructions (on the second column) obtained by L--model, ON--model with $\sigma = 0.1$ and PH--model with $k_S = 0.2$. 
In this example the PH--model seems to recognize the scene better than the ON--model. In some sense this is probably due to the fact that it has less parameters
so it is easier to tune to a real situation where the information on the parameters is not available.
We point out that this test is just an illustration of the fact that coupling SfS with additional informations (e.g. coming from distance sensors to fix boundary conditions) can be useful to describe a scene.

\begin{figure}[h!]
\centering
   \begin{tabular}{p{2.9cm}p{2.9cm}}
   \hline\noalign{\smallskip}
  \hspace{0.3cm} \textbf{Input}  &  \hspace{-0.35cm} \textbf{reconstruction} \\ 
  \hline \\
\hspace{-0.1cm}\includegraphics[width=0.24\textwidth]{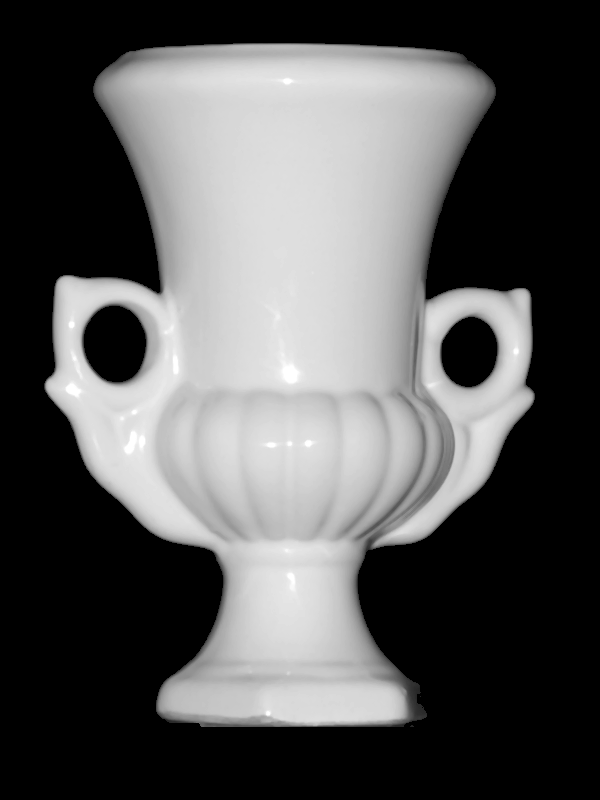} &\hspace{-0.2cm}
\hspace{-0.1cm}\includegraphics[width=0.27\textwidth]{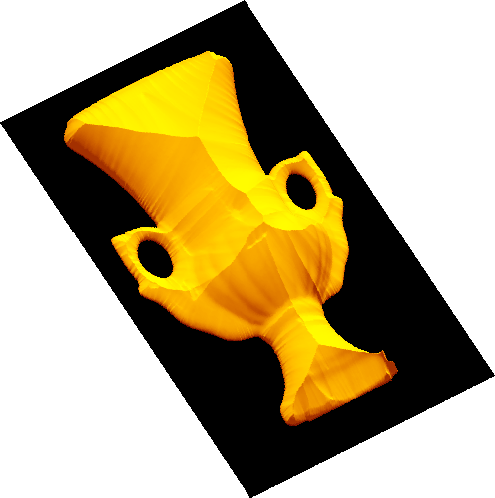} \\ 
\hspace{-0.1cm}\includegraphics[width=0.24\textwidth]{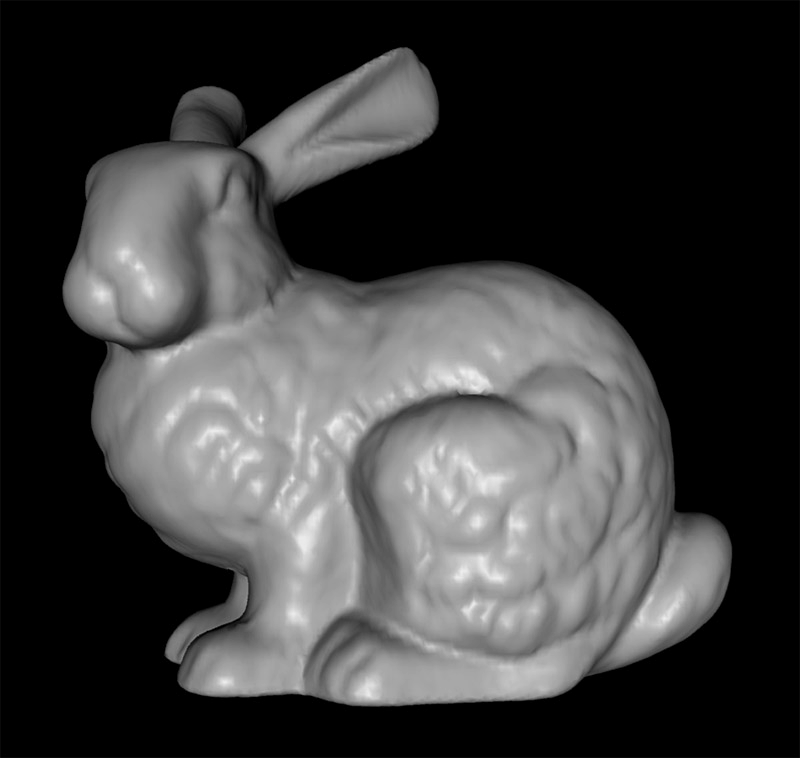} &\hspace{-0.2cm}
\vspace{-0.2cm} \hspace{-0.1cm}\includegraphics[width=0.27\textwidth]{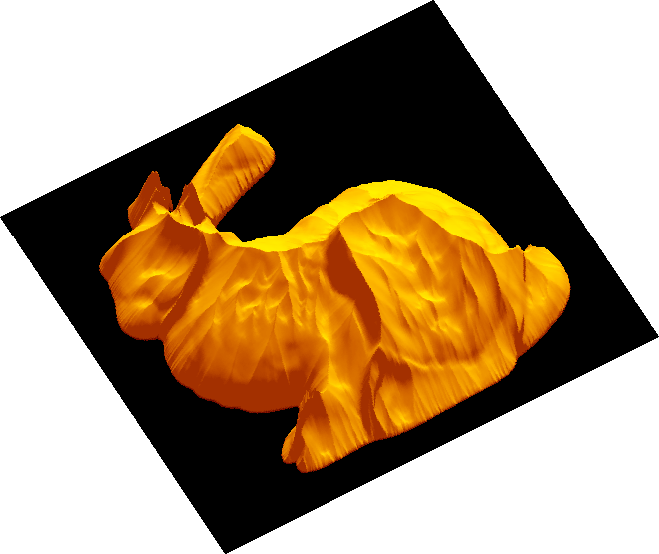} \\ 
\hspace{-0.1cm}\includegraphics[width=0.24\textwidth]{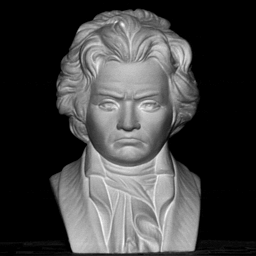} &\hspace{-0.2cm}
\vspace{-0.2cm} \hspace{-0.1cm}\includegraphics[width=0.27\textwidth]{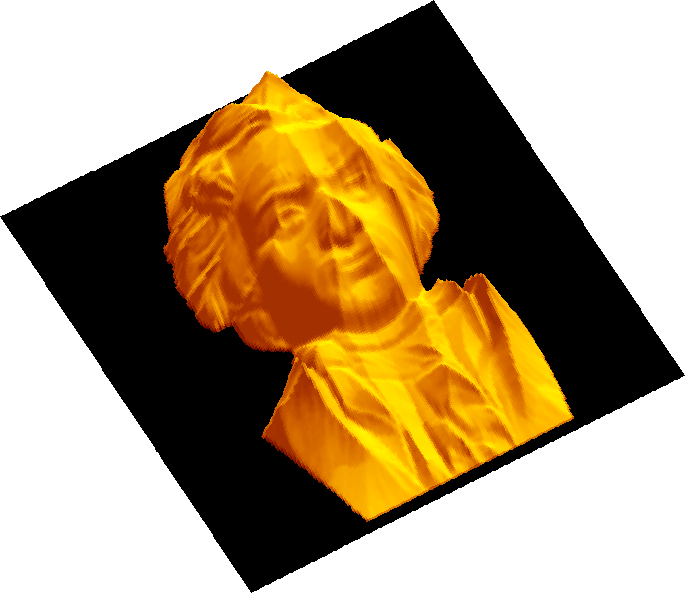} \\
  \hline 
\end{tabular}
 \caption{Experiments on real images: input images and the recovered shapes obtained by our approach. 
 On the first row: urn reconstructed by the L--model with $\vomega = (0,0,1)$. 
 On the second row: rabbit reconstructed by the ON--model with $\sigma = 0.2$, $\vomega = (0,0,1)$, $\vV = (0,0,1)$. 
On the third row: Beethoven reconstructed by the PH--model with $\vomega = (0.0168, 1.198, 0.9801)$, $k_S = 0.2$.}
 \label{fig:real_tests} 
\end{figure}

\paragraph{Test 10: Other tests on real images.}
In order to demonstrate the applicability of our proposed approach for real data, we added here other experiments conducted on a real urn, real rabbit and real Beethoven's bust. The input images and the related recovered shapes obtained by the three models studied under different light directions and parameters are shown in Fig. \ref{fig:real_tests}. As visible from Fig. \ref{fig:real_tests}, the results are quite good, even for pictures like the rabbit or the bust of Beethoven, which have many details, even if they still suffer for the concave/convex ambiguity typical of the SfS problem. 

\section{Conclusions}
In this paper we derived nonlinear partial differential equations of first order, i.e. Hamilton-Jacobi equations, associated to the non-Lambertian reflectance models ON--model and PH--model. We have obtained the model equations for all the possible cases, coupling vertical or oblique light source with vertical or oblique position of the observer.  This exhaustive  description has shown that these models lead to stationary Hamilton-Jacobi equations with a same structure and this allows for a unified mathematical formulation in terms of fixed point problem. This general formulation is interesting because we can switch on and off the different terms related to ambient, diffuse and specular reflection in a very simple way. As a result, this general model is very flexible to treat the various situations with vertical and oblique light sources.
Unfortunately, we have observed that none of these models is able to overcome the typical concave/convex ambiguity known for the classical Lambertian model.
Despite this limitations, the approach presented in this paper is able to improve the Lambertian model that is not suitable to deal with realistic images coming from medical or security application. The numerical methods presented here can be applied to solve the equations corresponding to the ON--model and the PH--model in all the situations although they suffer for the presence of several parameters related to the roughness or to the specular effect. 
Nonetheless, looking at the comparisons with other methods or models shown in this paper, one can see that the Semi-Lagrangian approach is competitive with respect to other techniques used, both in terms of CPU time and accuracy. 
As we have seen, in the complex nonlinear PDEs associated to non-Lambertian models the  parameters play a crucial role to obtain accurate results. In fact, varying the value of the parameters it is possible to improve the approximation with respect to the classical L--model. We can also say that for real images  the PH--model seems easier to tune, perhaps because we need to manage less parameters. \\
Focusing the attention on the tests performed with an oblique light source, we have to do some comments that are common to the PH--model and the ON--model. 
Several terms appear in these models and each of them gives a contribution to the roundoff error. Note that the accumulation of these roundoff errors makes difficult in the oblique case to obtain a great accuracy. A possible improvement could be the use of second order schemes, that release the link between the space and the time steps which characterizes and limits the accuracy for first order schemes. Another interesting direction would be to extend this formulation to other reflectance models (like e.g. the Ward's model) and/or to consider perspective projection also including an attenuation term which can help to resolve the concave/convex ambiguity or considering more than one input image with non-Lambertian models, that solves the well-known ambiguity (see \cite{TMDDb15} for a first step in this last direction considering the Blinn-Phong model).  These directions will be explored in future works. 

\subsection*{Acknowledgements}
The first author wishes to acknowledge the support obtained by the INDAM under the GNCS research project "Metodi ad alta risoluzione per problemi evolutivi fortemente nonlineari" 
and the hospitality obtained by Prof. Edwin Hancock of the University of York, Department of Computer Science, during her scholarship ``Borsa di perfezionamento all'estero" paid by Sapienza - Università di Roma.

\bibliographystyle{plain}

\begin{thebibliography}{10}

\bibitem{AF07}
A.~H. Ahmed and A.~A. Farag.
\newblock {S}hape from {S}hading under various imaging conditions.
\newblock In {\em IEEE International Conference on Computer Vision and Pattern
  Recognition CVPR'07}, pages X1--X8, Minneapolis, MN, June 2007.

\bibitem{AF06}
A.H. Ahmed and A.A. Farag.
\newblock A new formulation for {S}hape from {S}hading for non-{L}ambertian
  surfaces.
\newblock In {\em Proc. IEEE Conference on Computer Vision and Pattern
  Recognition (CVPR)}, pages 1817--1824, 2006.

\bibitem{BY94}
S.~Bakshi and Y.-H. Yang.
\newblock Shape-from-{S}hading for non-{L}ambertian surfaces.
\newblock In {\em Proc. IEEE International Conference on Image Processing
  (ICIP)}, pages 130--134, 1994.

\bibitem{Barles94}
G.~Barles.
\newblock {\em Solutions de viscosité des équations de {H}amilton-{J}acobi}.
\newblock Springer Verlag, 1994.

\bibitem{BAC09}
S.~Biswas, G.~Aggarwal, and R.~Chellappa.
\newblock Robust estimation of albedo for illumination-invariant matching and
  shape recovery.
\newblock {\em Pattern Analysis and Machine Intelligence, IEEE Transactions
  on}, 31(5):884--899, May 2009.

\bibitem{BCK92}
M.J. Brooks, W.~Chojnacki, and R.~Kozera.
\newblock Impossible and ambiguous shading patterns.
\newblock {\em International Journal of Computer Vision}, 7(2):119--126, 1992.

\bibitem{BH85}
M.J. Brooks and B.K.P. Horn.
\newblock Shape and source from shading.
\newblock In {\em Proceedings of the Ninth International Joint Conference on
  Artificial Intelligence (vol. II)}, pages 932--936, Los Angeles, CA, USA,
  1985.

\bibitem{Bruss81}
A.~Bruss.
\newblock The image irradiance equation: Its solution and application.
\newblock Technical report ai-tr-623, Massachusetts Institute of Technology,
  June 1981.

\bibitem{CF96}
F.~Camilli and M.~Falcone.
\newblock An approximation scheme for the maximal solution of the {S}hape from
  {S}hading model.
\newblock In {\em Proceedings of the IEEE International Conference on Image
  Processing (vol. I)}, volume~1, pages 49--52, Lausanne, Switzerland, Sep
  1996.

\bibitem{CG00}
F.~Camilli and L.~Gr\"une.
\newblock Numerical approximation of the maximal solutions for a class of
  degenerate {H}amilton-{J}acobi equations.
\newblock {\em SIAM Journal on Numerical Analysis}, 38(5):1540--1560, 2000.

\bibitem{CS99}
F.~Camilli and A.~Siconolfi.
\newblock Maximal subsolutions of a class of degenerate {H}amilton-{J}acobi
  problems.
\newblock {\em Indiana University Mathematics Journal}, 48(3):1111--1131, 1999.

\bibitem{CFF13}
E.~Carlini, M.~Falcone, and A.~Festa.
\newblock A brief survey on semi-lagrangian schemes for image processing.
\newblock In Michael Breuss, Alfred Bruckstein, and Petros Maragos, editors,
  {\em Innovations for Shape Analysis: Models and Algorithms}, pages 191--218.
  Ed. Springer, 2013.
\newblock ISBN: 978-3-642-34140-3.

\bibitem{Cha94}
A.~Chambolle.
\newblock A uniqueness result in the theory of stereo vision: coupling {S}hape
  from {S}hading and binocular information allows unambiguous depth
  reconstruction.
\newblock {\em Annales de l'Istitute Henri Poincarè}, 11(1):1--16, 1994.

\bibitem{Courteille_ICPR2004}
F.~Courteille, A.~Crouzil, J.-D. Durou, and P.~Gurdjos.
\newblock Towards {S}hape from {S}hading under realistic photographic
  conditions.
\newblock In {\em Proc. of the 17th International Conference on Pattern
  Recognition}, volume~2, pages 277--280. Cambridge, UK, 2004.

\bibitem{CCDG07}
F.~Courteille, A.~Crouzil, J.-D. Durou, and P.~Gurdjos.
\newblock Shape from shading for the digitization of curved documents.
\newblock {\em Machine Vision and Applications}, 18(5):301--316, 2007.

\bibitem{DD00}
P.~Daniel and J.-D. Durou.
\newblock From deterministic to stochastic methods for {S}hape from {S}hading.
\newblock In {\em Proceedings of the Fourth Asian Conference on Computer
  Vision}, pages 187--192, Taipei, Taiwan, 2000.

\bibitem{JVD51}
J.~Van Diggelen.
\newblock A photometric investigation of the slopes and heights of the ranges
  of hills in the maria of the moon.
\newblock {\em Bulletin of the Astronomical Institute of the Netherlands},
  11(423):283--290, 1951.

\bibitem{DurouEMMCVPR2009}
J.-D. Durou, J.-F. Aujol, and F.~Courteille.
\newblock Integrating the normal field of a surface in the presence of
  discontinuities.
\newblock In {\em Energy Minimization Methods in Computer Vision and Pattern
  Recognition (EMMCVPR)}, volume 5681, pages 261--273, 2009.

\bibitem{DFS08}
J.-D. Durou, M.~Falcone, and M.~Sagona.
\newblock Numerical methods for {S}hape from {S}hading: a new survey with
  benchmarks computer vision and image understanding.
\newblock {\em Elsevier}, 109(1):22--43, 2008.

\bibitem{Falcone97}
M.~Falcone.
\newblock Numerical solution of {D}ynamic {P}rogramming equations.
\newblock {\em In Bardi, M., Dolcetta, I.C., Optimal control and viscosity
  solutions of Hamilton-Jacobi-Bellman equations}, pages 471--504, 1997.

\bibitem{FF14}
M.~Falcone and R.~Ferretti.
\newblock {\em Semi-Lagrangian Approximation Schemes for Linear and
  {H}amilton-{J}acobi Equations}.
\newblock SIAM, 2014.

\bibitem{FS97}
M.~Falcone and M.~Sagona.
\newblock An algorithm for the global solution of the {S}hape from {S}hading
  model.
\newblock In {\em Proceedings of the Ninth International Conference on Image
  Analysis and Processing (vol. I), vol. 1310 of Lecture Notes in Computer
  Science}, pages 596--603, Florence, Italy, 1997.

\bibitem{FSS03}
M.~Falcone, M.~Sagona, and A.~Seghini.
\newblock A scheme for the {S}hape--from--{S}hading model with “black shadows”.
\newblock {\em Numerical Mathematics and Advanced Applications}, pages
  503--512, 2003.

\bibitem{FDSB04}
H.~Fassold, R.~Danzl, K.~Schindler, and H.~Bischof.
\newblock Reconstruction of archaeological finds using shape from stereo and
  shape from shading.
\newblock In {\em Proc. 9th Computer Vision Winter Workshop, Piran, Slovenia},
  pages 21--30, 2004.

\bibitem{FP12}
P.~Favaro and T.~Papadhimitri.
\newblock A closed-form solution to uncalibrated photometric stereo via diffuse
  maxima.
\newblock In {\em Proceedings of the IEEE Conference on Computer Vision and
  Pattern Recognition}, pages 821--828, 2012.

\bibitem{FC88}
R.T. Frankot and R.~Chellappa.
\newblock A method for enforcing integrability in {S}hape from {S}hading
  algorithms.
\newblock {\em IEEE Transactions on Pattern Analysis and Machine Intelligence},
  10(4):439--451, 1988.

\bibitem{Grumpe2014}
A.~Grumpe, F.~Belkhir, and C.~W{\"o}hler.
\newblock Construction of lunar {DEM}s based on reflectance modelling.
\newblock {\em Advances in Space Research}, 53(12):1735 --1767, 2014.

\bibitem{Horn_PhD1970}
B.K.P. Horn.
\newblock {\em {S}hape from {S}hading: A Method for Obtaining the Shape of a
  Smooth Opaque Object From One View}.
\newblock PhD thesis, Massachusetts Institute of Technology, 1970.

\bibitem{H75}
B.K.P. Horn.
\newblock Obtaining {S}hape from {S}hading information.
\newblock In P.H.~Winston (Ed.), editor, {\em The Psychology of Computer
  Vision}, chapter~4, pages 115--155. McGraw-Hill, 1975.

\bibitem{HB86}
B.K.P. Horn and M.J. Brooks.
\newblock The variational approach to {S}hape from {S}hading.
\newblock {\em Computer Vision, Graphics and Image Processing}, 33(2):174--208,
  1986.

\bibitem{HB89}
B.K.P. Horn and M.J. Brooks.
\newblock {\em {S}hape from {S}hading (Artificial Intelligence)}.
\newblock The MIT Press, 1989.

\bibitem{IH81}
Katsushi Ikeuchi and Berthold~K.P. Horn.
\newblock Numerical {S}hape from {S}hading and occluding boundaries.
\newblock {\em Artificial Intelligence}, 17(1--3):141--184, 1981.

\bibitem{IR95}
H.~Ishii and M.~Ramaswamy.
\newblock Uniqueness results for a class of {H}amilton-{J}acobi equations with
  singular coefficients.
\newblock {\em Communications in Partial Differential Equations},
  20:2187--2213, 1995.

\bibitem{JBBG12}
{Y.-C}. Ju, M.~Breu{\ss}, A.~Bruhn, and S.~Galliani.
\newblock {S}hape from {S}hading for rough surfaces: Analysis of the
  {O}ren-{N}ayar model.
\newblock In {\em Proc. British Machine Vision Conference (BMVC)}, pages
  104.1--104.11. BMVA Press, 2012.

\bibitem{JTBBK13}
{Y.-C}. Ju, S.~Tozza, M.~Breu\ss, A.~Bruhn, and A.~Kleefeld.
\newblock Generalised {P}erspective {S}hape from {S}hading with {O}ren-{N}ayar
  {R}eflectance.
\newblock In {\em Proceedings of the 24th British Machine Vision Conference
  (BMVC 2013)}, pages 42.1--42.11, Bristol, UK, 2013. BMVA Press.

\bibitem{KO01}
J.~Kain and D.N. Ostrov.
\newblock Numerical {S}hape from {S}hading for discontinuous photographic
  images.
\newblock {\em International Journal of Computer Vision}, 44(3):163--173, 2001.

\bibitem{Kozera91}
R.~Kozera.
\newblock Existence and uniqueness in photometric stereo.
\newblock {\em Appl. Math. Comput.}, 44(1):103, 1991.

\bibitem{Kozera97}
R.~Kozera.
\newblock Uniqueness in {S}hape from {S}hading revisited.
\newblock {\em Journal of Mathematical Imaging and Vision}, 7(2):123--138,
  1997.

\bibitem{Kru75}
S.~N. Kruzkov.
\newblock The generalized solution of the hamilton-jacobi equations of eikonal
  type i.
\newblock {\em Math. USSR Sbornik}, 27:406--446, 1975.

\bibitem{LRT93}
P.L. Lions, E.~Rouy, and A.~Tourin.
\newblock {S}hape-from-{S}hading, viscosity solutions and edges.
\newblock {\em Numerische Mathematik}, 64(3):323--353, 1993.

\bibitem{Lohse2006}
V.~Lohse, C.~Heipke, and R.~L. Kirk.
\newblock Derivation of planetary topography using multi-image
  shape-from-shading.
\newblock {\em Planetary and Space Science}, 54(7):661--674, 2006.

\bibitem{MF13}
R.~Mecca and M.~Falcone.
\newblock Uniqueness and approximation of a photometric {S}hape-from-{S}hading
  model.
\newblock {\em SIAM J. Imaging Sciences}, 6(1):616--659, 2013.

\bibitem{MT13}
R.~Mecca and S.~Tozza.
\newblock Shape reconstruction of symmetric surfaces using photometric stereo.
\newblock In Michael Breuss, Alfred Bruckstein, and Petros Maragos, editors,
  {\em Innovations for Shape Analysis: Models and Algorithms}, pages 219--243.
  Ed. Springer, 2013.
\newblock ISBN: 978-3-642-34140-3.

\bibitem{OD97}
T.~Okatani and K.~Deguchi.
\newblock Shape reconstruction from an endoscope image by {S}hape from
  {S}hading technique for a point light source at the projection center.
\newblock {\em Computer Vision and Image Understanding}, 66(2):119--131, 1997.

\bibitem{Oliensis91b}
J.~Oliensis.
\newblock {S}hape from {S}hading as a partially well-constrained problem.
\newblock {\em Computer Vision, Graphics and Image Processing: Image
  Understanding}, 54(2):163--183, 1991.

\bibitem{Oliensis91a}
J.~Oliensis.
\newblock Uniqueness in {S}hape from {S}hading.
\newblock {\em International Journal of Computer Vision}, 6(2):75--104, 1991.

\bibitem{ON_CVPR93}
M.~Oren and S.K. Nayar.
\newblock Diffuse reflectance from rough surfaces.
\newblock In {\em Proc. IEEE Conference on Computer Vision and Pattern
  Recognition (CVPR)}, pages 763--764, 1993.

\bibitem{ON_SIGGRAPH94}
M.~Oren and S.K. Nayar.
\newblock Generalization of {L}ambert's reflectance model.
\newblock In {\em Proc. International Conference and Exhibition on Computer
  Graphics and Interactive Techniques (SIGGRAPH)}, pages 239--246, 1994.

\bibitem{ON_ECCV94}
M.~Oren and S.K. Nayar.
\newblock Seeing beyond {L}ambert's law.
\newblock In {\em Proc. European Conference on Computer Vision (ECCV)}, pages
  269--280, 1994.

\bibitem{ON_IJCV95}
M.~Oren and S.K. Nayar.
\newblock Generalization of the {L}ambertian model and implications for machine
  vision.
\newblock {\em International Journal of Computer Vision}, 14(3):227--251, 1995.

\bibitem{Phong75}
B.~T. Phong.
\newblock Illumination for computer generated pictures.
\newblock {\em Communications of the ACM}, 18(6):311--317, 1975.

\bibitem{PCF06}
E.~Prados, F.~Camilli, and O.~Faugeras.
\newblock A viscosity solution method for {S}hape-from-{S}hading without image
  boundary data.
\newblock {\em ESAIM: Mathematical Modelling and Numerical Analysis},
  40(2):393--412, 2006.

\bibitem{PF03}
E.~Prados and O.~Faugeras.
\newblock A mathematical and algorithmic study of the lambertian sfs problem
  for orthographic and pinhole cameras.
\newblock Rapport de recherche 5005, Institut National de Recherche en
  Informatique et en Automatique, Sophia Antipolis, France, November 2003.

\bibitem{PF_ICCV2003}
E.~Prados and O.~Faugeras.
\newblock Perspective {S}hape from {S}hading and viscosity solutions.
\newblock In {\em Proc. IEEE International Conference on Computer Vision
  (ICCV)}, pages 826--831, 2003.

\bibitem{PF_CVPR2005}
E.~Prados and O.~Faugeras.
\newblock {S}hape from {S}hading: a well-posed problem?
\newblock In {\em Proc. IEEE Conference on Computer Vision and Pattern
  Recognition (CVPR)}, volume~2, pages 870--877, 2005.

\bibitem{PFR}
E.~Prados, O.~Faugeras, and E.~Rouy.
\newblock {S}hape from {S}hading and viscosity solutions.
\newblock In Anders Heyden, Gunnar Sparr, Mads Nielsen, and Peter Johansen,
  editors, {\em Computer Vision ‚ ECCV 2002}, volume 2351 of {\em Lecture Notes
  in Computer Science}, pages 790--804. Springer Berlin Heidelberg, 2002.

\bibitem{QLD15}
Y.~Qu{\'{e}}au, F.~Lauze, and J.{-}D. Durou.
\newblock Solving uncalibrated photometric stereo using total variation.
\newblock {\em Journal of Mathematical Imaging and Vision}, 52(1):87--107,
  2015.

\bibitem{Ragheb_3DPVT2004}
H.~Ragheb and E.R. Hancock.
\newblock Surface normals and height from non-lambertian image data.
\newblock In {\em Proc. International Symposium on 3D Data Processing,
  Visualization and Transmission (3DPVT)}, pages 18--25, 2004.

\bibitem{RH05}
H.~Ragheb and E.R. Hancock.
\newblock Surface radiance correction for {S}hape from {S}hading.
\newblock {\em Pattern Recognition}, 38(10):1574--1595, 2005.

\bibitem{Rindfleisch_PE1966}
T.~Rindfleisch.
\newblock Photometric method for lunar topography.
\newblock {\em Photogrammetric Engineering}, 32(2):262--277, 1966.

\bibitem{RT92}
E.~Rouy and A.~Tourin.
\newblock A viscosity solutions approach to {S}hape-from-{S}hading.
\newblock {\em SIAM Journal on Numerical Analysis}, 29(3):867--884, June 1992.

\bibitem{Sagona01}
M.~Sagona.
\newblock {\em Numerical methods for degenerate Eikonal type equations and
  applications}.
\newblock PhD thesis, Dipartimento di Matematica dell'Università di Napoli
  ""Federico II'', Napoli, Italy, November 2001.

\bibitem{Samaras_TPAMI2003}
D.~Samaras and D.~Metaxas.
\newblock Incorporating illumination constraints in deformable models for
  {S}hape from {S}hading and light direction estimation.
\newblock {\em IEEE Transactions on Pattern Analysis and Machine Intelligence},
  25(2):247--264, 2003.

\bibitem{SH05}
W.A.P. Smith and E.R. Hancock.
\newblock Recovering facial shape and albedo using a statistical model of
  surface normal direction.
\newblock In {\em 10th {IEEE} International Conference on Computer Vision
  {(ICCV} 2005), 17-20 October 2005, Beijing, China}, pages 588--595, 2005.

\bibitem{SH06}
W.A.P. Smith and E.R. Hancock.
\newblock Estimating facial albedo from a single image.
\newblock {\em {IJPRAI}}, 20(6):955--970, 2006.

\bibitem{SKWMM11}
W.~Steffen, N.~Koning, S.~Wenger, C.~Morisset, and M.~Magnor.
\newblock Shape: A 3d modeling tool for astrophysics.
\newblock {\em Visualization and Computer Graphics, IEEE Transactions on},
  17(4):454--465, April 2011.

\bibitem{Strat79}
T.M. Strat.
\newblock A numerical method for {S}hape from {S}hading from a single image.
\newblock Master's thesis, Department of Electrical Engineering and Computer
  Science, Massachussetts Institute of Technology, Cambridge, MA, USA, 1979.

\bibitem{Szeliski91}
R.~Szeliski.
\newblock Fast {S}hape from {S}hading.
\newblock {\em Computer Vision, Graphics and Image Processing: Image
  Understanding}, 53(2):129--153, 1991.

\bibitem{Tankus2005}
A.~Tankus and N.~Kiryati.
\newblock Photometric stereo under perspective projection.
\newblock In {\em Proc. IEEE International Conference on Computer Vision
  (ICCV)}, volume~1, pages 611--616, 2005.

\bibitem{TSY05}
A.~Tankus, N.~Sochen, and Y.~Yeshurun.
\newblock {S}hape-from-{S}hading {U}nder {P}erspective {P}rojection.
\newblock {\em International Journal of Computer Vision}, 63(1):21--43, 2005.

\bibitem{Tozza14}
S.~Tozza.
\newblock {\em Analysis and {A}pproximation of {N}on-{L}ambertian
  {S}hape-from-{S}hading {M}odels}.
\newblock PhD thesis, Dipartimento di Matematica, Sapienza - Università di
  Roma, Rome, Italy, January 2015.

\bibitem{ToFal_Dagstuhl}
S.~Tozza and M.~Falcone.
\newblock A comparison of non-lambertian models for the {S}hape--from-{S}hading
  problem.
\newblock In Michael Breuss, Alfred Bruckstein, Petros Maragos, and Stefanie
  Wuhrer, editors, {\em New Perspectives in Shape Analysis}. Ed. Springer.
\newblock To appear.

\bibitem{TF14}
S.~Tozza and M.~Falcone.
\newblock A {S}emi-{L}agrangian {A}pproximation of the {O}ren-{N}ayar {PDE} for
  the {O}rthographic {S}hape--from--{S}hading {P}roblem.
\newblock In Sebastiano Battiato and José Braz, editors, {\em Proceedings of
  the 9th International Conference on Computer Vision Theory and Applications
  (VISAPP)}, volume~3, pages 711--716. SciTePress, 2014.

\bibitem{TMDDb15}
S.~Tozza, R.~Mecca, M.~Duocastella, and A.~Del Bue.
\newblock {Direct differential Photometric-Stereo shape recovery of diffuse and
  specular surfaces}.
\newblock {\em International Journal of Mathematical Imaging and Vision
  (JMIV)}.
\newblock To appear.

\bibitem{VBLW09}
O.~Vogel, M.~Breu\ss, T.~Leichtweis, and J.~Weickert.
\newblock Fast {S}hape from {S}hading for phong-type surfaces.
\newblock In Xue-Cheng Tai, Knut M{\o}rken, Marius Lysaker, and Knut-Andreas
  Lie, editors, {\em Scale Space and Variational Methods in Computer Vision},
  volume 5567 of {\em Lecture Notes in Computer Science}, pages 733--744.
  Springer, 2009.

\bibitem{VBW08}
O.~Vogel, M.~Breu\ss, and J.~Weickert.
\newblock Perspective {S}hape from {S}hading with non-{L}ambertian reflectance.
\newblock In {\em Proc. DAGM Symposium on Pattern Recognition}, pages 517--526,
  2008.

\bibitem{VVBW09}
O.~Vogel, L.~Valgaerts, M.~Breu\ss, and J.~Weickert.
\newblock Making shape from shading work for real-world images.
\newblock In Joachim Denzler, Gunther Notni, and Herbert S\"u{\ss}e, editors,
  {\em Pattern Recognition}, volume 5748 of {\em Lecture Notes in Computer
  Science}, pages 191--200. Springer Berlin Heidelberg, 2009.

\bibitem{wikipedia_Agamemnon}
Wikipedia.
\newblock Mask of {A}gamemnon.
\newblock Website.
\newblock http://en.wikipedia.org/wiki/Mask\_of\_Agamemnon.

\bibitem{WON98}
Lawrence~B. Wolff, Shree~K. Nayar, and Michael Oren.
\newblock Improved diffuse reflection models for computer vision.
\newblock {\em International Journal of Computer Vision}, 30(1):55--71, 1998.

\bibitem{Wolff_JOSAA1994}
L.B. Wolff.
\newblock Diffuse-reflectance model for smooth dielectric surfaces.
\newblock {\em Journal of the Optical Society of America A}, 11(11):2956--2968,
  1994.

\bibitem{Woo80}
R.~J. Woodham.
\newblock Photometric method for determining surface orientation from multiple
  images.
\newblock {\em Optical Engineering}, 19(1):134--144, 1980.

\bibitem{WH99}
P.L. Worthington and E.R. Hancock.
\newblock New constraints on data-closeness and needle map consistency for
  {S}hape-from-{S}hading.
\newblock {\em IEEE Transactions on Pattern Analysis and Machine Intelligence},
  21(12):1250--1267, 1999.

\bibitem{WNJ10}
C.~Wu, S.~Narasimhan, and B.~Jaramaz.
\newblock A multi-image {S}hape-from-{S}hading framework for near-lighting
  perspective endoscopes.
\newblock {\em International Journal of Computer Vision}, 86:211--228, 2010.

\bibitem{Yoon_IJCV2010}
K.-J. Yoon, E.~Prados, and P.~Sturm.
\newblock Joint estimation of shape and reflectance using multiple images with
  known illumination conditions.
\newblock {\em International Journal of Computer Vision}, 86:192--210, 2010.

\bibitem{ZTCS99}
R.~Zhang, P.-S. Tsai, J.E. Cryer, and M.~Shah.
\newblock {S}hape from {S}hading: a survey.
\newblock {\em IEEE Transactions on Pattern Analysis and Machine Intelligence},
  21(8):690--706, 1999.

\bibitem{ZC91}
Qinfen Zheng and R.~Chellappa.
\newblock Estimation of illuminant direction, albedo, and shape from shading.
\newblock In {\em Computer Vision and Pattern Recognition, 1991. Proceedings
  CVPR '91., IEEE Computer Society Conference on}, pages 540--545, Jun 1991.

\end{thebibliography}

\end{document}